\documentclass[final]{siamltex}

\usepackage{geometry}
\usepackage{mathptmx}      % use Times fonts if available on your TeX system
\usepackage{amsfonts}
\usepackage{graphicx}
\usepackage{epstopdf}

\usepackage{amsmath}
\usepackage{textcomp}
\usepackage{amssymb}
\usepackage[numbers,sort&compress]{natbib}
\newcommand{\figpath}{./}

\newcommand{\ie}{{\it i.e.}~}
\newcommand{\eg}{{\it e.g.}~}

\newcommand{\cf}{{\it cf.}~}

\newcommand{\x}{\ensuremath{\mathbf{x}}}

\newcommand{\Real}{\ensuremath{\mathbb{R}}}
\newcommand{\RealN}[1]{\ensuremath{\mathbb{R}^{#1}}}

\newcommand{\Int}{\ensuremath{\mathbb{Z}}}

\newcommand{\CircleN}[1]{\ensuremath{\mathbb{S}^{#1}}}

\newcommand{\StableSet}[1]{\ensuremath{W^{s}(#1)}}

\newcommand{\Isyn}{\ensuremath{I_\mathrm{syn}}}
\newcommand{\gSyn}{\ensuremath{g_\mathrm{syn}}}
\newcommand{\vSyn}{\ensuremath{V_\mathrm{syn}}}
\newcommand{\aSyn}{\ensuremath{\alpha_\mathrm{syn}}}
\newcommand{\bSyn}{\ensuremath{\beta_\mathrm{syn}}}

\newcommand{\gSyni}{\ensuremath{g^{i}_\mathrm{syn}}}
\newcommand{\gSyne}{\ensuremath{g^{e}_\mathrm{syn}}}
\newcommand{\MP}{\ensuremath{\mathcal{M}_P}} % Periodic orbits
\newcommand{\MS}{\ensuremath{\mathcal{M}_S}} % Saddle points
\newcommand{\MQ}{\ensuremath{\mathcal{M}_Q}} % Quiescent fixed points

\title{Dissecting the Phase Response of a Model Bursting Neuron\thanks.This work was supported by NIH CRCNS grant 1R01NS050943, DOE grant DE-FG02-93ER25164, and NSF FIBR grant 0425878 subcontract SA4554-10295PG.}

\author{William Erik Sherwood\thanks{Center for BioDynamics, 
              Boston University, 
              111 Cummington Street ,
              Boston, MA 02215 ({\tt wesher@bu.edu}).}
        \and John Guckenheimer\thanks{Mathematics Department, 
              Cornell University, 
              565 Malott Hall, 
              Ithaca, NY 14853 ({\tt jmg16@cornell.edu}).}}

\begin{document}

\maketitle

\begin{abstract}
We investigate the phase response properties of the Hindmarsh-Rose model of neuronal bursting using burst phase response curves (BPRCs)  computed with an infinitesimal perturbation approximation and by direct simulation of synaptic input. The resulting BPRCs have a significantly more complicated structure than the usual Type I and Type II PRCs of spiking neuronal models, and they exhibit highly timing-sensitive changes in the number of spikes per burst that lead to large magnitude phase responses. We use fast-slow dissection and isochron calculations to analyze the phase response dynamics in both weak and strong perturbation regimes.
\end{abstract}

\begin{keywords} 
phase response, isochron, multiple time-scales, neuronal model, bursting
\end{keywords}

\begin{AMS}
92B25, 37N25, 37C27, 37C37, 37M05
\end{AMS}

\pagestyle{myheadings}
\thispagestyle{plain}
\markboth{W. E. SHERWOOD AND J. GUCKENHEIMER}{PHASE RESPONSE OF A BURSTING NEURAL MODEL}

\section{Introduction}
\label{sec:introduction}

The activity of synaptically interacting neurons is commonly modeled using networks of coupled oscillators. Among the most analytically tractable of such models are ``phase oscillator'' models that reduce the dynamics of each component neuron to a limit cycle and assume that the coupling between oscillators depends only upon the oscillators' positions along these limit cycles \cite{Cohen:1982, Kopell:1988, Ermentrout:1990,  Hoppensteadt:1997}. The states of the network are understood to be functions of the phase differences between neurons, and the interactions between neurons act to change these phase relationships but are presumed not to perturb any oscillator from its limit cycle. To understand the behavior of the network, one first maps how the phase of an individual oscillator changes in response to input, and then studies the behavior of these maps when coupled together according to the architecture of the network. 

This approach has been particularly useful for building and analyzing models of models of central pattern generators (CPGs), localized, autonomous neuronal networks that produce patterned rhythmic output underlying such behaviors as circulation, digestion, respiration, and locomotion. The corresponding mathematical theory has been developed largely in the context of simple oscillator models and weak coupling, yet the theory has held up well when applied to model networks of coupled spiking neurons \cite{Canavier:1997, Ermentrout:1998, Goel:2002, Ermentrout:2006}. But neurons exhibit a wider range of complex, nonlinear behaviors than just spiking, such as bursting, the periodic alternation between periods of spiking and silence. In many CPGs, such as the respiratory CPG of the mammalian preB\"{o}tzinger complex \cite{Butera:1999, Butera:1999a} and the digestive CPG of the crustacean stomatogastric ganglion~\cite{Harris-Warrick:1992}, bursting neurons play central roles and may be strongly coupled to other components of the network. What are the phase resetting properties of bursting neurons? How do they respond to strong perturbations?

In this paper, we investigate these questions for a particular model of neuronal bursting. We organize our study as follows: We first summarize the aspects of multiple time scale dynamics and geometric singular perturbation theory relevant for the mathematical analysis of neuronal bursting activity, and we introduce the Hindmarsh-Rose model, which we employ as a canonical model of bursting throughout the paper. We then present the terminology of phase response and introduce the concept of isochrons, and we discuss  infinitesimal and direct methods of calculating phase response curves. Next, we present phase response curves for the Hindmarsh-Rose model, calculated using both infinitesimal and direct methods, and we compare the features of these curves with typical phase response curves for spiking neural models. In the remainder of the paper, we first analyze of the shape of the infinitesimal phase response curves using fast-slow dissection and isochron portraits, and we then examine the phase response to strong perturbations. 
 
\section{Modeling neuronal bursting}
\label{sec:modeling_bursting}

Before describing neuronal dynamics in mathematical terms, we first define some neuroscience terminology which we will use throughout the paper: Electrically excitable cells such as neurons generically display three modes of behavior: Quiescence, tonic spiking, and bursting. In {\it quiescence}, a cell maintains a stable resting membrane potential, but if perturbed (via temporary injection of current, for example), the membrane potential makes a large excursion, or action potential, rising rapidly far above the resting potential, then quickly falling somewhat below the resting potential, before reestablishing its equilibrium value. {\it Tonic spiking} is marked by continual firing of action potentials, normally at a fixed rate, and may occur endogenously or in the presence of a steady stimulus, such as an external current. {\it Bursting} is characterized by relatively slow, periodic alternation between periods of spiking (active state) and periods of electrical inactivity (quiescent state). A relatively elevated membrane potential, one above the normal equilibrium value for a cell, is said to be {\it depolarized}, and the act of raising membrane voltage is called {\it depolarization}. Recovery to a normal voltage equilibrium following hyperpolarization is called {\it repolarization}. A relatively depressed membrane potential, below the cell's normal equilibrium value, is said to be {\it hyperpolarized}, and lowering membrane voltage is called {\it hyperpolarization}. Inputs to a neuron which tend to raise the membrane voltage and depolarize the cell are called {\it excitatory}, while those which tend to lower the membrane voltage and hyperpolarize the cell are termed {\it inhibitory}.

\subsection{Multiple time-scales and bursting}
\label{sec:multiple_time_scales}
Bursting is a phenomenon involving dynamics on multiple time-scales; as such, it is amenable to analysis using the tools of geometric singular perturbation theory \cite{Fenichel:1979, Jones:1995}. The canonical form of multiple time-scale or {\it fast-slow} systems is

%%%%%%%%%%%%%%%%%%
%	Fast system equations
%%%%%%%%%%%%%%%%%%  
\begin{eqnarray}
\label{eqn:fastx} \dot{x} &=& f(x,y,\lambda)\\
\label{eqn:fasty} \dot{y} &=& \epsilon g(x,y,\lambda)
\end{eqnarray}

\noindent where $x \in \mathbb{R}^{m}, y \in \mathbb{R}^{n}$, $\lambda$ denotes the of parameters, and $\epsilon \ll 1$ is a small positive parameter explicitly denoting the separation of time-scales between the {\it fast variables} $x$ and the {\it slow variables} $y$. Equations (\ref{eqn:fastx})--(\ref{eqn:fasty}) are called the {\it fast system}. Rescaling time, we may rewrite (\ref{eqn:fastx})--(\ref{eqn:fasty}) to obtain  the {\it slow system}:

%%%%%%%%%%%%%%%%%%
%	Slow system equations
%%%%%%%%%%%%%%%%%%   
\begin{eqnarray}
\label{eqn:slowx} \epsilon x' &=& f(x,y,\lambda)\\
\label{eqn:slowy} y' &=&  g(x,y,\lambda)
\end{eqnarray}

\noindent Here $'$ denotes differentiation with respect to rescaled time. For $\epsilon \neq 0$ the two systems are equivalent. Letting $\epsilon \rightarrow 0$ in the fast system, we obtain the family of {\it fast subsystems}: 

%%%%%%%%%%%%%%%%%%
%	Fast subsystem equations
%%%%%%%%%%%%%%%%%%  
\begin{eqnarray}
\label{eqn:fastsubx} \dot{x} &=& f(x,y,\lambda)\\
\label{eqn:fastsuby}\dot{y} &=& 0
\end{eqnarray}

\noindent For each $y$, the fast subsystem is an $m$-dimensional differential equation in $x$, with $y$ acting as additional parameters. The singular limit $\epsilon = 0$ of the slow system produces a differential-algebraic system:

%%%%%%%%%%%%%%%%%%
%	Slow subsystem equations
%%%%%%%%%%%%%%%%%%  
\begin{eqnarray}
\label{eqn:slowsubx} 0 &=& f(x,y,\lambda)\\
\label{eqn:slowsuby} y' &=&  g(x,y,\lambda)
\end{eqnarray}

\noindent Equation (\ref{eqn:slowsubx}) defines an $n$-dimensional manifold, called the {\it critical manifold}, which is the union of equilibria of the fast subsystems. If the equilibria (\ref{eqn:fastsubx}) are hyperbolic, then the critical manifold is normally hyperbolic. Equation (\ref{eqn:slowsuby}) implicitly defines a vector field, the {\it slow flow}, on the critical manifold at regular points of its projection onto the $y$ coordinate hyperplane. Fenichel theory characterizes the relationship between the full system and slow flow \cite{Fenichel:1972, Fenichel:1975, Fenichel:1977}, establishing, in particular, the existence of a {\it slow manifold}, an invariant manifold distance $O(\epsilon)$ away from the critical manifold on compact regions of regular points. The slow flow on the critical manifold provides an $O(\epsilon)$ approximation of the trajectories of the full system. Thus trajectories of the full system closely follow stable structures (fixed points or limit cycles) in the fast subsystems, with transitions occurring near bifurcations of these structures. Later in the paper, we will invoke the Fenichel theorems in order to justify our approximation of the isochrons of a bursting neuronal model.

Rinzel and Lee first applied the tools of geometric singular perturbation theory to analyze bursting \cite{Rinzel:1987,Rinzel:1987a}; it has since become a standard approach. The basic idea of {\it fast-slow dissection} is to identify the phase space variable(s) associated with the slowest membrane current and treat it as a quasi-static parameter, thereby obtaining a family of fast subsystems as the parameter varies. Over the course of a bursting cycle in the full system, the slow variable traces out a periodic trajectory, giving a range over which it varies when considered as a quasi-static parameter. The active spiking state of the full system corresponds to the burst trajectory closely tracking a family of stable limit cycles in the set of fast subsystems, while the quiescent state corresponds to the burst trajectory staying close to a curve of (hyperpolarized) stable fixed points. As the slow variable varies, bifurcations in the fast subsystem create, destroy, and change stability of the fixed points and limit cycles being tracked by the full system trajectory, prompting rapid switching between full system states corresponding to qualitatively different behaviors (spiking vs. quiescence). In order to produce bursting, the evolution of the slow variable must be cyclic, such that the sequence of fast subsystem bifurcations repeats periodically. The pairs of fast subsystem bifurcations marking the initiation and termination of the active spiking state (destruction of the stable fixed points and limit cycles, respectively) determine important burst characteristics, such as spike frequency at burst onset and spike frequency adaptation, and have been used to construct taxonomies for bursting models \cite{Rinzel:1987a, Izhikevich:2000}.

%%%%%%%%%%%%%%%%%%
%
%	HR INTRO
%
%%%%%%%%%%%%%%%%%%

\subsection{Hindmarsh-Rose model}
\label{sec:HR_model}
In this paper we focus on a particular model of neuronal bursting, the the Hindmarsh-Rose (HR) equations \cite{Hindmarsh:1984}, which augment with a third slow variable the planar Fitzhugh-Nagumo model  of the action potential \cite{Fitzhugh:1961}, itself a reduction of the Hodgkin-Huxley model of the action potential in the squid giant axon. The HR equations were originally introduced to model the intrinsic alternation between spiking and quiescence reported in isolated neurons of the pond snail {\it Lymnea} and the R15 neuron of the mollusc {\it Aplysia}. Though strictly speaking not a biophysical neuronal model, the HR model captures the essential features of both spiking and bursting, and it is commonly used as a `minimal' or `prototypical' model of bursting. The HR equations are: 

%%%%%%%%%%%%%%%%%%
%	Full HR equations
%%%%%%%%%%%%%%%%%%  
\begin{eqnarray}
\dot{V} &=& n - aV^3 + bV^2 - h + I  \label{eqn:HR_full_V}\\ 
\dot{n} &=& c - dV^2- n \label{eqn:HR_full_n}\\
\dot{h} &=& r(\sigma(V - V_0) - h) \label{eqn:HR_full_h}
\end{eqnarray}
\noindent We refer to (\ref{eqn:HR_full_V})--(\ref{eqn:HR_full_h}) as the {\it full (HR) system}, and we exclusively employ the parameter set $a = 1$, $b=3$, $c=1$, $d=5$, $r=0.001$, $\sigma=4$, $V_0 = -1.6$, $I = 2$, which produces a bursting periodic orbit with 9 spikes per burst. 

\begin{figure}[ht]
 \begin{minipage}{3in}
     \begin{center}
    \includegraphics[width= 3in]{\figpath/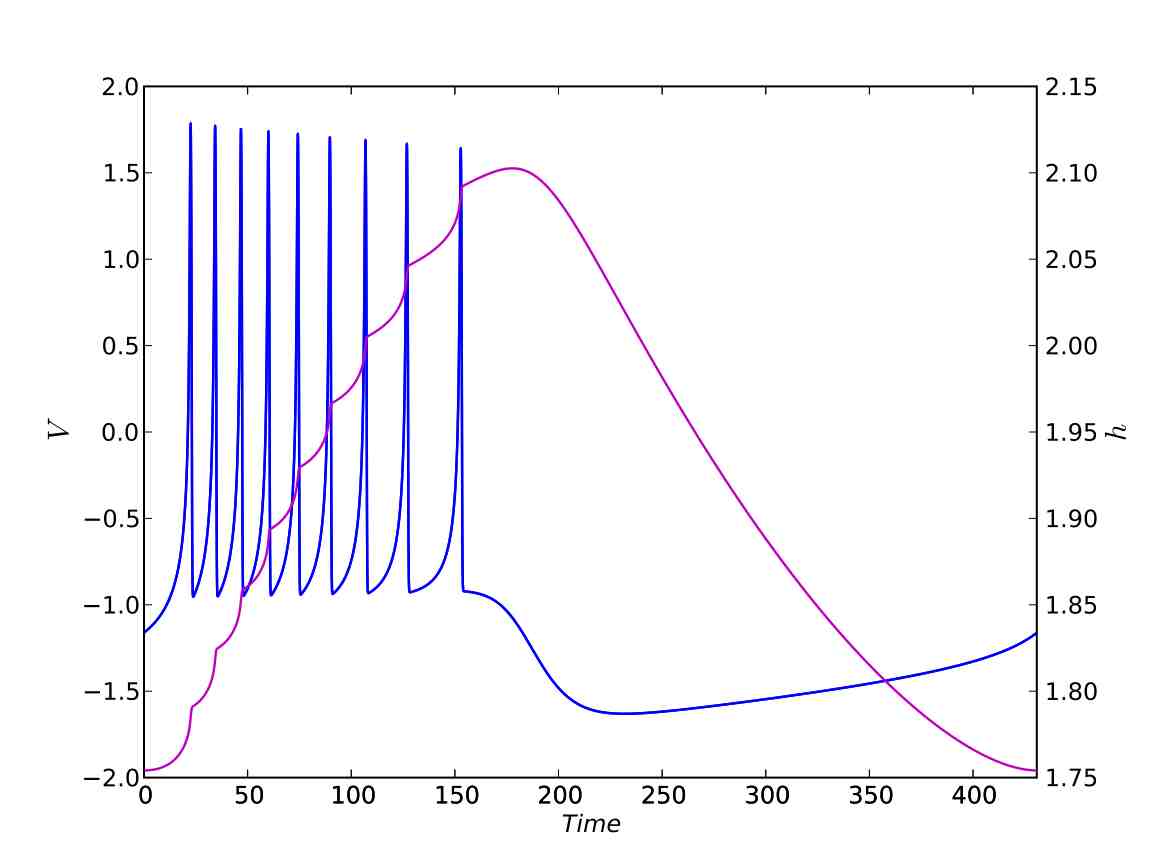}
       (a)
    \end{center}
  \end{minipage}%
  \begin{minipage}{3in}
      \begin{center}
    \includegraphics[width= 3in]{\figpath/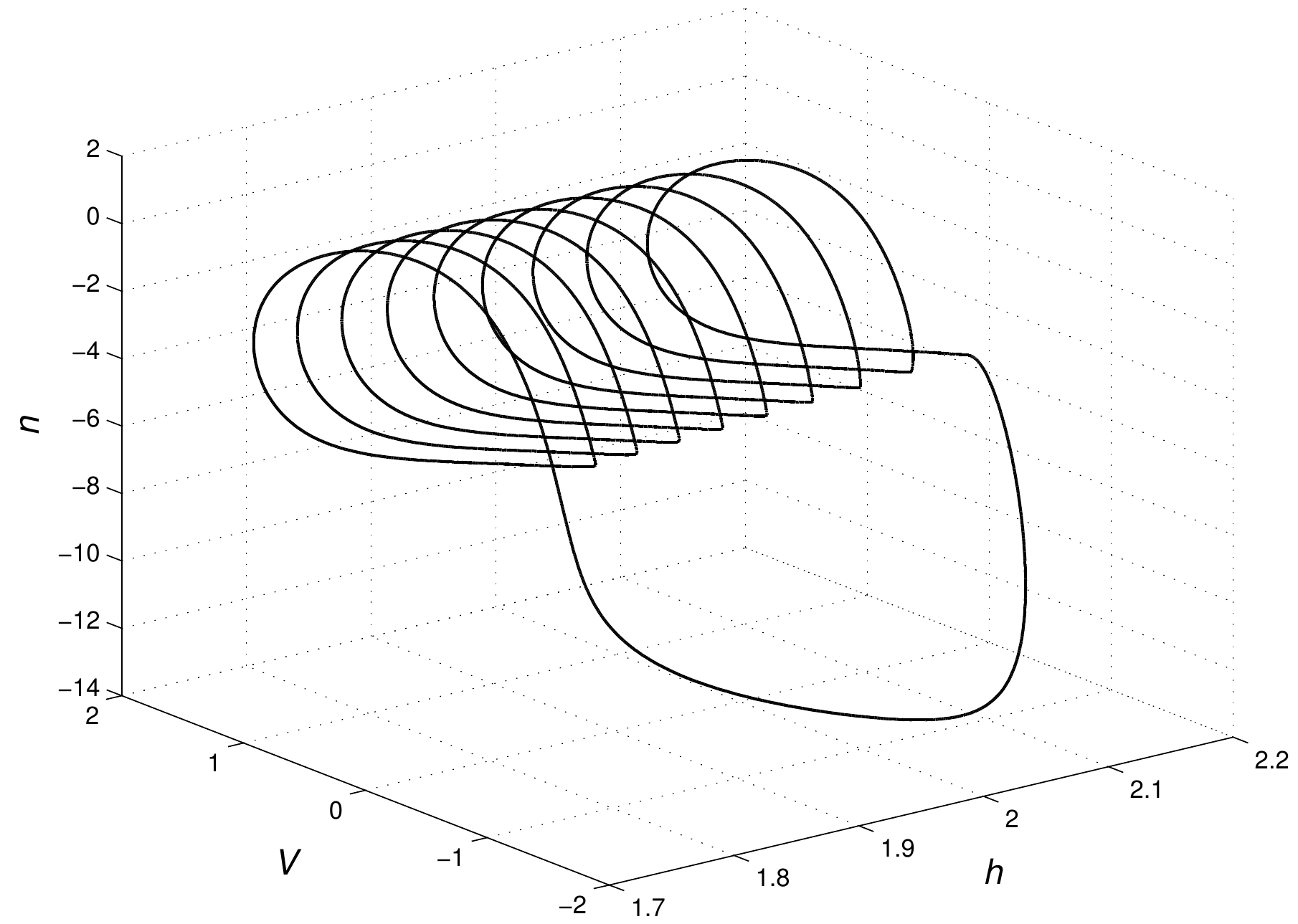}
        (b)
    \end{center}
  \end{minipage}
\caption[HR Burst Trajectory]{Full system burst trajectories for HR model. (a) $V$ and $h$, the slow variable, vs. time. (b) Three dimensional ``reference burst trajectory'' }
\label{fig:HR_Burst_Traj}
\end{figure}

In the HR system, the variable $V$ represents the neuron's membrane voltage and $I$ may be interpreted as an applied current. The parameter $r$ sets the separation of time-scales between $V$ and  the slow `recovery' variable $h$, whose dynamics determines the alternation between spiking and quiescence.  Figure \ref{fig:HR_Burst_Traj} (a) shows the temporal evolution of $V$ and $h$ over the course of one burst in the full HR system; Figure \ref{fig:HR_Burst_Traj} (b) depicts a burst orbit in all three phase space variables. Spiking begins near the point at which $h$ reaches it minimum value and ends just prior to the point at which $h$ is maximum.

%%%%%%%%%%%%%%%%%%
%	HR burst trajectory 2D,3D
%%%%%%%%%%%%%%%%%%  

\noindent Setting equation (\ref{eqn:HR_full_h}) equal to zero and treating $h$ as a quasi-static parameter, we obtain the {\it fast (HR) subsystem} described by (\ref{eqn:HR_full_V}--\ref{eqn:HR_full_n}).\footnote{Though they did not employ geometric singular perturbation theory explicitly, Hindmarsh and Rose did invoke the wide separation of time-scales to justify treating $h$ as a quasi-static parameter in order to perform a phase plane analysis of bursting in the full HR system.} The bifurcation structure of the HR fast subsystem determines its bursting behavior; as described below, the HR model belongs to the category of saddle-node/homoclinic bursters \cite{Izhikevich:2000}, also (traditionally) termed Type I or `square-wave' bursters \cite{Rinzel:1987a}. 

\begin{figure}[ht]
 \begin{minipage}{3in}
     \begin{center}
    \includegraphics[width= 3in]{\figpath/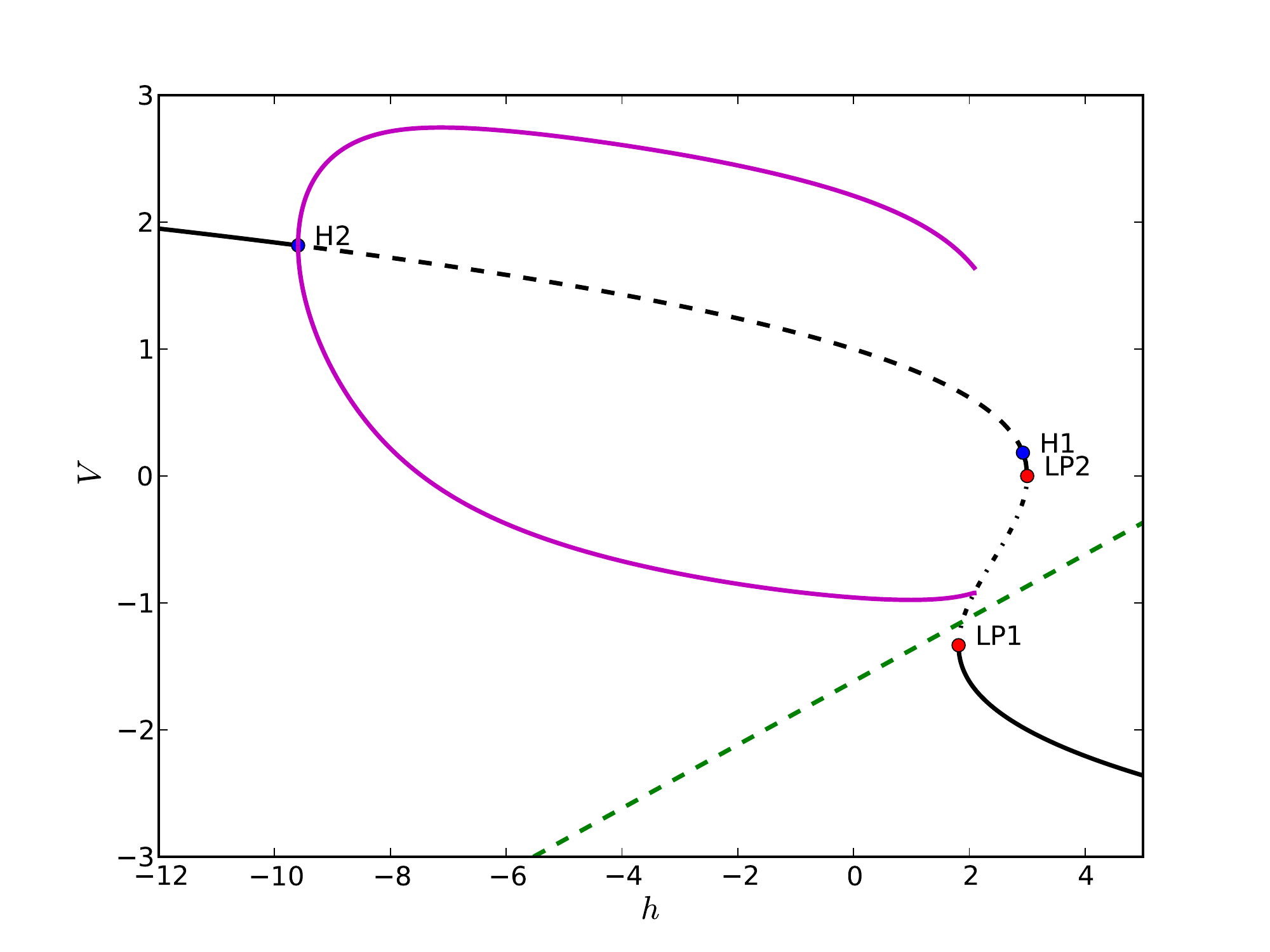}
       (a)
    \end{center}
  \end{minipage}%
  \begin{minipage}{3in}
      \begin{center}
    \includegraphics[width= 3in]{\figpath/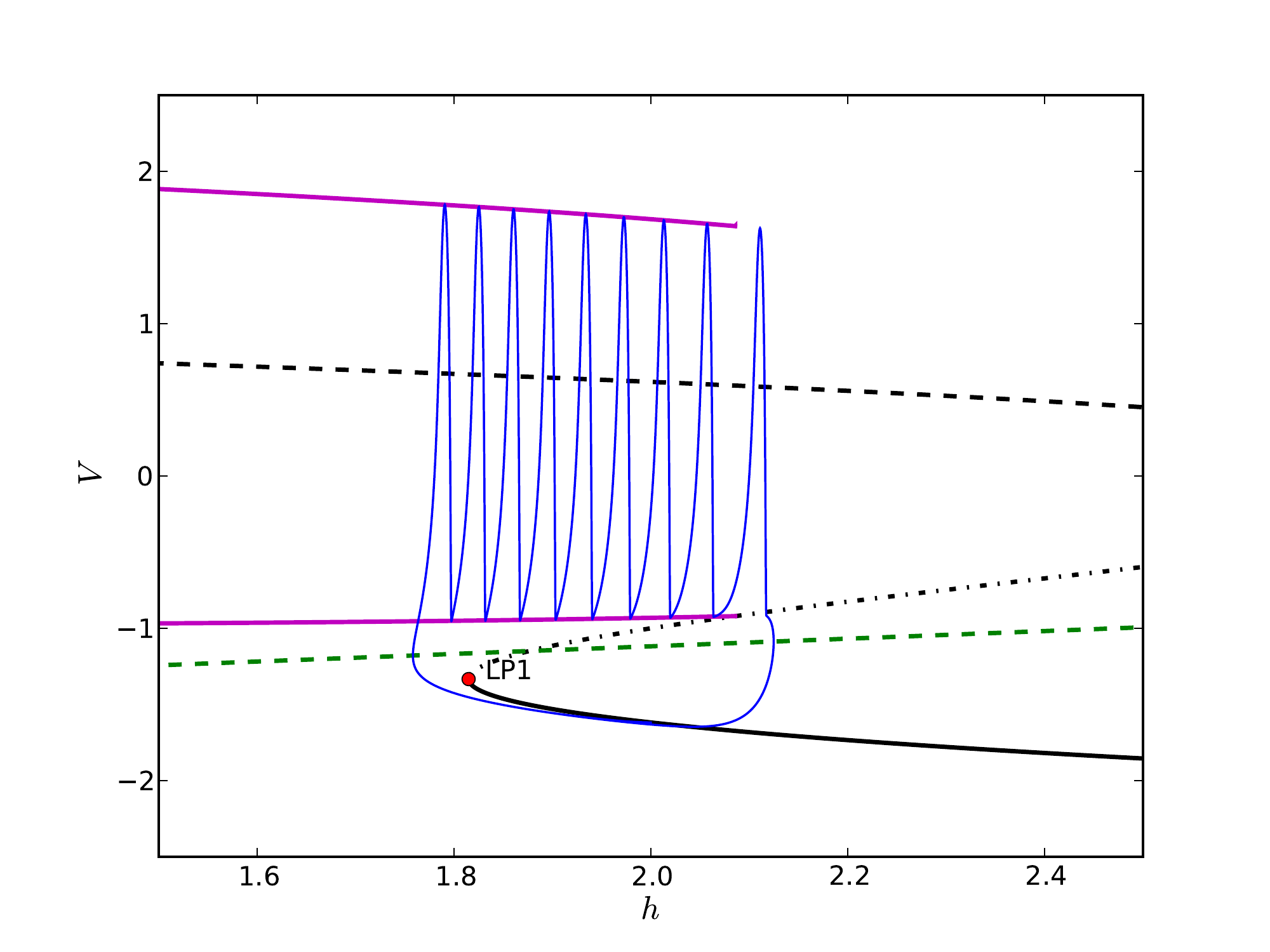}
        (b)
    \end{center}
  \end{minipage}
\caption[HR Fast Subsystem Bifurcations]{Bifurcation diagram for the HR fast subsystem. Solid black lines indicate stable fixed points, dashed black lines show unstable fixed points. Solid magenta lines show the maximum and minimum $V$ values of stable limit cycles. The dashed green line indicates the $h$-nullcline of the full HR system projected onto the $(h,V)$-plane. Saddle-node bifurcation points: LP1, LP2. Hopf bifurcation points: H1, H2.  (a) Full bifurcation diagram. (b) Close-up of $h$ range for bursts in the full subsystem. An $(h, V)$-plane projection of a full system trajectory is superimposed in blue.}
\label{fig:HR_fast_bifns}
\end{figure}

Figure (\ref{fig:HR_fast_bifns}) shows the bifurcation diagram (computed with AUTO \cite{AUTO} via the PyCont module of PyDSTool, a Python-based toolkit for simulation and analysis of dynamical systems \cite{PyDSTool:2007}) for the HR fast subsystem, using $h$ as the continuation parameter. For large negative $h$ values, there exists a lone stable fixed point, and at $h \approx  -9.5931404587$ there is a supercritical Hopf bifurcation (denoted H2 in Figure \ref{fig:HR_fast_bifns} (a)) where the fixed point loses stability and a stable limit cycle is created. At large positive $h$ values, there again exists a lone fixed point, though at lower (negative) $V$ values than the fixed point that exists for large, negative $h$. At $h = 3.0$ a pair of fixed points (one stable, the other unstable) is born at the saddle-node bifurcation denoted LP2 in Figure \ref{fig:HR_fast_bifns}. These fixed points have higher $V$ values than the other stable fixed point; the lower, unstable member of the pair is a saddle. As $h$ decreases to 2.92647388572, the supercritical Hopf bifurcation point denoted H1, the upper stable fixed point loses stability and a small stable limit cycle emerges. This small periodic orbit is quickly destroyed in a homoclinic bifurcation with the lower saddle point. At $h \approx 1.81481481481$, denoted LP1,  the lower stable fixed point merges with the saddle, and the two lower fixed points disappear, leaving only the upper, unstable fixed point and the stable limit cycle (from H2) around it. The stable limit cycle from H2 is destroyed in a homoclinic bifurcation, merging with the saddle point at $h \approx 2.08560088198$. 

%%%%%%%%%%%%%%%%%%
%	HR fast subsystem bifurcations
%%%%%%%%%%%%%%%%%%  

\noindent The upper unstable fixed point, saddle point, lower stable fixed point, and stable periodic orbit from H2 are the only fast subsystem structures which play a role in full system bursting. As shown in Figure \ref{fig:HR_Burst_Traj} (b), $h$ varies between 1.75415439813 and 2.10256601768 in the full HR system over the course of a burst cycle. Figure \ref{fig:HR_fast_bifns} (b) shows a close-up of this region of the bifurcation diagram along with a burst trajectory in the full system projected onto the $(h, V)$-plane. The line of stable fixed points at lower $V$, denoted \MQ, corresponds to the resting membrane voltage during quiescence; also shown is the $h$-nullcline of the full system: In the full system, $h$ decreases as it tracks the lower fixed points, as they lie below the $h$-nullcline. The disappearance of these fixed points at the saddle-node bifurcation LP1 corresponds to the start of the next active segment and the onset of spiking, when the full system ceases tracking the line of low voltage stable fixed points. Since the full system now lies above the $h$-nullcline, $h$ increases as the full system follows the family of stable periodic orbits. The homoclinic bifurcation near the maximum value of $h$ marks the end of the active segment, when the full system trajectory leaves the vicinity of the family of (now extinguished) periodic orbits and follows the line of low voltage stable fixed points. 

%%%%%%%%%%%%%%%%%%%%%%%
%	HR 3D reference burst w/fast subsystem
%%%%%%%%%%%%%%%%%%%%%%% 
\begin{figure}[!ht]
 %\begin{minipage}{3.25in}
     \begin{center}
	\includegraphics[width= 3.5in]{\figpath/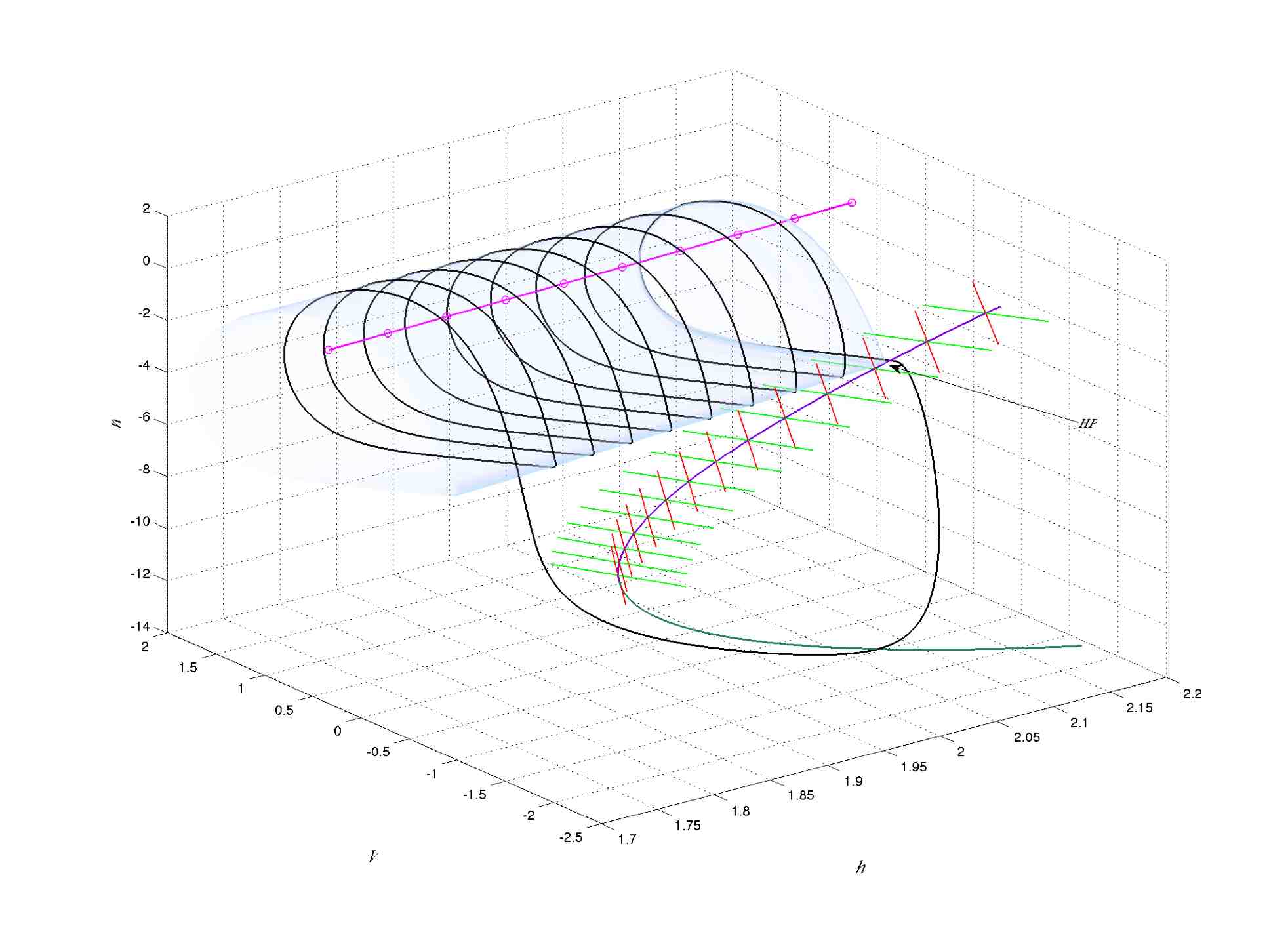}
	%\
	    \end{center}
  %\end{minipage}
\caption[Fast Subsystem 3D Portrait]{Fast subsystem structures and burst trajectory of the HR model. The curve of hyperpolarized, stable fixed points of the fast subsystem, \MQ, is drawn in dark green; the line of saddles, \MS, is drawn in purple, along with the stable (green) and unstable (red) eigendirections. The pale blue tube is the manifold of stable limit cycles in the family of fast subsystems, \MP; the full system burst trajectory hugs its outer surface. The magenta line represents the line of unstable fixed points interior to the fast subsystem limit cycles. The homoclinic point of the fast subsystem, HP,  is indicated by the arrow. }
\label{fig:HR_ref_3D}
\end{figure}

This sequence is recapitulated in Figure \ref{fig:HR_ref_3D}, which depicts in three dimensions the configuration of fast subsystem structures along a burst trajectory of the full system. During the active spiking segment of the burst, the full system trajectory hugs the outer surface of the manifold of  fast subsystem stable limit cycles, which we denote $\MP$. At the homoclinic point in the fast subsystem, the stable limit cycle merges with the stable and unstable manifolds of the saddle, while in the full system, the burst trajectory crosses the line of saddles, denoted \MS, along its stable manifold and jumps down to the curve of stable, hyperpolarized fixed points \MQ, thus terminating spiking and beginning the quiescent segment of the burst cycle.

\section{Phase response}
\label{sec:phase_response}
We consider the definition and computation of phase response in more detail before turning to the particulars of the Hindmarsh-Rose case. Our discussion draws substantially upon the exposition in \cite{Sherwood:2008, Brown:2004, Govaerts:2006, Guillamon:2009}. 

\subsection{Phase and isochrons}
\label{sec:phase_and_isochrons}
Let the evolution of an autonomous dynamical system $F$ be governed by the system of ordinary differential equations

%%%%%%%%%%%%%%%%%%
%	Dynamical system
%%%%%%%%%%%%%%%%%%  
\begin{equation}
\dot{x} = f(x),\ x \in U \subseteq \RealN{n},\ n \geq 2
\end{equation}

\noindent with associated (smooth) flow $\Phi(x,t),\ \Phi:\RealN{n} \times \Real \rightarrow \Real$, and let $\Gamma$ be a stable hyperbolic limit cycle for $F$ with period $T$. As $\Gamma$ is a one dimensional manifold, the natural phase space of $F$ restricted to $\Gamma$ is diffeomorphic to the the circle $\CircleN{1} = \Real/\Int$, and thus we may parameterize $\Gamma$ in terms of a single scalar variable, or {\it phase}, $\theta \in [0,1), \theta := t/T$, such that $\Gamma$ has period 1 under this parameterization, \ie $\Gamma(\theta) = \Gamma(\theta + 1)$.

%%%%%%%%%%%%%%%%%%
%	Phase parameterization
%%%%%%%%%%%%%%%%%%  
\begin{eqnarray}
\Gamma: \CircleN{1} &\rightarrow& \RealN{n}, \\
\theta &\mapsto& \Gamma(\theta)\nonumber
\end{eqnarray}

Each point $x \in \Gamma$ is associated with a unique {\it phase} $\theta$; a natural labeling of phases on $\Gamma$ is obtained by fixing a reference point $x_{0} \in \Gamma$ and measuring each point's temporal distance to $x_{0}$ under $\Phi$: or $x \in \Gamma$, if $\Phi(x_{0},t) = x$ we write $\theta(x) = t/T\mod 1$. This idea of phase represented by $\theta$ is strictly meaningful only on $\Gamma$, but we can extend the concept to include points in the basin of attraction of $\Gamma$, denoted $B_{\Gamma}$.

%%%%%%%%%%%%%%%%%%
%	Isochron definition
%%%%%%%%%%%%%%%%%%  
\begin{definition}\label{def:isochron} For $x \in \Gamma$ with phase $\theta = \theta(x)$, we say $y \in B_{\Gamma}$ has {\it asymptotic phase} $\vartheta(y) = \theta$ if 
\begin{equation}
\lim_{t \rightarrow \infty} \|\Phi(x,t) - \Phi(y,t)\| = 0.
\end{equation} 
The set of all points $y \in B_{\Gamma}$ such that $\vartheta(y) = \theta(x) = \theta$ is called the {\it isochron} of phase $\theta$ or the {\it isochron} at point $x$, denoted $\varpi(x)$.
\end{definition} 
%When associating the scalar value $\theta(\x)$ with $\y$, we call it the {\it asymptotic}  {\it phase} of $\y$ and denote it by $\vartheta(\y)$. For each $\x \in \Gamma$, we define the {\it isochron} $\varpi(\x)$ of $\x$ as the set of all points in the basin of attraction of $\Gamma$ having the same asymptotic phase as the actual phase of $\x$, {\it i.e.} $\varpi(\x) = \{ \y \in B: \vartheta(\y) = \theta(\x)\}$. For $\x \in \Gamma, \vartheta(x) = \theta(x)$.

\

\noindent This formulation of the notions of asymptotic phase and isochrons is due to Winfree \cite{Winfree:1974, Winfree:2001}. Equivalently, isochrons may be considered as cross-sections of $\Gamma$ (manifolds intersecting $\Gamma$ transversely at a single point) having a first return time equal to $T$. If $\Gamma$ is hyperbolic, their existence follows from a corollary of the Invariant Manifold Theorem stating that for each $\x \in \Gamma, \StableSet{x}$ is a cross-section of $\Gamma$, diffeomorphic to $\RealN{n-1}$. Each isochron is invariant under $\Phi$. Furthermore, $\cup_{x \in \Gamma} \StableSet{x}$ is an open neighborhood of $\Gamma$ and its stable manifold \cite{Guckenheimer:1975} The set of isochrons foliates $B_{\Gamma}$, and we may thus speak sensibly about phase for any point in $B_{\Gamma}$. %(Definitions for asymptotic phase and isochrons for unstable hyperbolic limit cycles may be obtained immediately by reversing the direction of time in definition (\ref{def:isochron}). For planar systems, asymptotic phase and isochrons can be defined for nonhyperbolic limit cycles as well; for the purposes of this paper, we consider only stable hyperbolic limit cycles.)

%We can also rewrite the dynamical system to emphasize phase evolution, by choosing $\theta(\cdot)$ appropriately, so that $F$ restricted to $U$ has the simple form $\frac{d\theta(\x)}{dt} = \omega$, where $\omega = 1/T$. This is accomplished using the chain rule: $\frac{d\theta(\x)}{dt} = \omega = \frac{d\theta}{d\x}(\x)\cdot f(\x)$. There is a unique solution, up to translation by a constant,  for $\theta(\cdot)$ on $U$.

\subsection{Phase response curves}
\label{sec:phase_response_curves}
Suppose we temporarily (smoothly) perturb $F$ (denote the perturbed dynamical system by $\hat{F}$) at point $x \in \Gamma$ having phase $\theta(x_0)$, so that at the end of the perturbation $x_0$ is mapped to point $\hat{F}(x_0) = y_0 \in B_{\Gamma}$. If the perturbation displaces points only along the limit cycle $\Gamma$, then the effect of $\hat{F}$ applied at time $t_0$ (equivalently, at phase $\theta(x_0) = t_0/T \mod 1$) at point $x_0 \in \Gamma$ is to increase or decrease (for one cycle) the period of $\Gamma$, or equivalently, to delay or advance the next time at which the perturbed trajectory returns to $x_0$. Then the {\it phase shift} (also {\it phase difference}, {\it phase resetting} or {\it phase response}) at point $x_0$ is given by

%%%%%%%%%%%%%%%%%%
%	Direct phase diff equations
%%%%%%%%%%%%%%%%%%  
\begin{equation}
\label{eqn:direct_phase_diff}
\Delta\theta(x_0) = (\hat{T} - T)/T = (\hat{t}_0 - t_0)/T
\end{equation}
\noindent where $\hat{T}$ is the subsequent period and $\hat{t}_0$ is the subsequent return time to $x_0$, respectively. Note that positive $\Delta\theta(x_0)$ indicates phase delay and negative $\Delta\theta(x_0)$ indicates phase advancement.

Typical perturbations do not displace points solely along the limit cycle, but rather points on the limit cycle are mapped into the limit cycle's basin of attraction. More generally, then, the asymptotic phase of $\hat{F}(x_0)$ is $\vartheta(y_0)$ following the perturbation. Systematically perturbing every point of $\Gamma$ to obtain a locus of points $\hat{\Gamma} = \cup_{x \in \Gamma} \hat{F}(x)$, we obtain a mapping of phases 

%%%%%%%%%%%%%%%%%%
%	Phase mapping
%%%%%%%%%%%%%%%%%%  
\begin{eqnarray}
\varphi: [0,1) &\rightarrow& [0,1)\\ 
\varphi(\theta(x))  &=& \vartheta(\hat{F}(x)) 
\end{eqnarray}

\noindent We call $\varphi$ the {\it phase response curve (PRC)} of $\Gamma$ for $\hat{F}$ \cite{Winfree:2001}. The shape of the PRC depends on the characteristics of $F$, $\Gamma$, and $\hat{F}$. 

For a weak perturbation of brief duration occurring at time $t^{*}$, we may write the perturbed system $\hat{F}$ as 

%%%%%%%%%%%%%%%%%%
%	Dirac perturbation
%%%%%%%%%%%%%%%%%%  
\begin{equation}
\dot{x} = f(x) + I\delta(t-t^{*})
\end{equation}

\noindent where $I = (I_1, \ldots, I_n) \in \RealN{n}$ is the impulse applied by the perturbation and $\delta(t)$ is the Dirac delta function (\cf \cite{Guillamon:2009, Govaerts:2006}). For $\| I \| \ll 1$, the phase response imparted by the impulse is usually approximated by the {\it infinitesimal PRC}: in the limit of infinitesimal impulse $\| I \| \rightarrow 0$, the phase response for $x \in \Gamma$ is given by

%%%%%%%%%%%%%%%%%%
%	Infinitesimal PRC equations
%%%%%%%%%%%%%%%%%%  
\begin{equation}
\Delta\theta(x) =  \langle I , \nabla\vartheta(x)\rangle = \langle (I_1, \ldots, I_n), (\frac{\partial\vartheta}{\partial x_1}(x), \ldots, \frac{\partial\vartheta}{\partial x_n}(x))\rangle
\end{equation}

\noindent where $\langle \cdot, \cdot \rangle$ denotes the inner product. Most commonly in neuronal modeling, the impulse is taken to be purely a perturbation in the voltage variable, so that $I$ has only one nonzero coordinate, \eg if $x_1$ corresponds to the voltage variable $V$, then $I = (I_1, 0, \ldots, 0)$ so that 

%%%%%%%%%%%%%%%%%%
%	Inner product PRC
%%%%%%%%%%%%%%%%%%  
\begin{equation}
\Delta\theta(x) = I_1\frac{\partial\vartheta}{\partial V}(x).
\end{equation}

There are two classifications of PRC shapes in widespread use. The first is phenomenological \cite{Hansel:1995, Ermentrout:1996}: PRCs that are wholly positive or wholly negative are called Type I; regardless of timing, perturbations can have only one effect on the oscillator's phase (delay or advance, depending on sign). Type II PRCs have both positive and negative portions, so that the oscillator's phase may advance or delay depending on the timing of the perturbation. This classification of PRCs concords with behavioral categories for models of excitable membranes, which are in turn associated with bifurcation structures in the models \cite{Rinzel:1998, Ermentrout:1996}. Specifically, Type I PRCs are associated with excitable membranes that show arbitrarily low frequency oscillations at the onset of tonic spiking (saddle-node on invariant circle bifurcation). Type II PRCs are associated with excitable membranes for which the onset of tonic spiking occurs only at a fixed minimum (non-zero) frequency (\eg subcritical Hopf bifurcation).  

The second classification scheme is topological \cite{Winfree:2001}: Recall $\varphi: \CircleN{1} \rightarrow \CircleN{1}$ and consider the graph $G$ of $\varphi$ on the $(\theta, \hat{\theta})$-torus $\CircleN{2}$. If $G$ has degree 0 then the PRC is also said to be Type 0. Otherwise, $G$ must have degree greater than zero; if the degree is 1, the PRC is called Type 1. Oscillators with Type 0 PRCs are occasionally said to show `strong' or `even' phase response; those with Type 1 PRCs are sometimes referred to as having `odd' or `weak' phase response. Oscillators with Type 1 PRCs may exhibit large changes in their responses even as the phase of perturbation varies only slightly. This feature of the topological classification captures the large excursion and phase-sensitivity characteristics of biological oscillators subjected to strong perturbations. Both classification schemes assume that the oscillator remains in the basin of attraction of the original limit cycle after perturbation; otherwise, the PRC is undefined.

In the neuronal context, $\hat{F}$ typically represents the action of a single spike from a pre-synaptic neuron on the voltage and gating of ion channels in the membrane of the post-synaptic neuron. Action potentials characteristically have short duration compared to the post-synaptic periodic orbit. They comprise a rapid depolarizing upswing in voltage followed almost immediately by a hyperpolarizing downswing and then (usually slower) repolarization back towards equilibrium voltage. For small magnitude synaptic coupling, $\hat{F}$ is often simulated as an instantaneous increment or decrement of the voltage variable of the post-synaptic neuron, ignoring any change in the state of variables representing ion channel gating. The range of perturbation strengths for which this approximation is valid depends strongly on the isochron geometry of the system in the basin of attraction of the perturbed limit cycle, and the validity of the approximation is frequently left unexamined. 

Alternatively, additional terms describing the release of neurotransmitter by the synapse, neurotransmitter effects on membrane components, etc. may be incorporated into the post-synaptic neuronal model, in order to more fully simulate the effect of the synaptically mediated pre-synaptic voltage change on the post-synaptic neuron. Calculation of the phase response to such a perturbation requires integration of the unperturbed system $F$ up to the time of perturbation (incoming spike), integration of the perturbed system $\hat{F}$ for the duration of the perturbation, and then further integration of the original system. The phase response is then measured according to formula (\ref{eqn:direct_phase_diff}), and the phase response curve is constructed by applying the perturbation on a sufficiently dense, representative subset of times in $[0,T)$. We call phase response curves calculated by such methods {\it direct phase response curves}.\footnote{Most mathematical or computational studies of PRCs assume that perturbed trajectories return to the periodic orbit within one cycle and use the timing of a phase marker (\eg the peak of the action potential for a spiking neuron, or the first crossing of a voltage threshold after quiescence for a bursting neurons) on the next cycle to measure the phase response (\cf \cite{Acker:2003, Bose:2004, Brown:2004, Demir:1997}, for example). This procedure is also the norm for experimental studies of neuronal phase response (\cf \cite{Pinsker:1977, Prinz:2003, Tateno:2007}).}

The infinitesimal PRC method assumes instantaneous relaxation back to the original limit cycle after perturbation; after direct perturbation, the perturbed trajectory may take a considerable length of time, perhaps several times the original period, to return to the original limit cycle (or close enough to be considered as having returned). In the latter case, we may define, for $n \geq 1$, the $n$-th order phase shift for a perturbation applied to point $x_0$ at time $t_0$:

%%%%%%%%%%%%%%%%%%
%	nth phase shift
%%%%%%%%%%%%%%%%%%  
\begin{equation}
\label{eqn:nth_phase_shift}
\Delta_n\theta(x_0) = (\hat{t}^n_0 - t_0)/nT
\end{equation}
\noindent where $\hat{t}^n_0$ is the $n$-th crossing time (after $t_0$) of a hyperplane normal to $\Gamma$ at $x_0$. The $n$-th order phase response curve is constructed from the $n$-th order phase shifts for each point on $\Gamma$. For $n=1$ (and sufficiently weak perturbations), this definition essentially matches the definition of $\varphi$ above. The shapes of the $n$-th order PRCs, in comparison to the shape of the first order PRC, provide a measure of the long term persistence of phase shifts. 
%After perturbation, the trajectory of the model neuron reapproaches its original limit cycle, eventually getting close enough to be considered as having returned to the unperturbed limit cycle. In every case, regardless of model, perturbation type, or perturbation strength, the perturbed trajectory returned to the original limit cycle within three burst cycles. 

For a bursting neuron, the phase response curve records the shift in the timing of the onset of the active segment of a neuron's next burst cycle due to perturbation at a particular phase of its current burst cycle. This approximates the change in burst timing caused by a single incoming spike from a pre-synaptic neuron, as distinct from the change in timing of the next spike within the active segment of the next burst cycle. In a tonically spiking neuronal model, $\Gamma$ represents the periodic firing of a single action potential, and thus the change in timing of the onset of the next spike is precisely the relevant quantity in measuring phase response. For a bursting neuron, however, it is possible for a perturbation to change the timing of the next spike (of the perturbed burst cycle) without changing the timing of the next full burst cycle, \eg by shortening the next interspike interval but lengthening the quiescent portion of the burst cycle by an equal amount. We refer to phase response curves (both infinitesimal and direct) which measure the change in timing of the next burst cycle as {\it burst phase response curves (BPRCs)} to avoid confusion with PRCs that measure the shift in spike timing for oscillators representing tonically spiking neurons. 

It is possible for the perturbations used in direct BPRC calculations to add or delete spikes to the perturbed burst or subsequent bursts. The {\it spike number response curve (SNRC)} tallies the number of spikes in the perturbed burst for each phase $\theta$ at which it is perturbed.  For $n \geq 1$, the  $n$-th order SNRC counts the number of spikes occurring during the $n$-th burst cycle after perturbation. The zeroth order SNRC measures the change in spike number for the burst being perturbed. (Since the SNRC calculation is made in the context of perturbations far from the limit cycle, the definition of `burst cycle' requires some explanation; see the following subsection.) Changes in SNRC value mean that the perturbation added or deleted spikes from burst, thus indicating large deviations from the original limit cycle due to the perturbation.

\subsection{BPRC computation}
\label{sec:BPRC_computation}
Here we detail the calculations of burst phase response curves specific to our example Hindmarsh-Rose model. These computational details concern three sets of calculations: finding the periodic orbit of the reference burst trajectory, calculating the infinitesimal PRC, and computing direct burst PRCs.

\subsubsection{Reference periodic orbit}
\label{sec:reference_orbit}
 The full HR system (\ref{eqn:HR_full_V}--\ref{eqn:HR_full_n}) possesses a single stable periodic orbit, corresponding to a complete burst oscillation, which we term the {\it reference (burst) orbit}. The reference orbit  was found using multiple shooting with automatic differentiation \cite{Guckenheimer:2000}, as implemented in the ADMC++ automatic differentiation package for MATLAB \cite{Phipps:2003}, which produced initial conditions close to the true periodic orbit, with error less than $10^{-12}$. All other BPRC calculations were performed with PyDSTool and with AUTO via PyDSTool's PyCont module. Integration of the full HR system were performed using either a variable time-step fifth-order Runge-Kutta solver (Dormand-Prince-853, \cite{Hairer:1987}) with eighth-order dense output for the non-stiff systems or a fifth-order variable-time step implicit solver with eighth-order dense output (Radau5, \cite{Hairer:1991}) for stiff systems. The relative and absolute error tolerances were $10^{-12}$ for both solvers.

When calculating PRCs directly (see below), it is necessary to fix a reference point  on the unperturbed periodic orbit specifying the start (and end) of one oscillatory cycle and corresponding to time $t = 0$, phase $\theta = 0$. For burst PRC calculations, it is also important to designate reference points specifying the start and end of the active spiking segment of a burst cycle. These specifications are typically implemented using event functions, say $g(x), x \in \RealN{n}$, such that $g(x^{*}) = 0$ (and optionally, $g'(x^{*}) > 0$ or $g'(x^{*}) < 0$) at the desired reference point. There are several reasonable choices of event functions to specify initiation and termination of active spiking (see \eg \cite{Tien:2008}). We chose the minima and maxima of the slow variable $h$ as markers of active spiking initiation and termination, respectively; the start of active spiking was also taken to be the start of the full burst oscillation. Events were found numerically with an error tolerance of $10^{-10}$ or smaller.

\subsubsection{Infinitesimal BPRCs}
\label{sec:infinitesimal_BPRC}
We calculated the infinitesimal BPRC for the full HR model via the adjoint method of Ermentrout and Kopell \cite{Ermentrout:1990, Ermentrout:1991}, following the numerical implementation of Govaerts and Sautois \cite{Govaerts:2006}. The BPRC calculated this way corresponds to infinitesimal excitatory perturbations; inhibitory linear BPRCs are obtained by reflecting across the $\theta$-axis. The algorithmic parameters for AUTO were: collocation points, 4; number of intervals, 500; maximum step size, $10^{-4}$, initial step size, $10^{-6}$; maximum number of steps, 15. The resolution of the BPRC was 2500 points. The exact location of each point on the BPRC was determined automatically by continuation routines in AUTO, and the resulting curve was quite smooth. %Note that the use of a continuation routine to find the adjoint PRC means that the calculated curve is for a slightly different parameter value of the continuation parameter (always the applied current) than for the model parameters used to generate the reference burst trajectory. However, the resulting PRC was always for a parameter value within $10^{-5}$ of the reference parameters, so the correspondence with the true PRC for the reference parameter model should be very close.

\subsubsection{Direct PRCs}
\label{sec:direct_prcs}
We calculated direct BPRCs for the full HR model using a method we term {\it spike injection}, which closely models the biophysics of synaptic input to a postsynaptic neuron due to a single spike in a presynaptic neuron. The voltage profile of {\it stereotypical spike} is used as input to a differential equation model of synaptic transmission, and the output of the synapse model is incorporated into the equations governing the membrane voltage of the postsynaptic neuron. (For the HR model, the stereotypical spike was taken to be the third spike of the HR model's reference burst.) The resulting perturbation to the model postsynaptic neuron has a brief, nonzero duration and a distinct shape similar to the voltage profile of a single action potential. Unlike the simulated perturbation in the infinitesimal PRC calculation, whether spike injection perturbation acts to increase or decrease the voltage variable of the model postsynaptic neuron depends on the current state of postsynaptic neuron.

A synaptic current term $\Isyn$ was added to equation for $\dot{V}$:

%%%%%%%%%%%%%%%%%%
%	HR V with synaptic term
%%%%%%%%%%%%%%%%%%  
\begin{equation}
\dot{V} = n - aV^3 + bV^2 - h + I  + \Isyn \label{eqn:HR_V_syn}
\end{equation} 

\noindent The form of the $\Isyn$ term modeled the activity of a synapse with graded release of neurotransmitter \cite{Destexhe:1994, Destexhe:1998}, in the absence of any activity-modulated facilitation or depression:

%%%%%%%%%%%%%%%%%%
%	synapse equations
%%%%%%%%%%%%%%%%%%  
\begin{eqnarray}
\Isyn &=& \gSyn s(V-\vSyn) \\
\label{eqn:dot_s}
\dot{s} &=& \aSyn T_{\infty}(V_\mathrm{pre})(1-s) - \bSyn s \\
T_{\infty}(V_\mathrm{pre}) &=& (1+\exp(-(V_\mathrm{pre} - V_p)/K_p))^{-1}
\end{eqnarray}

\noindent Here $s$ represents the level of neurotransmitter released into the synaptic cleft and actively affecting the post-synaptic cell. The rate of neurotransmitter release is given by Equation~(\ref{eqn:dot_s}), and depends on the concentration of neurotransmitter already released as well as the voltage of the presynaptic cell, \ie the instantaneous strength of the incoming spike that constitutes the perturbation.

The parameter $\gSyn$ represents the synaptic coupling conductance, which is the usual measure of synaptic strength. The incoming spike was taken to be the voltage profile of the model's {\it stereotypical spike} $V_{\mathrm{spike}}(t)$, aligned with the chosen start time of the perturbation. Prior to the time of perturbation, $\gSyn = 0$ and $s = 0$. During the perturbation, the vector field defining the neuron model, including the synaptic input equations, was integrated with $V_{\mathrm{pre}}=V_{\mathrm{spike}}(t)$ and $\gSyn$ set to a fixed positive value.  If the numerical integration required voltage values for the stereotypical spike profile at times not given in the spike profile, these values were determined using linear interpolation. At the end of the stereotypical spike duration, the perturbation was considered complete and $g_{syn} = 0$ for the rest of the integration.

The synaptic reversal potential $\vSyn$ determined whether the action of the synaptic input was excitatory or inhibitory. For excitatory perturbations, $\vSyn = 1.37$, and for inhibitory perturbations,$\vSyn = -1.63$.  These parameter values for the HR model synapse were adapted from neurophysiological data for neonatal rat synapses \cite{Raastad:1998} and corresponded to voltages at 0 and 88.32 percent of the total spike height for the HR model's stereotypical spike. The half-activation point for the voltage-dependent neurotransmitter release, $V_{p}$, was also derived from neurophysiological data and corresponded to a voltage at 55 percent of the total spike height for the stereotypical spike. 

The response of the simulated synapse to the injected pre-synaptic action potential spike is sensitive to $K_{p}$, the steepness of the voltage response curve for neurotransmitter release, and the forward and backward rate constants for neurotransmitter release, $\aSyn$ and $\bSyn$. These parameters were chosen to ensure that the response to an injected spike was relatively fast and decayed before the typical interspike interval elapsed: $\aSyn$ was set to $110/ \mathrm{ISI}$, where ISI is the shortest interspike interval of the model's reference burst. $\bSyn$ was set to $0.1\cdot \aSyn$. The parameter $K_{p}$ was fixed such that it was in a biophysically reasonable range; specifically, the value was chosen such that the ratio of $K_{p}$ to stereotypical spike height for the HR model equal to the ratio of the original $K_{p}$ value \cite{Destexhe:1998}  to the stereotypical spike height of a similar, biophysically derived model of a bursting neuron in the preB\"{o}tzinger complex \cite{Butera:1999}. The specific parameter values were $\aSyn = 9.304, \bSyn = 0.9304, K_{p} = 0.227, V_{p} = 0.239$.

Spike injection perturbations were applied at different times $t \in [0,T]$ to the model neuron on its reference orbit $\Gamma$. The BPRCs were calculated by perturbing at 500 time points, with 60\% of the time points evenly spaced in the active segment of $\Gamma$ and the remainder evenly spaced in the quiescent segment. Direct BPRCs were calculated at strengths ranging over six orders of magnitude for both excitatory and inhibitory perturbations, \ie $\gSyn = 10^{-4} - 10$. This range encompassed functionally `very weak' to functionally `very strong' perturbations.

\section{Burst phase response}
\label{sec:burst_phase_response}
\begin{figure}[!ht]
     	\begin{center}
		\includegraphics[width= 3.5in]{\figpath/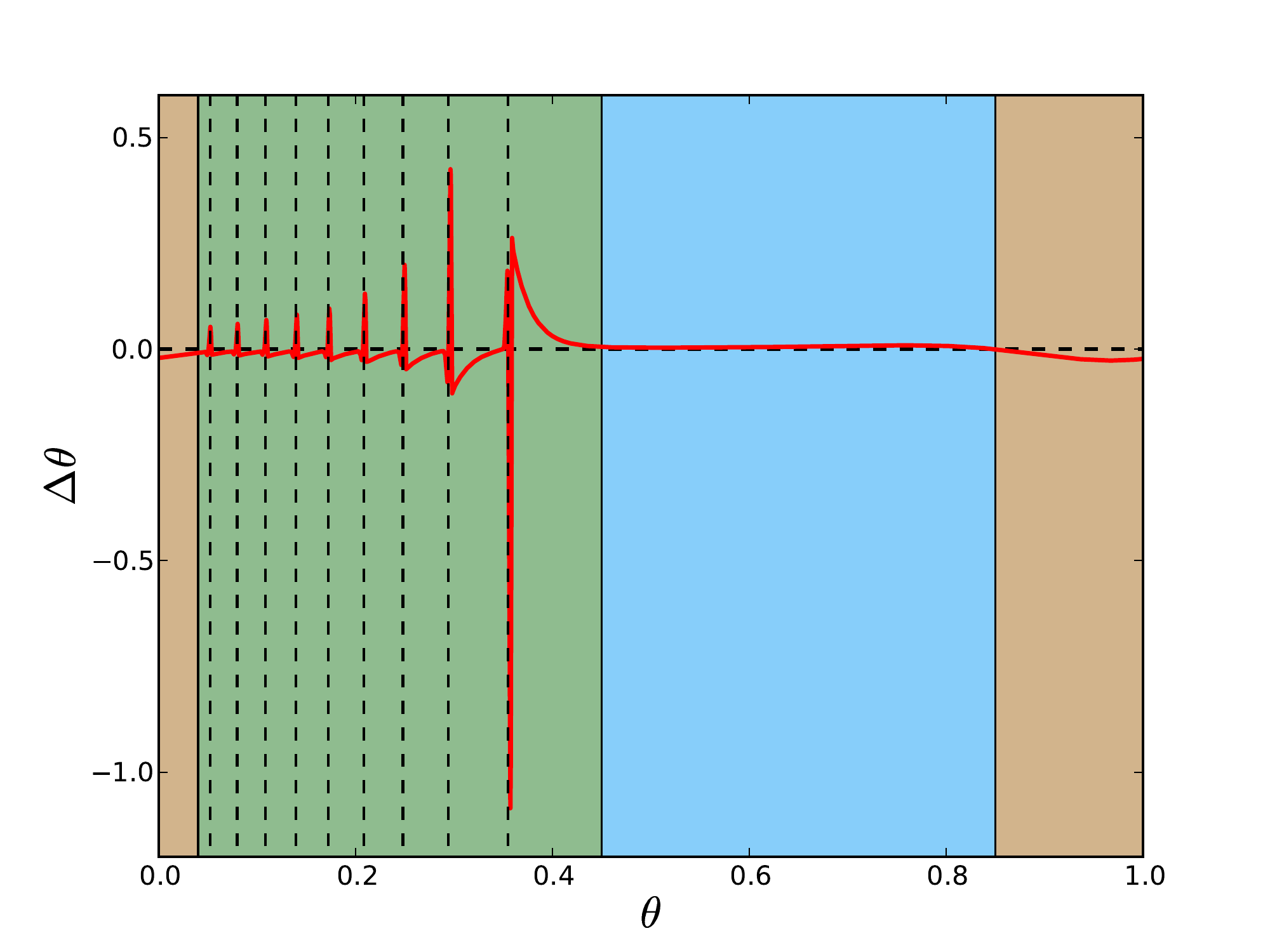}
	\end{center}
\caption[Infinitesimal BPRC]{Infinitesimal burst phase response curve for the full HR model. Positive $\Delta\theta$ indicates phase delay; negative $\Delta\theta$ indicates phase advancement. The dashed vertical lines indicate phases corresponding to the peaks of spikes in the reference burst orbit. Background shading denotes segments of distinct burst phase response behavior: Segment I, sea green; Segment II, light blue; Segment III, tan.}
\label{fig:AdjointPRC}
\end{figure}

Figure \ref{fig:AdjointPRC} shows the infinitesimal phase response curve of the full HR model, while Figures \ref{fig:HRIe_BPRCs} and \ref{fig:HRIi_BPRCs} depict direct excitatory and inhibitory BPRCs, respectively, computed for a range of $\gSyn$ values. In each of the figures, dashed vertical lines indicate the position in phase of individual spike maxima in the reference burst orbit. The direct excitatory and inhibitory BPRCs for $\gSyn < 0.01$ had shapes identical (up to scaling by a constant amount) to the BPRC presented in Figure \ref{fig:HRIe_BPRCs} (a) and are not shown.

%%%%%%%%%%%%%%%%%%
%	Infinitesimal PRC
%%%%%%%%%%%%%%%%%%  

%%%%%%%%%%%%%%%%%%
%	Excitatory BPRCs
%%%%%%%%%%%%%%%%%%  
\newcommand{\gwidth}{2.95in}
\newcommand{\mpwidth}{2.95in}
\begin{figure}[!ht]
  \begin{minipage}{\mpwidth}
      \begin{center}

    \includegraphics[width= \gwidth]{\figpath/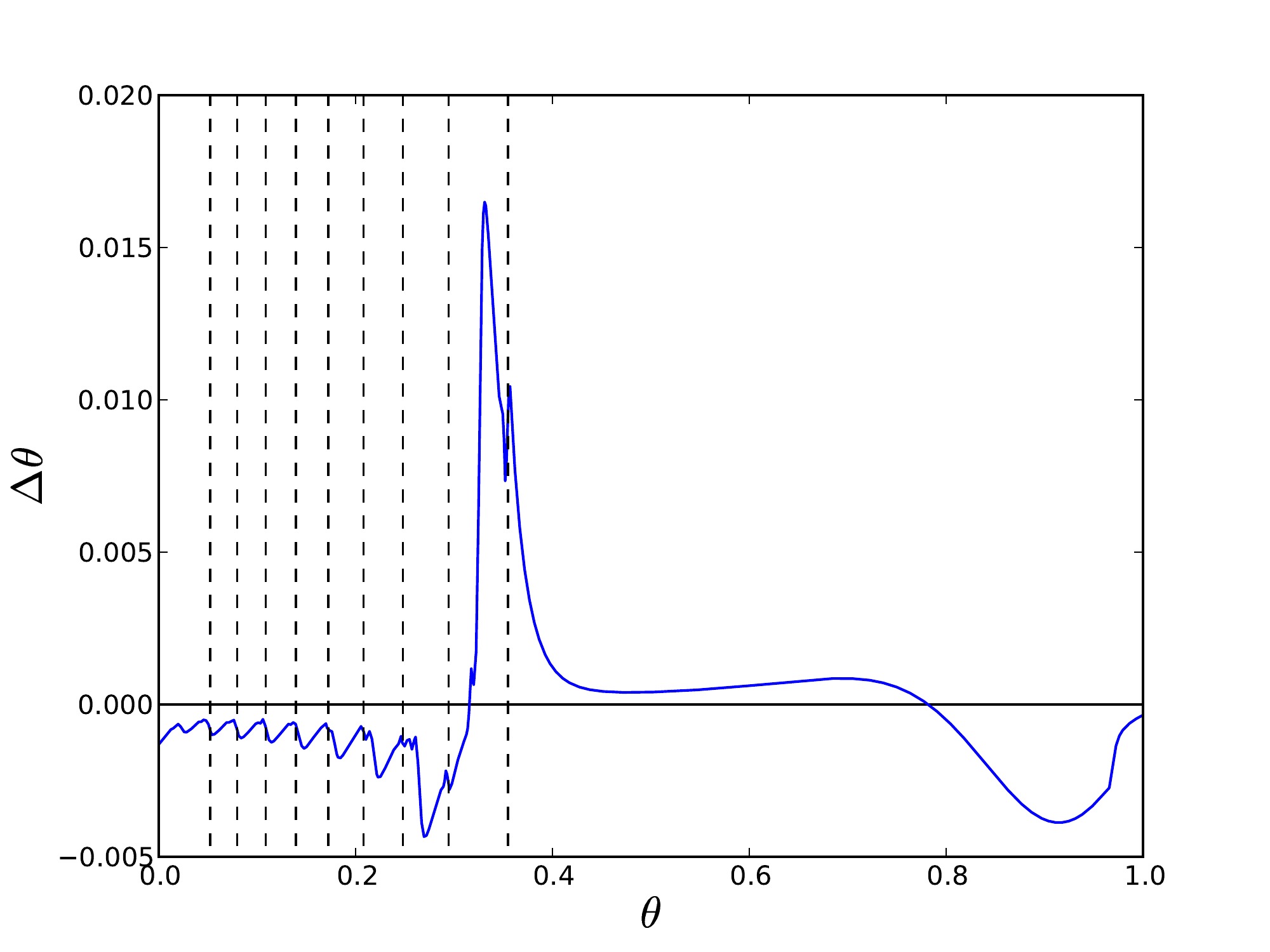}
        (a) $\gSyn = 0.01$
    \end{center}
  \end{minipage}%
  \begin{minipage}{\mpwidth}
      \begin{center}

    \includegraphics[width= \gwidth]{\figpath/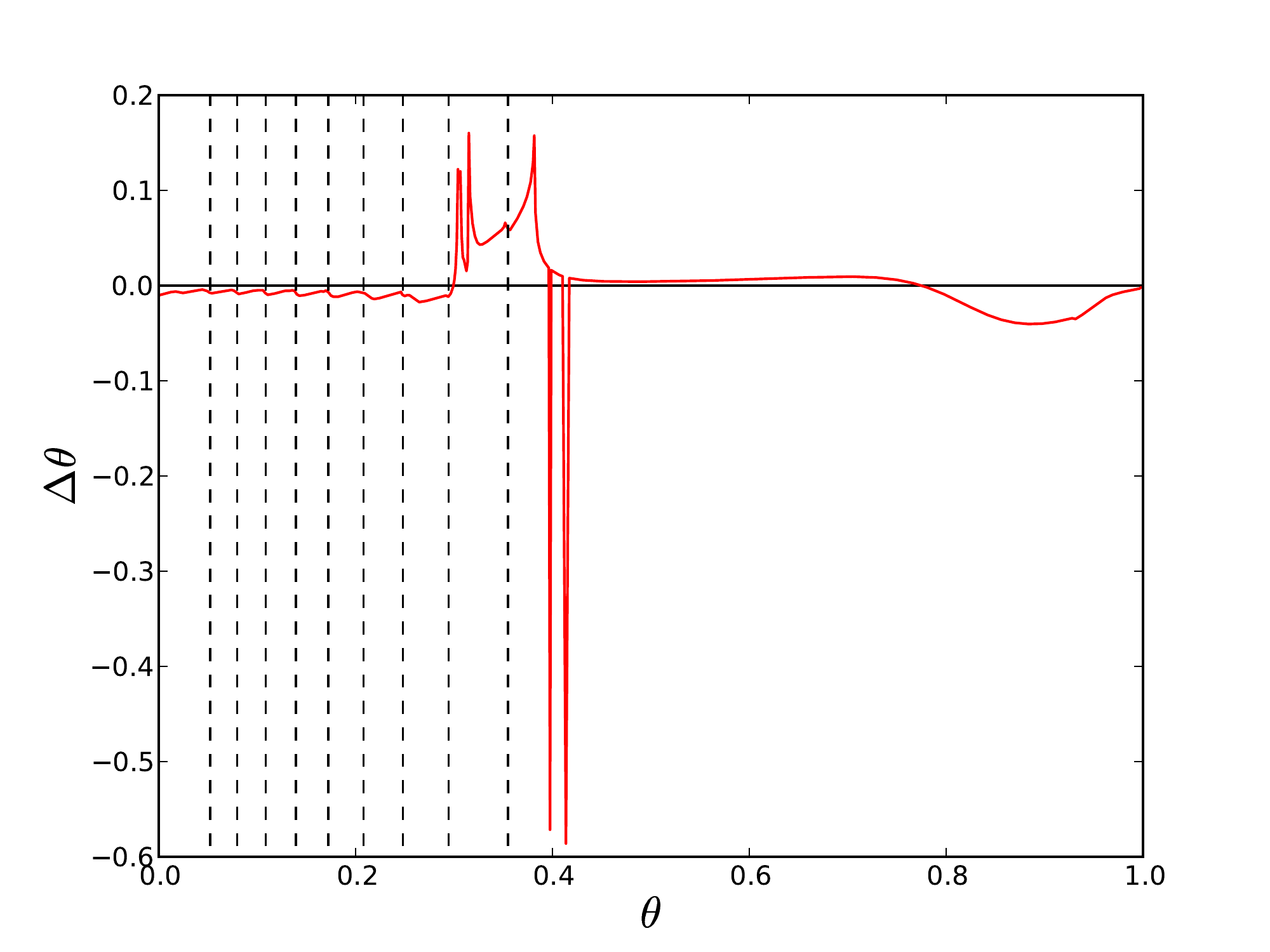}
        (b) $\gSyn = 0.1$
    \end{center}
  \end{minipage}
  \begin{minipage}{\mpwidth}
      \begin{center}

    \includegraphics[width= \gwidth]{\figpath/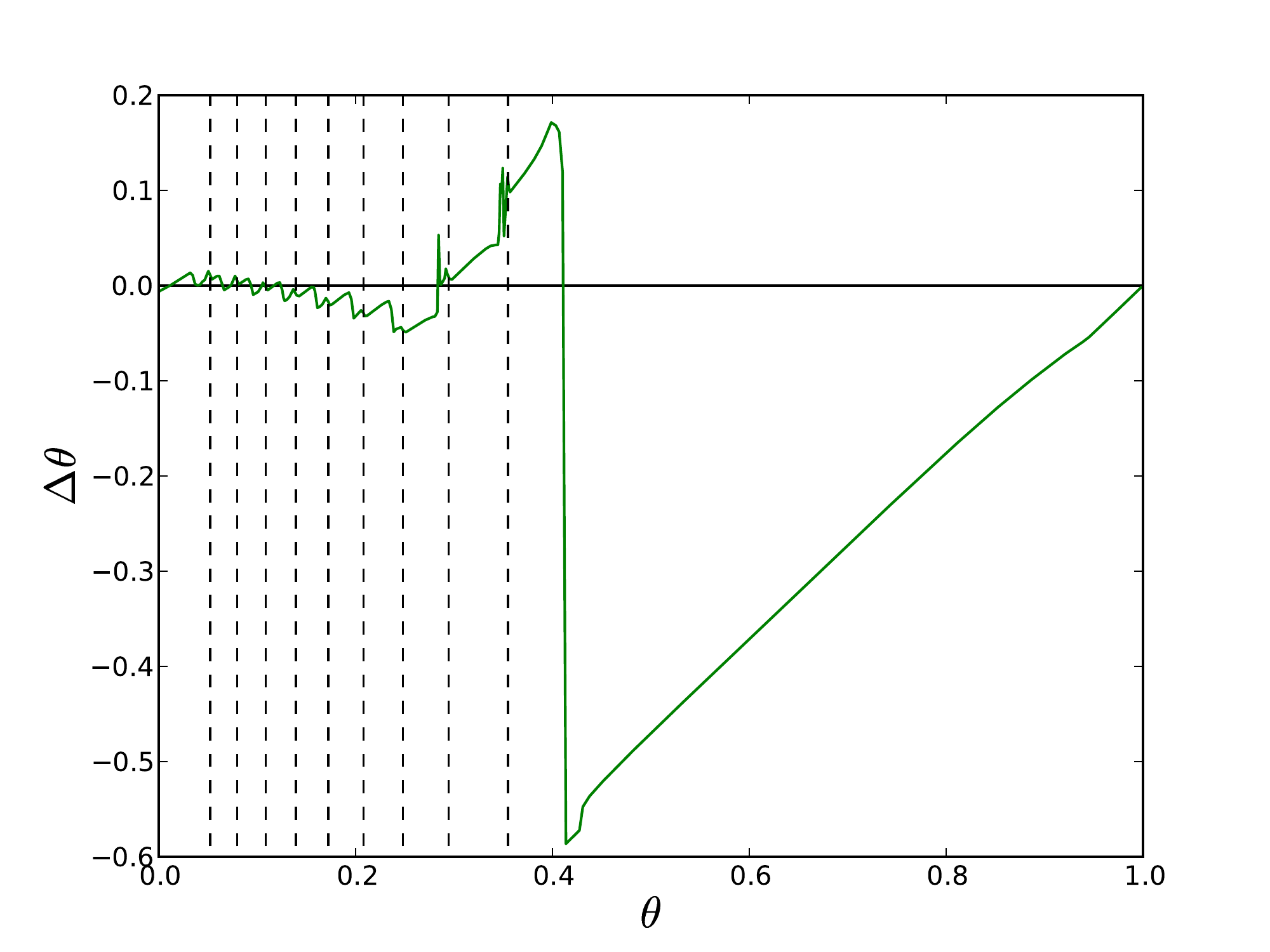}
        (c) $\gSyn = 1$
    \end{center}
  \end{minipage}%
  \begin{minipage}{\mpwidth}
      \begin{center}

    \includegraphics[width= \gwidth]{\figpath/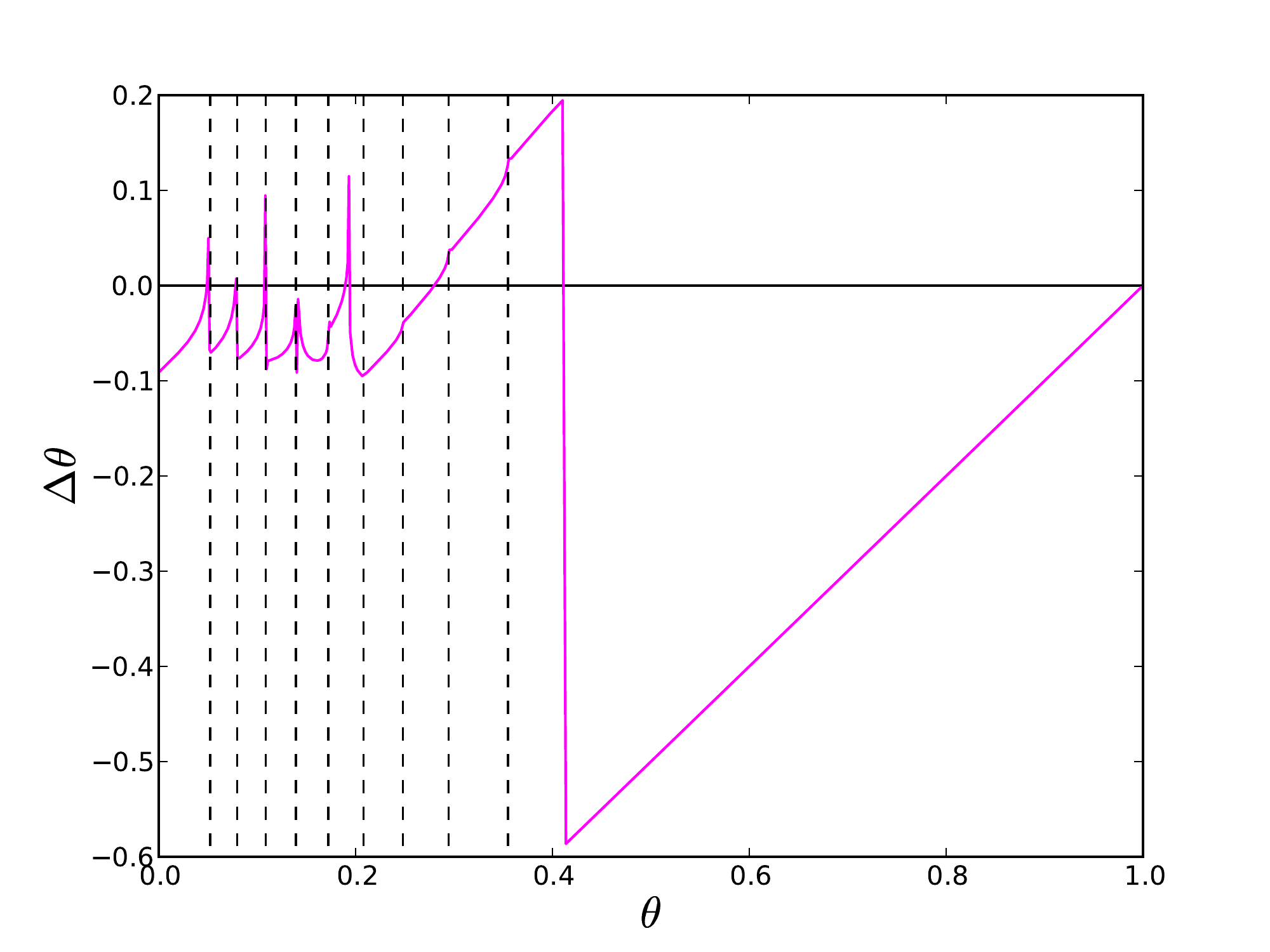}
        (d) $\gSyn = 10$
    \end{center}
  \end{minipage}
\caption[HRIe  Burst PRCs ]{Direct (spike injection) BPRCs for the full HR model, using an excitatory synapse, for a range of perturbation strengths ($\gSyn$ values). Dashed vertical lines indicate phases corresponding to the peaks of spikes in the reference burst orbit.}% From (a)--(f),  \gSyn\ ranges from $10^{-4}$ to $10$.}
\label{fig:HRIe_BPRCs}
\end{figure}

%%%%%%%%%%%%%%%%%%
%	Inhibitory BPRCs
%%%%%%%%%%%%%%%%%%  
\begin{figure}[!ht]
  \begin{minipage}{\mpwidth}
      \begin{center}

    \includegraphics[width= \gwidth]{\figpath/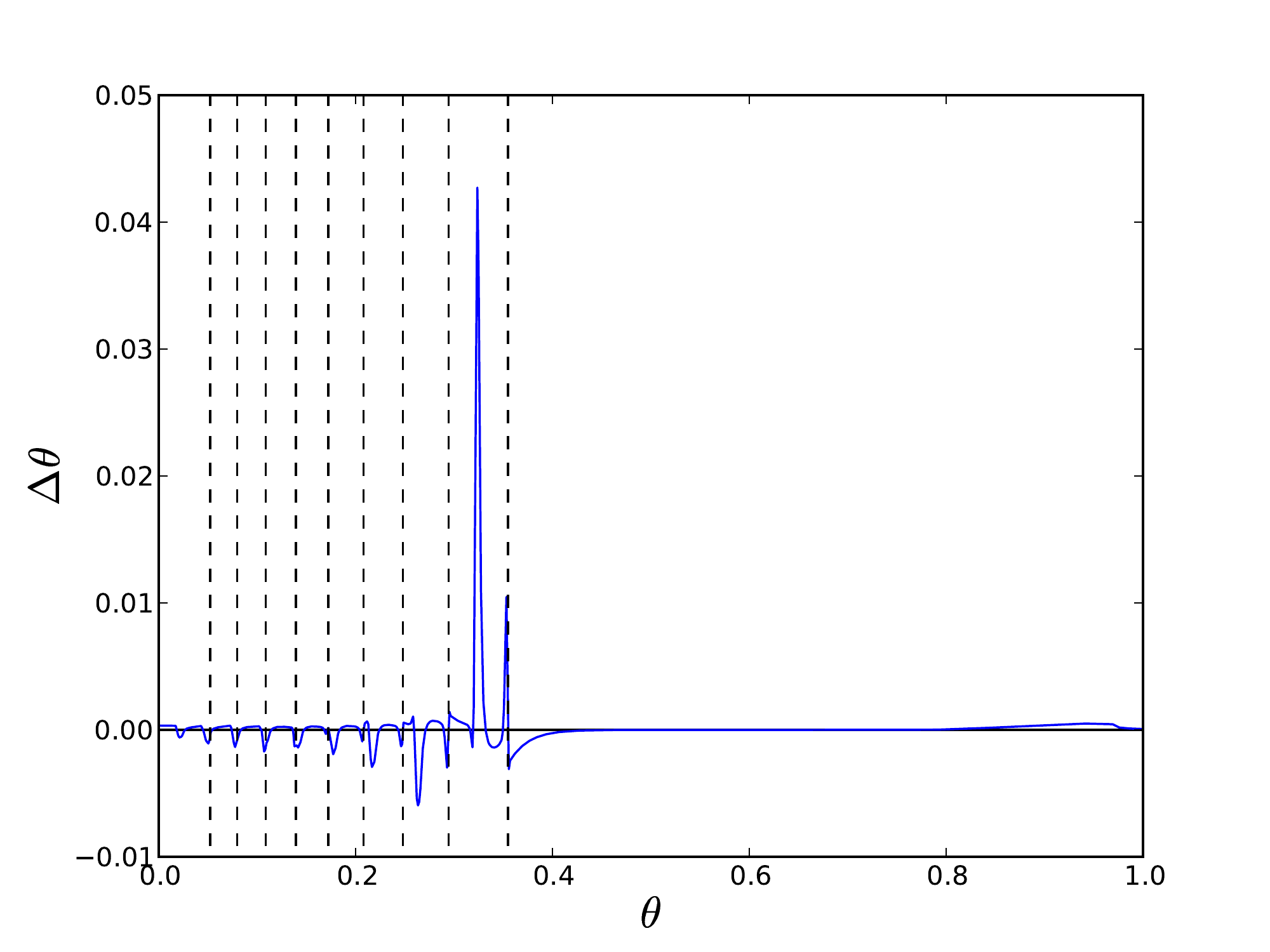}
        (a) $\gSyn = 0.01$
    \end{center}
  \end{minipage}%
  \begin{minipage}{\mpwidth}
      \begin{center}

    \includegraphics[width= \gwidth]{\figpath/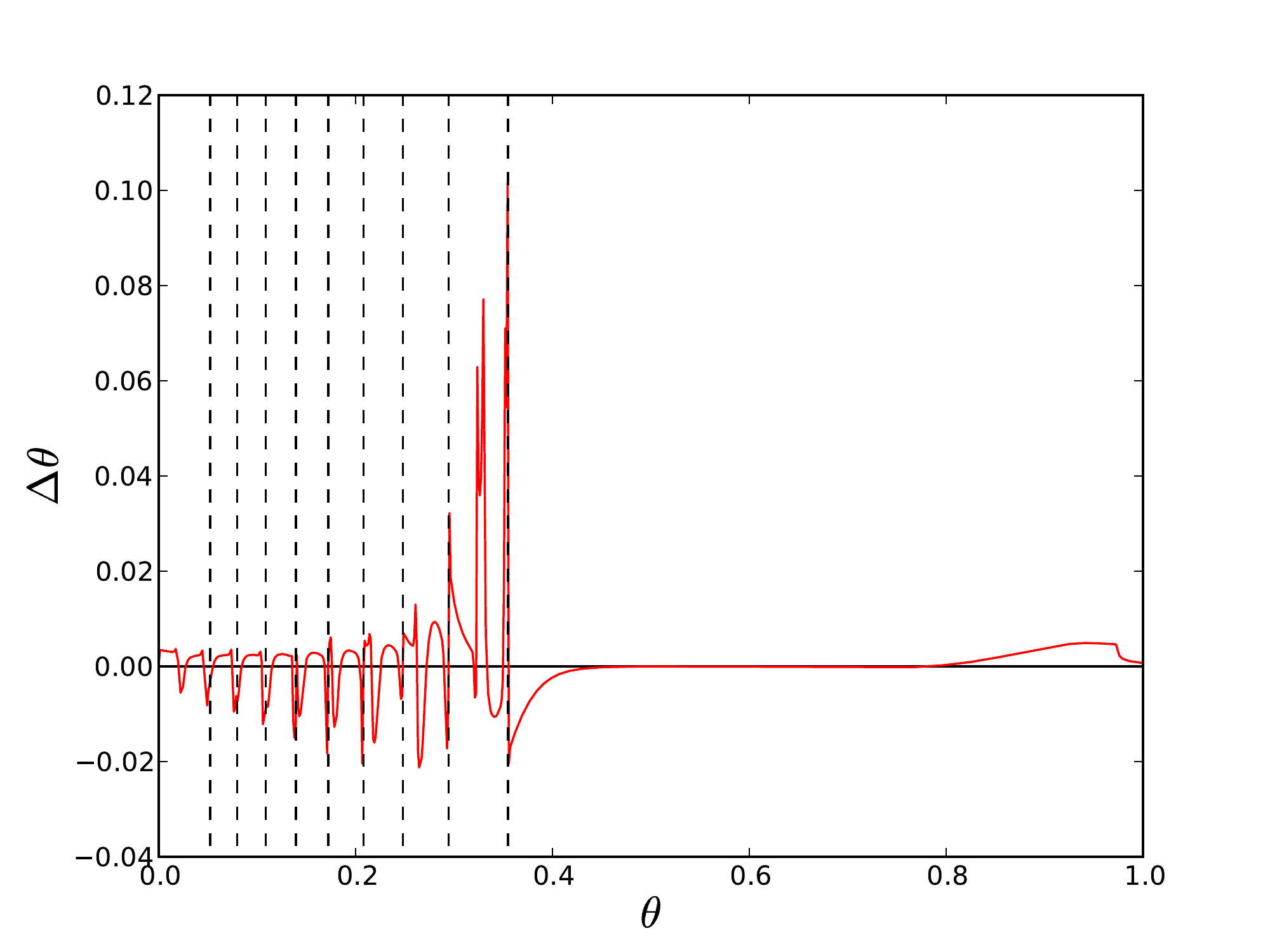}
        (b) $\gSyn = 0.1$
    \end{center}
  \end{minipage}
  \begin{minipage}{\mpwidth}
      \begin{center}

    \includegraphics[width= \gwidth]{\figpath/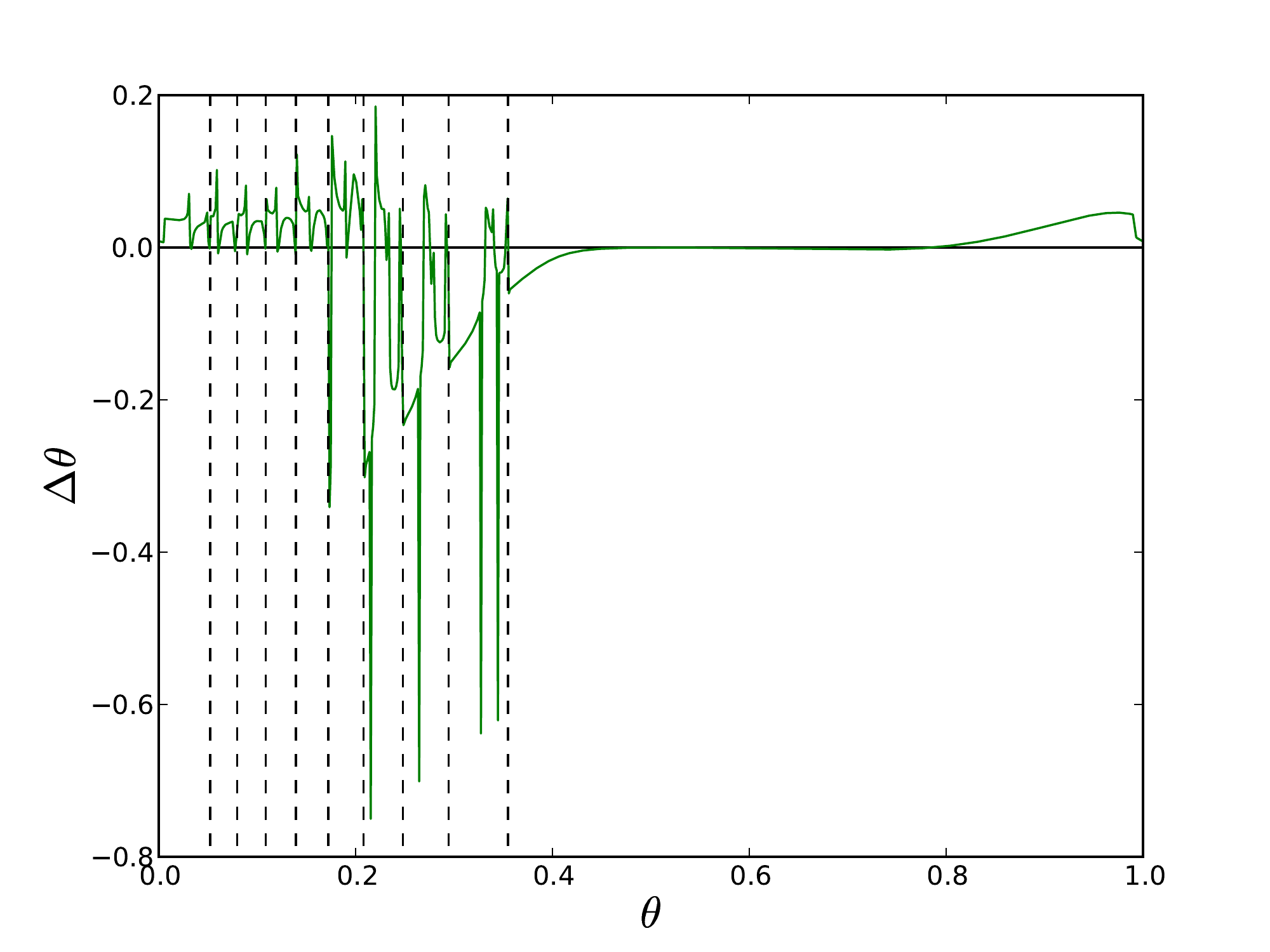}
        (c) $\gSyn = 1$
    \end{center}
  \end{minipage}%
  \begin{minipage}{\mpwidth}
      \begin{center}

    \includegraphics[width= \gwidth]{\figpath/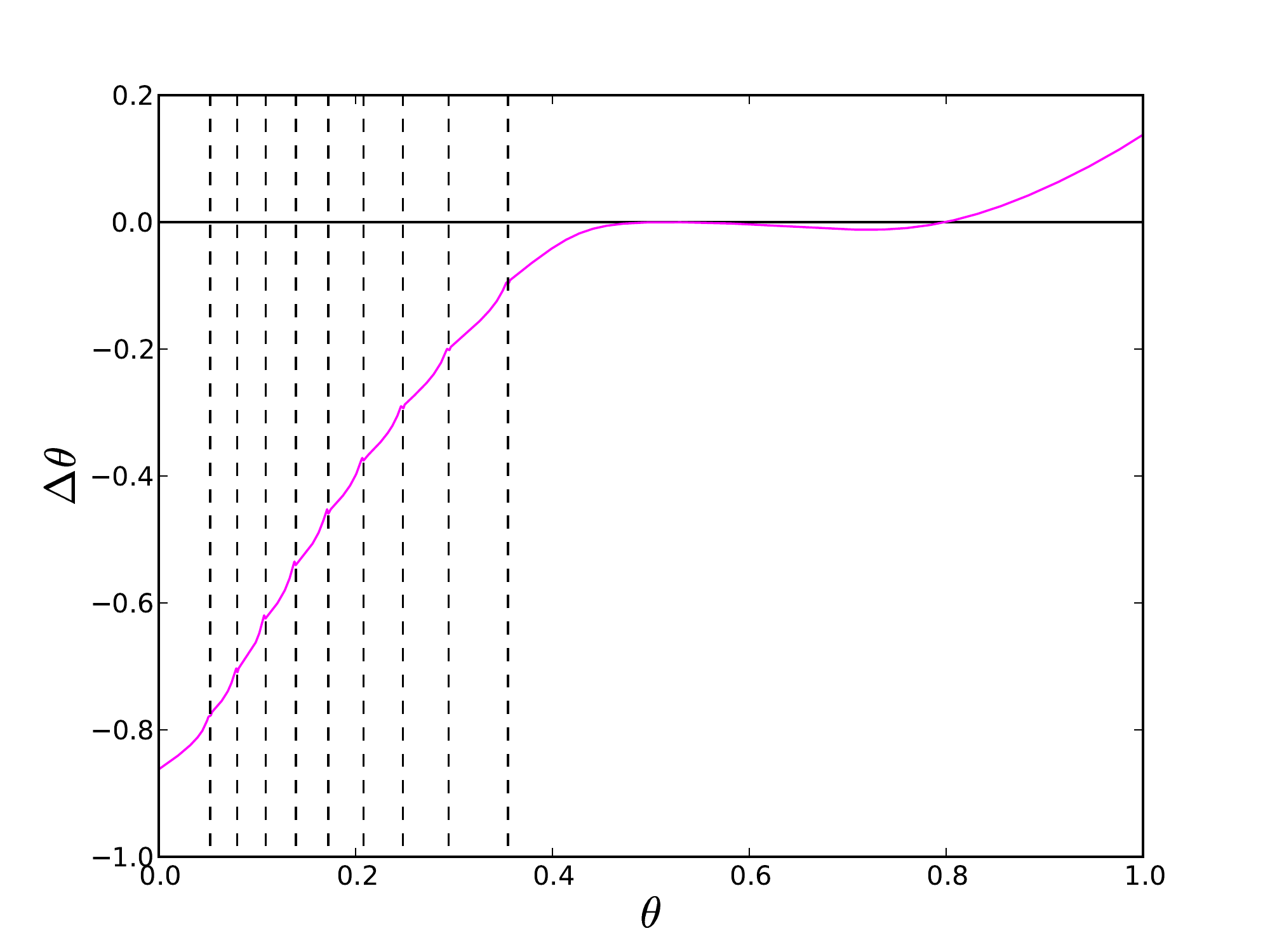}
        (d) $\gSyn = 10$
    \end{center}
  \end{minipage}
\caption[HRIi  Burst PRCs ]{Direct (spike injection) BPRCs for the full HR model, using an inhibitory synapse, for a range of perturbation strengths ($\gSyn$ values). Dashed vertical lines indicate phases corresponding to the peaks of spikes in the reference burst orbit.}% From (a)--(f),  \gSyn\ ranges from $10^{-4}$ to $10$.}
\label{fig:HRIi_BPRCs}
\end{figure}

%%%%%%%%%%%%%%%%%%
%	SNRCs
%%%%%%%%%%%%%%%%%%  
\begin{figure}[ht]
 \begin{minipage}{3in}
     \begin{center}
    \includegraphics[width= 3in]{\figpath/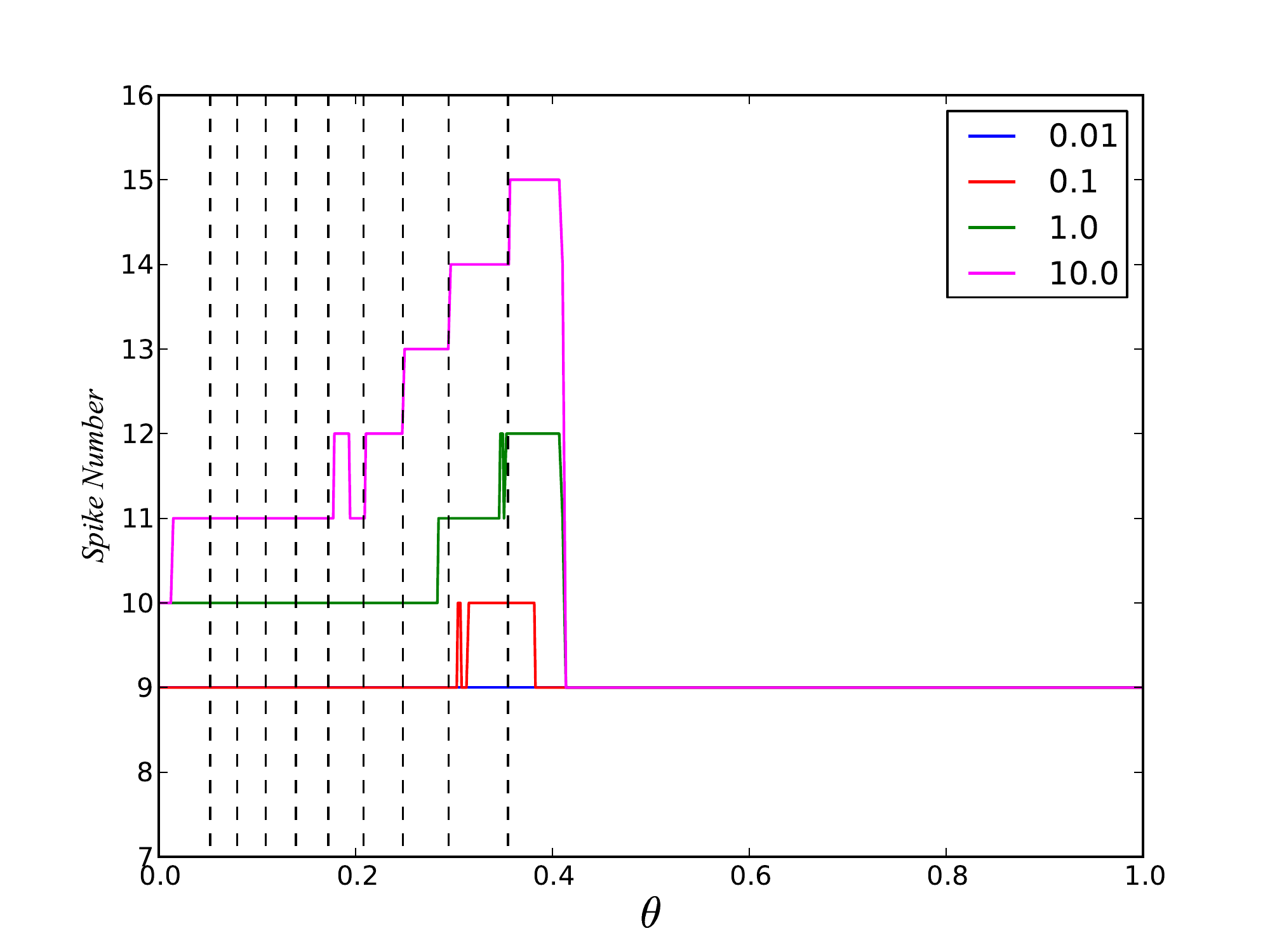}
       (a)
    \end{center}
  \end{minipage}%
  \begin{minipage}{3in}
      \begin{center}
    \includegraphics[width= 3in]{\figpath/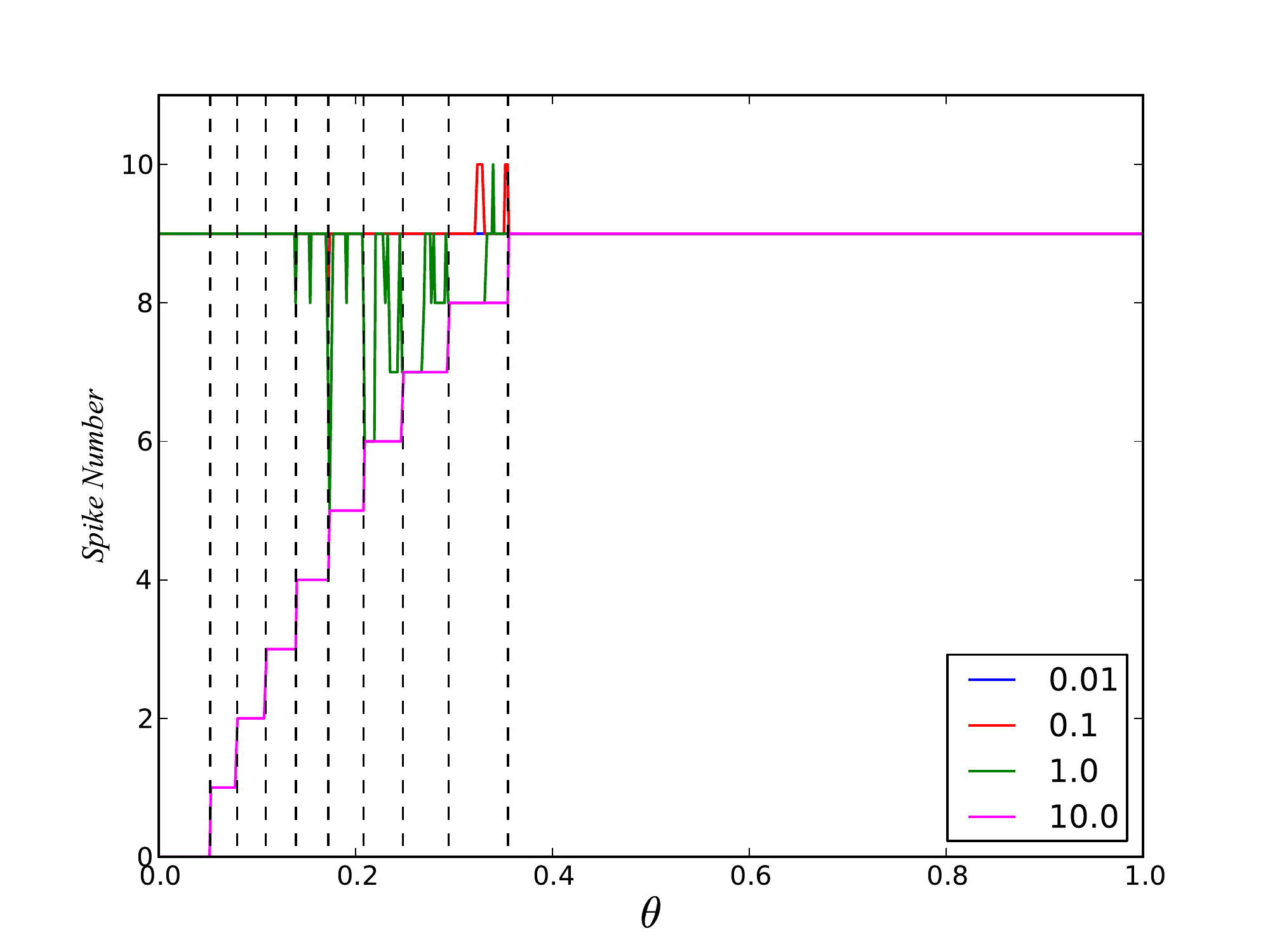}
        (b)
    \end{center}
  \end{minipage}
\caption[HR SNRCs]{Spike number response curves (zeroth order) for the HR model. Dashed vertical lines indicate phases corresponding to the peaks of spikes in the reference burst orbit. SNRC colors and perturbation strengths match those of the BPRCs depicted in Figures \ref{fig:HRIe_BPRCs} and \ref{fig:HRIi_BPRCs}. (a) Excitatory SNRCs. (b) Inhibitory SNRCs.}
\label{fig:HR_SNRCs}
\end{figure}

As mentioned above, phase response curves for spiking neuronal models generally fall into one of two categories: Type I, if they are either all positive or all negative and are associated with a saddle node of periodic orbits bifurcation, or Type II, if they have both positive and negative part and are associated with a subcritical Hopf bifurcation. Either type of PRC typically exhibits a single large peak (or trough) of phase response, often but not always aligned with the phase of the voltage maximum corresponding to the peak of the action potential. (See, for example, \cite{ Brown:2004b, Ermentrout:1996, Govaerts:2006}). None of the BPRCs calculated for the HR model possess such a simple configuration. 

Instead, each of BPRCs has a striking, relatively complex shape; all of them may be visually decomposed into three distinct segments of phase response: Segment I, an initial, `spiky' region of rapidly varying phase response, marked in sea green in Figure \ref{fig:AdjointPRC}; Segment II, a broad region of relatively flat or constant slope phase response, marked in light blue in Figure \ref{fig:AdjointPRC}; and Segment III, a relatively narrow region, marked in tan in Figure \ref{fig:AdjointPRC}. Segment I corresponds to the active spiking segment of the perturbed reference burst orbit and roughly spans the range of phases $[0.04, 0.45]$. Segment II corresponds to the quiescent, hyperpolarized segment of the perturbed reference burst orbit and roughly spans the range of phases $[0.45--0.85]$. Segment III spans the remaining phases $[0.85,0.04]$, and covers the segment of the perturbed reference burst orbit close to the initiation of a new round of active spiking. This decomposition is roughly the same for all of the BPRCs---infinitesimal or direct, excitatory or inhibitory, weak or strong perturbations---and is strongly related to the multiple time-scale dynamics of bursting, as we explore in the following two sections.

Within Segment I, there is a close association between the locations (in phase) of spike maxima in the reference burst orbit and the positions of local extrema of phase response, and the magnitude of phase response increases with proximity to the end of the segment and the termination of active spiking in the reference burst orbit. Also rather remarkable is the high degree of sensitivity to the timing of perturbations shown in the phase response in Segment I: As seen in nearly all of the BPRCs, a perturbation of set strength at a given phase $\theta_1$ may result in a significant phase delay, while perturbation at phase $\theta_2 > \theta_1$, $\| \theta_2 - \theta_1 \| \leq 0.01$, induces a large phase advance, and perturbation at phase $\theta_3 > \theta_2$,  $\| \theta_3 - \theta_2 \| \leq 0.01$, again results in substantial phase delay. What explains the distinctive shape of the BPRCs, and what causes this phase response sensitivity? 

Before turning to our analyses of its phase response dynamics, we make a few additional observations about the direct BPRCs we calculated for the HR model. There are broad, significant similarities between the infinitesimal and direct BPRCs, but we emphasize the various discrepancies between the shapes of BPRCs calculated via direct methods and the BPRC shapes one would expect from theoretical considerations and by extrapolation from the computed shape of the infinitesimal BPRC. We believe these differences merit further mathematical investigation as well as considered attention from a biological perspective.

Cursory inspection of Figures \ref{fig:AdjointPRC}, \ref{fig:HRIe_BPRCs} and \ref{fig:HRIi_BPRCs} reveals that even at low perturbation strengths, the shapes of both excitatory and inhibitory direct BPRCs are far from identical to the shape of the infinitesimal BPRC. Though the infinitesimal and both direct BPRCs are decomposable into three qualitatively different segments of phase response, there is some difference in the details of phase response within those segments. Within Segment I, the alignment of local extrema of phase response, \ie phase delays, with the phases of spikes in the reference burst orbit is very tight for the infinitesimal BPRC, but somewhat looser for both the excitatory and inhibitory direct BPRCs. Near spike phases in the excitatory direct BPRCs, the phase response is phase advancement (albeit locally reduced phase advancement), rather than phase delay as in the infinitesimal BPRC case. As shown in Figure \ref{fig:HRIe_BPRCs} (a), for excitatory perturbation at strengths $\gSyn \leq 0.01$, there is no reversal of the direction of phase response near the termination of active spiking at the end of Segment I. At higher strengths, as shown in Figure \ref{fig:HRIe_BPRCs} (c) and (d),  there is substantial phase advancement that occurs immediately after the last spike of the reference burst orbit, at the end of Segment I and beginning of Segment II. The magnitude of this phase advancement decreases almost linearly over the course of Segment II. 

If the perturbation modeled by spike injection were well-approximated by the infinitesimal BPRC, and assuming that the isochrons of the HR system were approximately rectilinear near the burst periodic orbit, one would expect that the inhibitory direct BPRCs would have a shape very similar to the infinitesimal BPRC, modulo reflection across the $\theta$-axis.
That is roughly the case for the inhibitory direct BPRCs at perturbation strengths $\gSyn \leq 0.01$, \eg Figure \ref{fig:HRIi_BPRCs} (a). This similarity includes the reversal of phase response at the end of Segment I, but this reversal occurs earlier, closer to the phase of the penultimate spike of the reference burst orbit, than for the infinitesimal BPRC. The greater similitude of the inhibitory, rather than excitatory, direct BPRCs to the infinitesimal BPRC is contrary to expectations, and it suggests that even for weak synaptic coupling, the phase resetting dynamics produced by realistic inputs to bursting neurons may be considerably more complicated than the linearized approximation of the infinitesimal BPRC would indicate.

For direct BPRCs computed for weaker perturbation strengths, \ie $\gSyn \leq 0.01$, differences from the infinitesimal BPRC may be attributed in part to the character of the spike injection perturbation, which has finite, nonzero duration and a shape distinct from a delta function pulse. Another factor is that the reversal potential of the model synapse determines whether the direct BPRC perturbation is inhibitory or excitatory. For either excitation or inhibition, whether an injected spike acts to increase or decrease the membrane voltage variable depends on the state of the model neuron, specifically the current value of the membrane voltage variable relative to the synaptic reversal potential, whereas for the infinitesimal BPRC an excitatory perturbation always acts to increase the membrane voltage variable by a small amount, and an inhibitory perturbation always has the opposite effect.

Deviation from the shape of the infinitesimal BPRC becomes more pronounced for both excitatory and inhibitory direct BPRCs as perturbation strength increases. As well as changing shape, the direct BPRCs for $\gSyn \geq 0.1$ exhibit very large phase resetting at very specific perturbation phases, \eg near $\theta\approx 0.4$ in Figure \ref{fig:HRIe_BPRCs} (b) (excitatory) and near $\theta\approx 0.2$ in Figure \ref{fig:HRIe_BPRCs} (c) (inhibitory). As perturbation strength increases, perturbation at more phases produce incommensurately large phase advancement or delay; as shown in Figure \ref{fig:HRIe_BPRCs} (d), for example, when $\gSyn = 10$ an inhibitory perturbation at any phase $\theta \in [0.0,0.2]$ results in a phase advance of nearly half a cycle or more. 

What accounts for these large phase resettings? The SNRCs presented in Figure \ref{fig:HR_SNRCs} provide the essential clues: large phase shifts are closely connected with changes  in the number of spikes in the perturbed burst. Increases in the number of spikes in the active spiking segment of the perturbed burst typically delays the onset of the next burst cycle; this is the usual effect of strong excitatory perturbations. Reductions in the number of spikes during the active spiking segment typically advances the onset of the next burst cycle, the usual effect of strong inhibitory perturbations. We note, however, that for specific values of $\gSyn$ and specific phases of perturbation, inhibitory perturbations may act to add spikes. Here again there may be great phase sensitivity, with a perturbation of a given strength at several nearby phases producing spike number changes at some phases and no change at others.

Thus there is a range of perturbation strengths for which direct BPRCs are fairly well approximated by the infinitesimal BPRC (approximately $\gSyn \leq 0.01$), and a range of perturbation strengths for which the direct BPRCs deviate strongly from the infinitesimal BPRC (approximately $\gSyn \geq 1$)---including, in particular, deviations such that the corresponding SNRCs are not flat, \ie the spike number changes. We call the former range the  {\it weak perturbation regime}, and refer to the latter range as the {\it strong perturbation regime}. In the following section, we examine in detail the reasons for the peculiar shape of BPRCs in the weak perturbation regime; we consider BPRCs in the strong perturbation regime, and mechanisms of spike number change, in Section \ref{sec:strong_perturbation}.

\section{Weak perturbation regime BPRCs}
\label{sec:weak_perturbation}
 
In our analysis of burst phase response in the weak perturbation regime, we focus on the infinitesimal BPRC and consider an excitatory perturbation that acts instantaneously to increase the value of the membrane voltage variable, $V$, by a small but nonzero amount. We do not specifically investigate the shapes of the direct BPRCs in the weak perturbation regime, but our analysis for the infinitesimal BPRC should be approximately correct for the direct BPRCs as well. We first address why the BPRCs visually comprise three distinct segments associated with either spiking or quiescence during the reference burst cycle, treating the quiescent segments (Segments II and III) and active spiking segment (Segment I) separately. To do this we employ the fast-slow dissection of the burst dynamics of the HR model. 

\subsection{Quiescent segment analysis}
\label{sec:quiescent_segment}

BPRC Segment II corresponds to the quiescent portion of the reference burst cycle, when $h$ recovers from 2.10256601768 to 1.75415439813 (see Figures \ref{fig:HR_Burst_Traj} and \ref{fig:HR_fast_bifns}).  Following the destruction of the family of stable periodic orbits in the fast subsystem at the homoclinic bifurcation, the burst trajectory in the full HR system tracks the line of low voltage stable fixed points.  These fixed points are the only remaining stable structures in the fast subsystem, having a large domain of attraction and relatively large, negative eigenvalues so that the rate of approach of a trajectory in the full system is quite rapid. Furthermore, the vector field for the full HR system lies nearly parallel to the $V$ direction in the fast subsystem. Thus small perturbations in the $V$ direction during this segment have little effect, since the trajectory is drawn back immediately to the stable fixed point, with only a small change in $h$ and therefore only a small change in phase. Hence the phase response in Segment II, as recorded in the infinitesimal and direct BPRCs of the weak perturbation regime, is negligible.

In Segment III, the perturbed neuron approaches the end of the quiescent portion of its current burst cycle and the onset of the subsequent round of spiking. In the fast subsystem, this corresponds to nearing the saddle-node bifurcation at which the low voltage fixed point being tracked by full system trajectory disappears. At that bifurcation point, the trajectory immediately switches to following the coexistent stable periodic orbit. As the saddle-node bifurcation point gets closer, the eigenvalues of the stable fixed point diminish, reducing its rate and domain of attraction so that small perturbations may knock the full system trajectory into the basin of attraction for the stable periodic orbit. Such an event marks the early advent of the next burst cycle, and hence Segment III shows a slight phase advancement.

\subsection{Active segment analysis}

Phase response in the active segment of the burst has a more complicated structure and appearance than the phase response during quiescence. The main features we wish to explain are (1) the close association of peaks of phase response to the spike maxima in the full system trajectory, (2) the increase in amplitude of these phase response peaks closer to the end of the active segment, and (3) the large magnitude and sign reversal of the phase response at the end of the active segment. Phase portraits of the planar fast subsystem and isochron calculations are the main tools we use in our analysis, concentrating on the fast subsystem at a few representative values of the slow variable $h$. The corresponding  phase portraits provide a cross-sectional view or `snapshot' of the dynamics and phase response characteristics of the full system at particular moments in the evolution of the full burst cycle.

\subsubsection{Fast subsystem phase portraits}
\label{sec:phase_portraits}

We investigate the phase portrait of the fast subsystem at three representative cross-sectional $h$ values, 1.8, 1.95, and 2.085, which correspond to the first, sixth and ninth (final) spikes in the full HR system, respectively. At $h=1.8$ (Figure \ref{fig:PP_Iso_18}),  the fast subsystem is at the beginning of the active segment, and there exists no stable fixed point, only the stable periodic orbit (corresponding to spikes in the full system) surrounding an unstable fixed point. The value $h=1.95$ (Figure \ref{fig:PP_Iso_195}) lies in the middle of the active spiking segment, after the stable fixed point and a saddle have emerged via the saddle-node bifurcation.  For $h=2.085$ (Figure \ref{fig:PP_Iso_2085}), the fast subsystem is very close to the homoclinic bifurcation, and the saddle and the periodic orbit nearly touch.  Phase breaks down in the limit cycle at the homoclinic point; the period of the orbit becomes infinite. Past the homoclinic point, the fast subsystem has no periodic orbit,  though the full HR system may still emit one spike. Phase is undefined in the fast subsystem, though the excitability of the full system means that perturbations from the stable fixed point may follow trajectories which track `ghosts' of the stable structures which existed in the fast subsystem at lower $h$ values. 

In each of the phase portraits, the cubic $V$-nullcline is drawn with a dotted red line, and the parabolic $n$ nullcline is drawn with a dotted green line. This parabolic-cubic nullcline configuration is different from the standard slow-fast system picture in which the slow variable has a linear nullcline and periodic trajectories alternate via fast jumps between slowly following one branch of the fast nullcline to the other. The left branches of the two nullclines in the HR fast subsystem lie very close together, so that trajectories proceed very slowly in their vicinity. Near the right branch of the $n$-nullcline, away from the fast nullcline, $V$ evolves very rapidly; this is the depolarizing region of the spike. %The stable eigendirections for the saddle and stable fixed point are drawn as dashed green lines, and the unstable eigendirection of the saddle is drawn as a dashed red line. At each labeled phase, we also draw a line normal to the vector field at that point on the periodic orbit.

We plot several points about the periodic orbit, equally spaced in phases, and we assign phase 0 to the point of maximum voltage---the spike peak. Trajectories near the periodic orbit proceed clockwise around or inside of it. The nonuniform spatial distribution of points on the periodic orbit that are equally spaced in phase shows how the speed of spiking trajectories varies  in different regions of phase space. The movement is most rapid where phase points are spread farthest apart, namely during hyperpolarization after the spike maximum, between phases 0 and 0.1. It is almost as fast during depolarization before the spike maximum, from phase 0.9 to phase 0, and it is slowest at low $V$ values, between phases 0.1 and 0.8. The speed of the hyperpolarized portion of the orbit relative to the spike slows dramatically as the homoclinic point approaches, as can be seen by noting the compression of phases 0.1 to 0.8 along the left arc of the periodic orbit for $h = 2.085$ in Figure \ref{fig:PP_Iso_2085}. The proximity of the saddle point to the periodic orbit has a stark retarding effect on the trajectory. In early and middle parts of the active segment, Figures \ref{fig:PP_Iso_18} and \ref{fig:PP_Iso_195}, the spike occurs in the phase interval 0.9 to 0.1, but at the end of the active segment, this has shrunk to the interval between 0.96 and 0.04. Hence the hyperpolarized interspike intervals increase in length as the active segment of the burst progresses.

\subsubsection{Fast subsystem isochrons}
\label{sec:isochrons}

After a perturbation of  the full HR system, if bursting is not permanently silenced, the perturbed trajectory has been moved to a point $y$ in the basin of attraction of the (full system) burst periodic orbit $\Gamma$. Then $y$ lies on an isochron $\varpi(x)$ for $x \in \Gamma$, and this isochron determines the phase shift due to the perturbation. During the active segment of the burst cycle, the unperturbed trajectory tracks the stable periodic orbit (denoted $\Gamma_{h}$) in the fast subsystems associated with a narrow range of $h$ values, and the isochrons of $\Gamma_{h}$ lie near related isochrons of $\Gamma$. The repositioning of the perturbed trajectory relative to the isochrons of $\Gamma$ is reflected in its shift relative to the isochrons of $\Gamma_{h}$, and so examination of the geometry of the isochrons of the fast subsystem for several representative $h$ values presents an informative picture with respect to phase response in the full HR model. In the following sections, we refer to the segment of an isochron $\varpi(x)$ (in the fast subsystem) that lies within the periodic orbit $\Gamma_{h}$ as the {\it inner isochron} at point $x$, and the segment lying outside $\Gamma_{h}$ as the {\it outer isochron} at $x$. %We calculate isochrons in the fast subsystems for cross-sections of the active segment at $h = 1.8$, 1.95, and 2.085. 

We emphasize that the notion of phase for the full and fast subsystems are not identical, and that the isochrons for the full and fast subsystems are different mathematical objects with different meanings for the behavior of the model. Isochrons in the full system are two-dimensional manifolds with asymptotic phases relative to the full burst cycle, whereas isochrons in the fast subsystem are one-dimensional curves with asymptotic phases associated with a periodic orbit approximating a single spike of the full burst cycle. That is, a phase $\theta_{h}$ (and its associated isochron) in the fast subsystem is not the same as $\theta$ in the full system, even if they have the same scalar value, since the former refers to phase in a smaller subset (a single spike) of the periodic orbit associated with the latter. However, phase shifts $\Delta\theta_{h}$ and $\Delta\theta$ are directly related, since advancement or delay along a single spike implies a similar shift along the full burst cycle. The fast subsystems we consider are obtained directly via fast-slow dissection, taking the singular limit of the full system for fixed slow variable values; by Fenichel theory, their isochron geometries represent a quite accurate approximation to the phase response dynamics of the full system. %Hence knowledge of the configuration of isochrons in the fast subsystem is useful for analyzing phase response in the full system.

\subsubsection{Isochron calculation}
\label{sec:isochron_calc}

Analytical solutions for isochrons can be calculated for some relatively simple systems \cite{Winfree:2001, Tatsumi:1983, Nomura:1994}, and closed form approximations to isochrons locally (very close to the periodic orbit) may be obtained for systems near bifurcation (near homoclinic points, in particular) by considering the specific normal forms of the bifurcation \cite{Brown:2004}. In general, however, isochrons must be numerically approximated. Backwards integration is the method typically used to find points lying near a particular isochron \cite{Josic:2006}. Let $\theta_{x} = \theta(x), x \in \Gamma$ be the phase for which we wish to calculate the isochron. If we choose a point $x_{0}$ such that $\|x - x_{0}\| < \epsilon, \epsilon > 0$, and $\epsilon$ sufficiently small, then integrating backwards for time $\tau$, we obtain $\Phi_{-\tau}(x_{0})$, which lies on approximately the same isochron as $\Phi_{-\tau}(x)$. Hence by choosing relative phase $\theta' \in [0,1]$, $y \in \Gamma$ such that $\theta(y) = \theta_{x} + \theta'$, and $y_{0}$ such that $\|y - y_{0}\| < \epsilon$, then integrating backwards for time $\tau' = \theta'/T$, we obtain $\Phi_{-\tau'}(y_{0})$, which lies near the isochron for $\Phi_{-\tau'}(y) = x$.

With this method, computing the isochron curve for a given phase $\theta_{x}$ is a matter of choosing a set of relative phases and integrating backwards from the associated test points. Considerable care must be taken with the backwards integration and the test points must be chosen judiciously in order to obtain visually comprehensible curves with reasonable efficiency. Although the fast subsystem may not be particularly stiff in comparison to the full system, backwards integration in the fast subsystem is quite susceptible to failure due to numerical instabilities. This necessitates a good choice of integration routine and careful control of step sizes and integration tolerances. We use the RADAU5 stiff integrator \cite{Hairer:1991} with an initial step size of $10^{-5}-10^{-11}$, maximum step size of $10^{-4}$, and relative and absolute tolerances of $10^{-12}$. The hyperbolicity of the stable periodic orbit implies that trajectories are attracted to it exponentially in forwards time and disperse from it exponentially in backwards time. Therefore we begin with a set $S$ of $K$ logarithmically spaced relative phases in an interval close to $\theta_{x}$.  For each relative phase in $S$, we integrate backwards to find its associated point on the isochron for $\theta_{x}$; we repeat the process until we obtain $N$ points on the isochron. 

The strong stability of the periodic orbit introduces a difficulty: trajectories leave the periodic orbit extremely rapidly in backwards time. The effect, in some cases, is to introduce significant gaps between the points calculated on an isochron; this is typically more problematic for outer isochrons. Our calculation method attempts to fill in test points adaptively to reduce the size of such gaps, but our heuristic methods do not always succeed. Superior methods for planar isochron calculation have been developed recently \cite{Guillamon:2009} which would likely have superior performance on the problem considered here. However, though our method produces somewhat crude sketches of the isochrons, the portraits we compute suffice to understand the isochron geometry of the HR fast subsystem. 

%For a given relative phase $\theta'$ and $\y$ such that $\theta(\y) = \theta_{x} + \theta'$, we let $\y_{0} = \y + \gamma\w$, where $\gamma = \pm \epsilon$ and $\w$ is the unit vector normal to the vector field at $\y$ pointing outward from $\Gamma$. We refer to the segment of the isochron calculated for $\gamma$ positive (negative) and thus lying outside (inside) $\Gamma$ as the {\it outer} ({\it inner}) isochron for phase $\theta_{x}$.
%interval $[\theta_{\mathrm{low}}, \theta_{\mathrm{high}}]$, with $\theta_{\mathrm{low}}$ close to $\theta_{x}$. Since the calculations may not be stable for every relative phase in  $[\theta_{\mathrm{low}}, \theta_{\mathrm{high}}]$ (this is particularly problematic for larger $\theta_{\mathrm{high}}$), we find the largest relative phase $\theta_{\mathrm{max}} \in S$ for which backwards integration succeeds. We then set $\theta_{\mathrm{high}} = \theta_{\mathrm{max}}$, and recalculate $S$ to have $K$ relative phases in the new interval, if necessary. For each relative phase in $S$, we integrate backwards to find its associated point on the isochron for $\theta_{x}$; we repeat the process until we obtain $N$ points on the isochron. We summarize our isochron calculation method in Algorithm \ref{alg:isochron}. 

%\input{\thischapter/Algorithms/GrowIsochron}

\subsubsection{Early active segment isochron portrait}
\label{subsec:early_active}

\begin{figure}[!ht]
 %\begin{minipage}{3.25in}
     \begin{center}
	\includegraphics[width= 3.5in]{\figpath/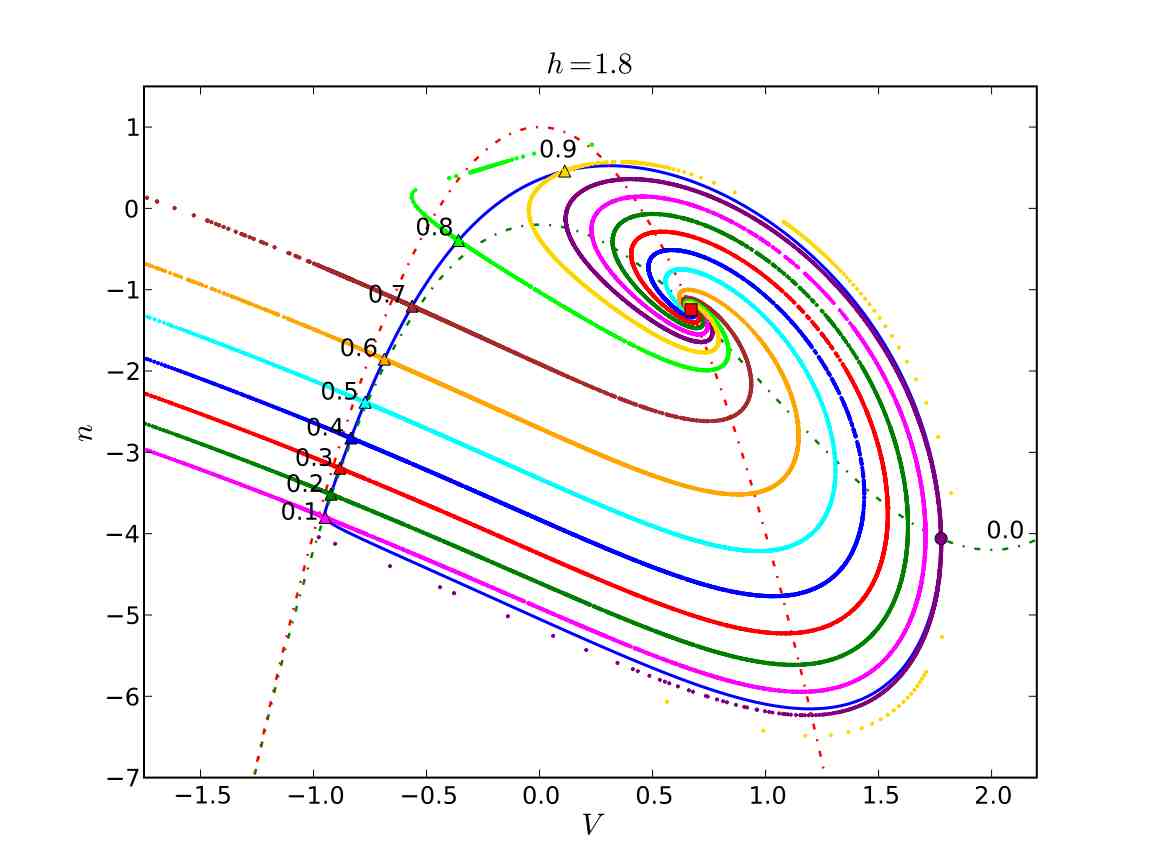}
	%\
	    \end{center}
  %\end{minipage}
\caption[Early Active Segment Isochron Portrait]{Phase portrait and isochrons in the HR fast subsystem near the beginning of the active segment of the full system, $h = 1.8$. }
\label{fig:PP_Iso_18}
\end{figure}

\begin{figure}[!ht]
      \begin{center}
	  \includegraphics[width= 3.5in]{\figpath/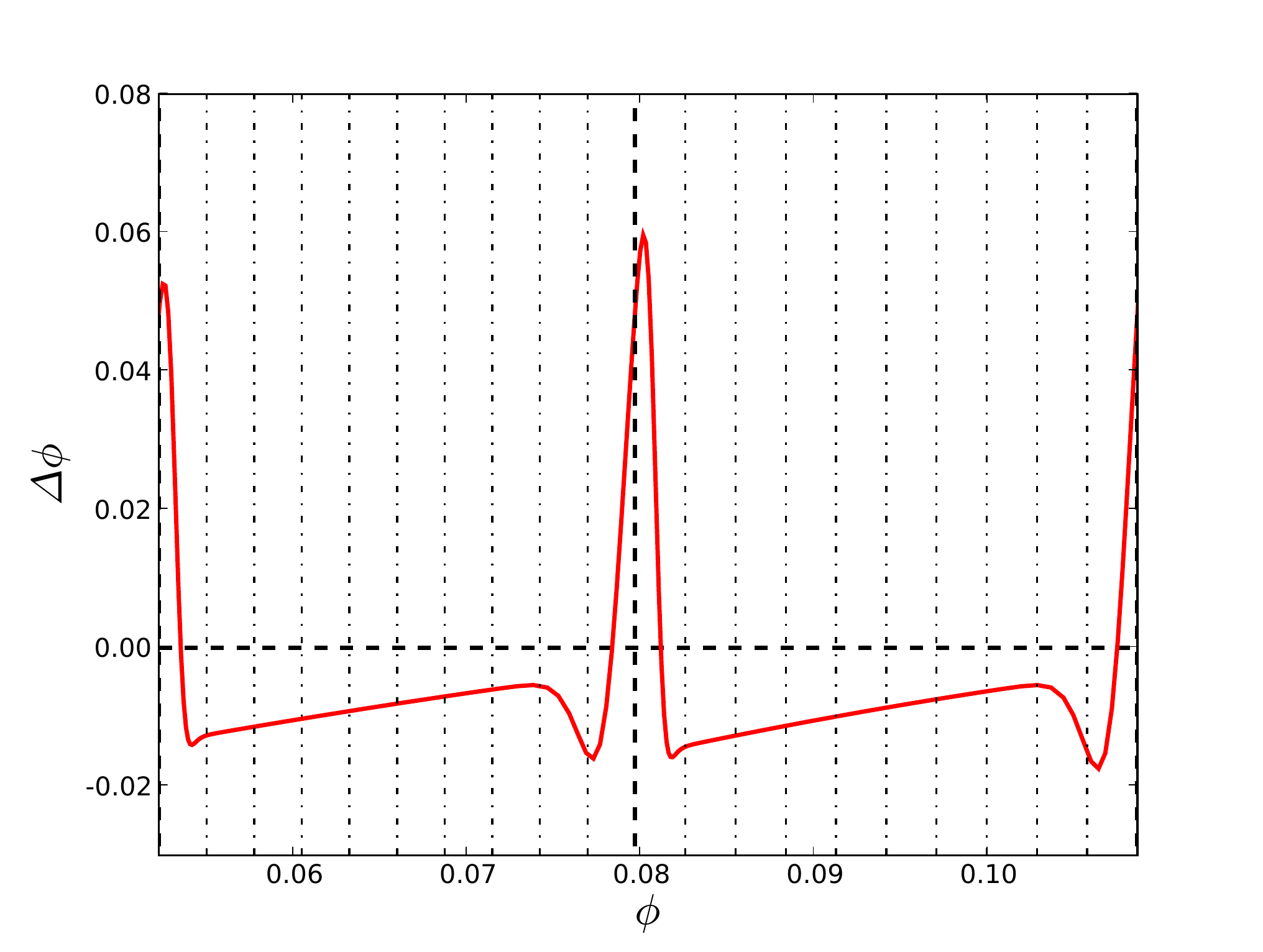}
    \end{center}
\caption[Early Active Segment Infinitesimal Burst Phase Response ]{Close-up of infinitesimal BPRC in the full HR model near the beginning of the active segment. The heavy dashed line marks the spike peak for the spike closest to $h=1.8$ slow variable value. Light dash-dotted lines denote equally spaced interspike phases in the full system.}
\label{fig:Iso_BPRC_18}
\end{figure}

Figure \ref{fig:PP_Iso_18} shows the fast subsystem phase portrait with isochrons at $h=1.8$, near the beginning of the active segment of the full burst cycle and the initial portion of Segment I of the burst phase response. For comparison, Figure \ref{fig:Iso_BPRC_18} shows the infinitesimal burst phase response in the full system in the interval corresponding to $h \in [1.78464559529, 1.85251757552]$, centered about the second spike of the burst. The dashed-dotted lines of Figure \ref{fig:Iso_BPRC_18} (and subsequent figures) mark equally spaced phases between spikes of the burst to provide a visual reference for comparing phase response in the full HR system with the phase geometry of the fast subsystem.

The outer isochrons for phases in the hyperpolarized portion of the spike, 0.1 to 0.7, extend linearly from the left edge of the periodic orbit and are nearly parallel to one another.%\footnote{In each of the fast subsystems we consider, the outer isochrons for each extend out nearly linearly to $V$ values below -200, but we focus on the area of the phase plane close to the periodic orbit.} 
The inner isochrons at the same phases also have extended parallel linear segments, but past the right branch of the $n$ nullcline, they begin to curve and spiral into the unstable fixed point, coming close to the inside depolarized edge of the orbit. The phases of the upswing of the spike, 0.8 to 0.0, have outer isochrons that bend to the right and curve around the depolarized side of the orbit; to the left of the right branch of the $n$ nullcline, they straighten and extend parallel to the 0.1 to 0.7 isochrons. The inner isochrons of the spike upswing also spiral into the unstable fixed point, but they remain away from the depolarized edge of the periodic orbit. The inner 0.0 isochron, corresponding to the spike peak, lies very close to the depolarized edge of the periodic orbit above the $V$ nullcline, and the outer isochron hugs the outer edge of the periodic orbit tightly. The isochrons for phases 0.9 to 0.08 are nearly tangent to the periodic orbit where they intersect it.

%%%%%%%%%%%%%%%%%%
%	Isochrons 1.8
%%%%%%%%%%%%%%%%%%  

%%%%%%%%%%%%%%%%%%
%	BPRC close-up 1.8
%%%%%%%%%%%%%%%%%%  

%%%%%%%%%%%%%%%%%%
%	Upper Isochrons 1.8
%%%%%%%%%%%%%%%%%%  
\begin{figure}[!ht]
 \begin{minipage}{3in}
     \begin{center}
	\includegraphics[width= 3in]{\figpath/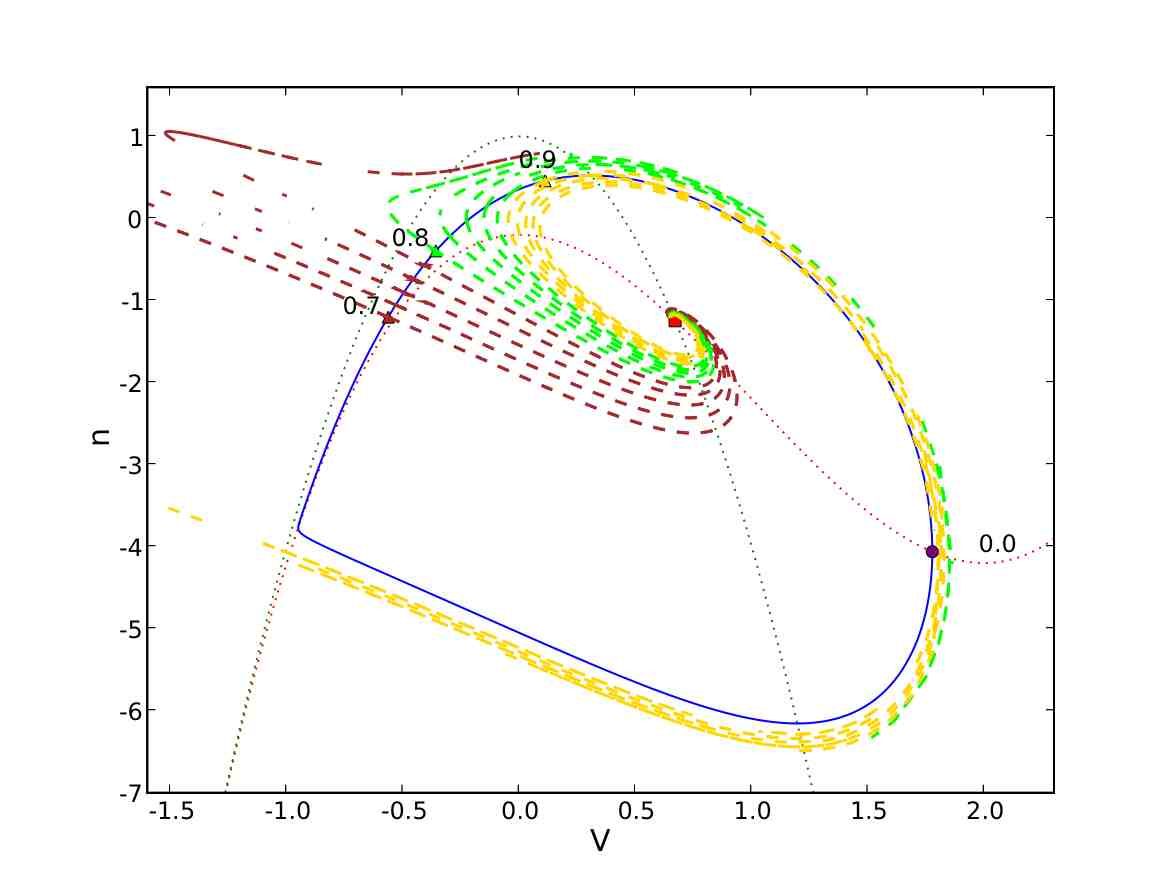}
	(a)
	    \end{center}
  \end{minipage}
  \begin{minipage}{3in}
      \begin{center}
	  \includegraphics[width= 3in]{\figpath/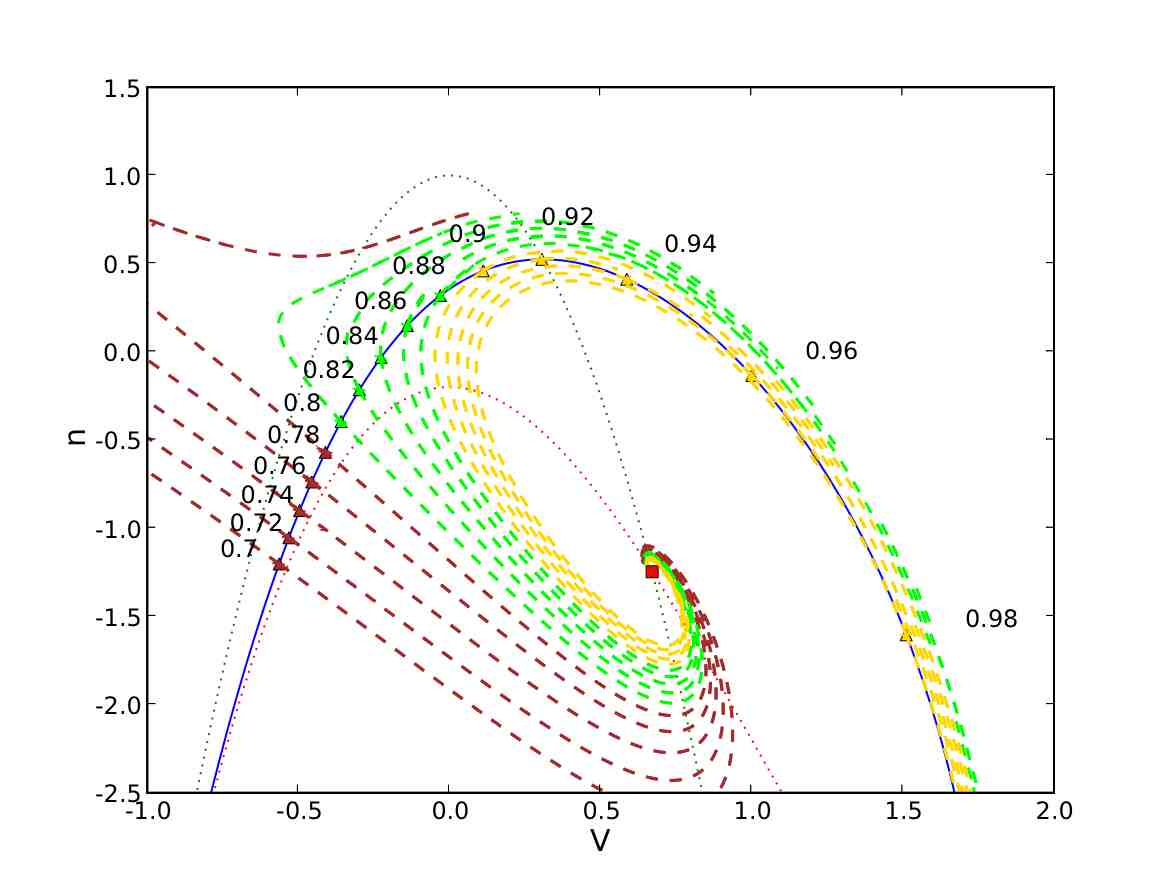}
	  (b)
    \end{center}
  \end{minipage}%
\caption[Early Active Segment Isochron Portrait (Late Phases)]{Phase portrait of the HR fast subsystem and isochrons for phases in $[0.7, 0.0)$, $h=1.8$. Isochrons with phases in $[0.7, 0.8)$ are colored brown; $[0.8, 0.9)$, green; $[0.9, 0.0)$, yellow. (a) Full orbit. (b) Close-up of `late phase' region.}
\label{fig:PP_Iso_18upper}
\end{figure}

%%%%%%%%%%%%%%%%%%
%	Lower Isochrons 1.8
%%%%%%%%%%%%%%%%%%  
\begin{figure}[!ht]
 \begin{minipage}{3in}
     \begin{center}
	\includegraphics[width= 3in]{\figpath/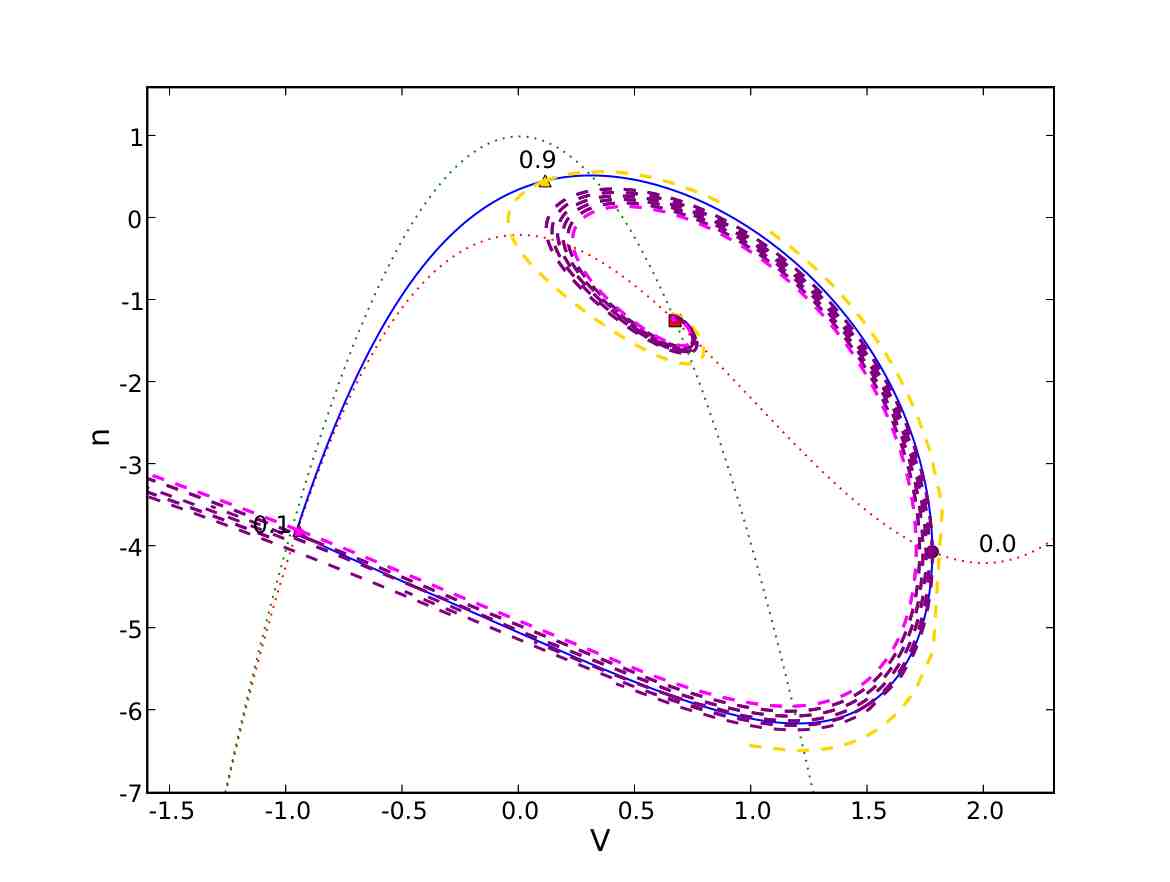}
	(a)
	    \end{center}
  \end{minipage}
  \begin{minipage}{3in}
      \begin{center}
	  \includegraphics[width= 3in]{\figpath/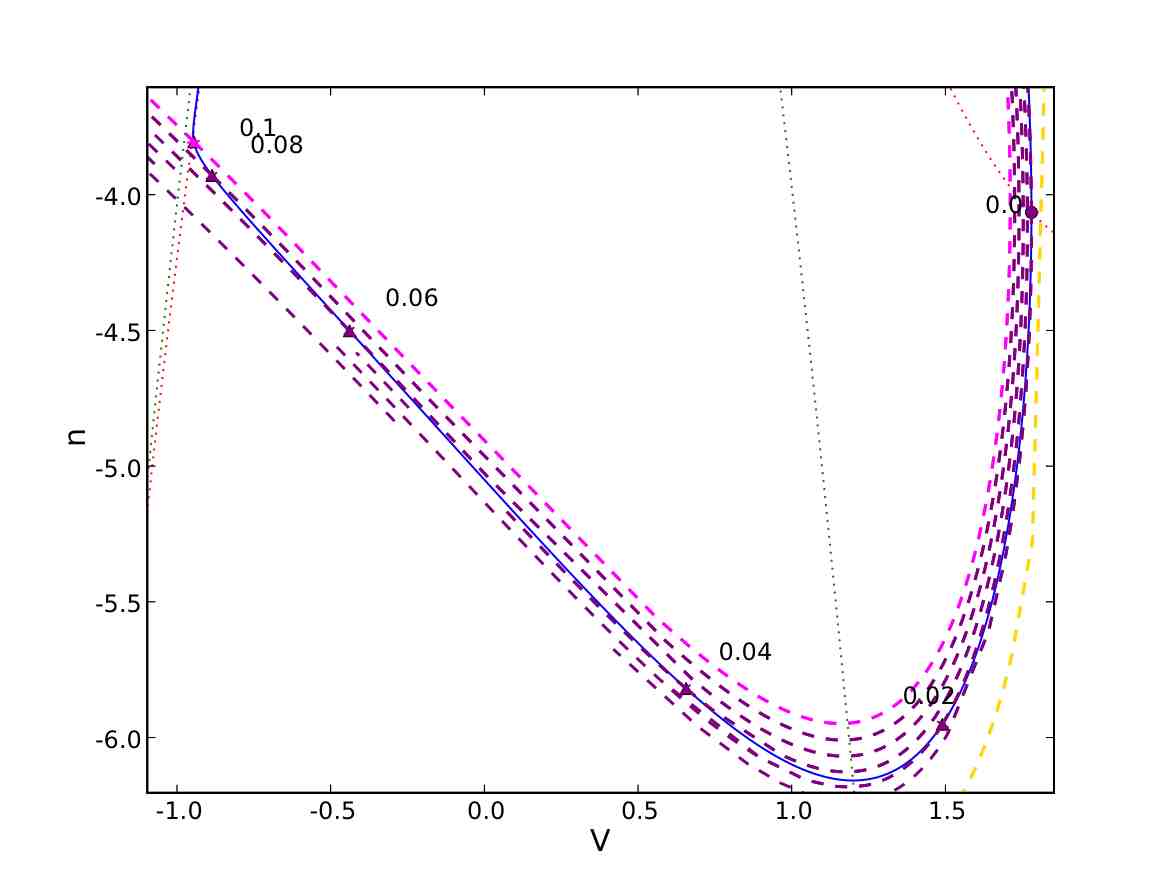}
	  (b)
    \end{center}
  \end{minipage}
\caption[Early Active Segment Isochron Portrait (Early Phases)]{Phase portrait of the HR fast subsystem and isochrons for phases in $[0.9, 0.1]$, $h=1.8$. Isochrons with phases in $[0.0, 0.1)$ are colored purple; the 0.9 isochron, yellow; the 0.1 isochron, pink. (a) Full orbit. (b) Close-up of early phase region.}
\label{fig:PP_Iso_18lower}
\end{figure}

The effect of a small excitatory perturbation, such as approximated in the infinitesimal BPRC method, is to move the trajectory away from the limit cycle in the positive $V$ direction. Imagine the action of such a perturbation as a short horizontal line segment in the phase plane with its left end point on the limit cycle. The right end point (`perturbation point')  marks the position of the trajectory after perturbation. As the left end point moves around the limit cycle, the intersection of the perturbation point with the isochrons forms a curve in the phase plane. Following this curve, we can trace the phase response over the course of a spike; it is a fast subsystem approximation to the burst phase response curve in the full system at corresponding $h$ values. 

When the left end point phase is between 0.1 and 0.82, the perturbation point intersects the straight portions of inner isochrons for larger phases, and so the phase is advanced (negative $\Delta\theta$). The inner isochrons are spaced further apart at larger phases in this range, so the magnitude of phase response for the same size perturbation decreases as the phase increases. Thus the infinitesimal PRC segment of Figure \ref{fig:Iso_BPRC_18} between $\theta \approx 0.055$ and $\phi \approx 0.075$ is negative and increasing. 

As shown in Figure \ref{fig:PP_Iso_18upper}, between 0.82 and 0.92, the perturbation point intersects curved portions of isochrons as they spiral into the unstable fixed point. In this region of phase space, the isochrons bunch together, so the magnitude of the phase advancement increases, hence we see the dip near $\phi \approx 0.075$ in Figure \ref{fig:Iso_BPRC_18}.

After the maximum $n$ value of the periodic orbit, starting just past phase 0.92, the perturbation point lies outside the orbit and intersects isochrons of lower phases as they wrap around the periodic orbit's depolarized edge. This configuration persists through the majority of the depolarized portion of the spike, until approximately phase 0.03 (see Figures \ref{fig:PP_Iso_18upper} (b) and \ref{fig:PP_Iso_18lower} (b)). Around phase 0.0, the isochrons lie nearly tangent to the periodic orbit where they intersect it, and significant portions of their inner and outer segments (away from the point of intersection) remain very close to the periodic orbit. Hence this is a region of phase delay, and there is a sharp peak in the magnitude of phase delay very close to phase 0.0, the spike maximum (\cf Figure \ref{fig:Iso_BPRC_18} near $\theta \approx 0.08$). 

Past phase 0.03, the perturbation point lies inside the periodic orbit again, so that perturbation advances the phase. Since the inner (and outer) isochrons between phase 0.0 and 0.1 lie very close the periodic orbit, the magnitude of the phase change in this region is relatively large, and the change from phase 0.0 is steep. The cycle of phase response repeats beyond phase 0.1.

\subsubsection{Late active segment isochron portrait}
\label{subsec:late_active}

\begin{figure}[!ht]
     \begin{center}
	\includegraphics[width= 3.5in]{\figpath/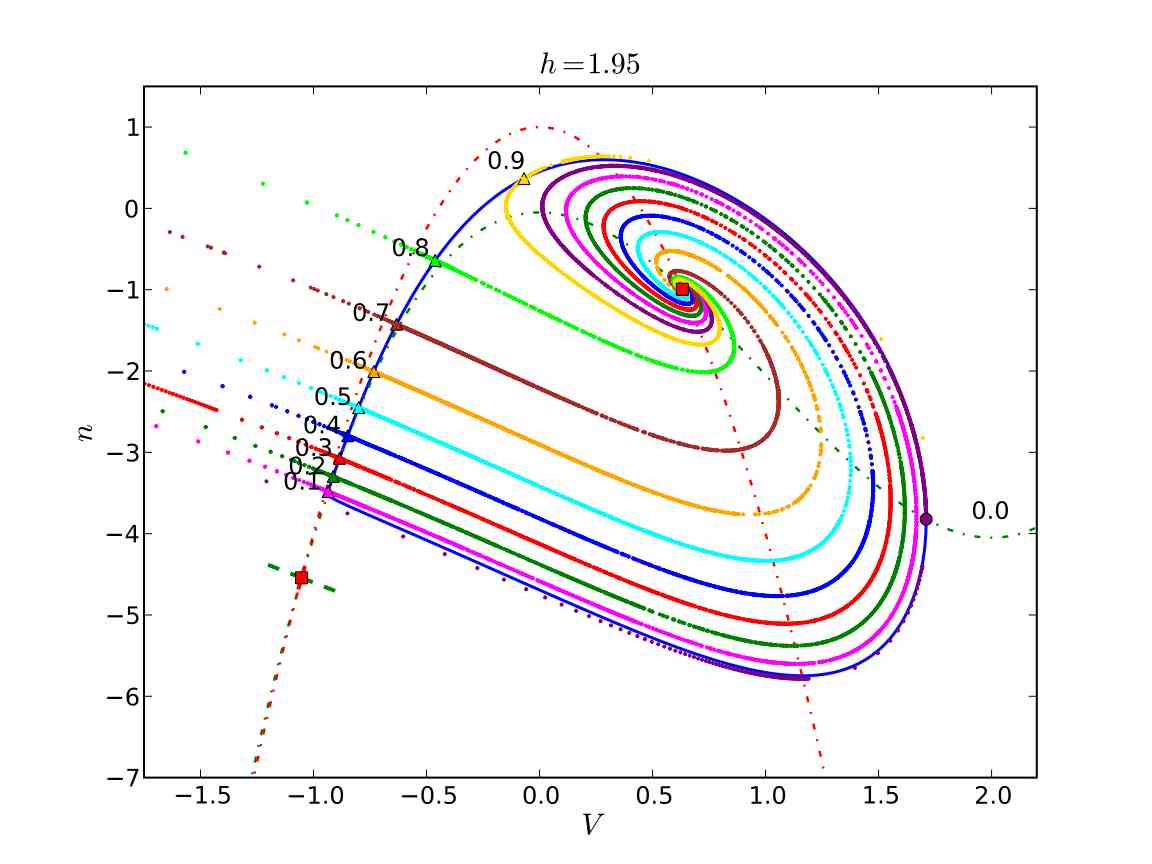}
	    \end{center}
\caption[Late Active Segment Isochron Portrait]{Phase portrait and isochrons in the HR fast subsystem latter half of the active segment of the full system, $h = 1.95$.}
\label{fig:PP_Iso_195}
\end{figure}

\begin{figure}[!ht]
      \begin{center}
	  \includegraphics[width= 3.5in]{\figpath/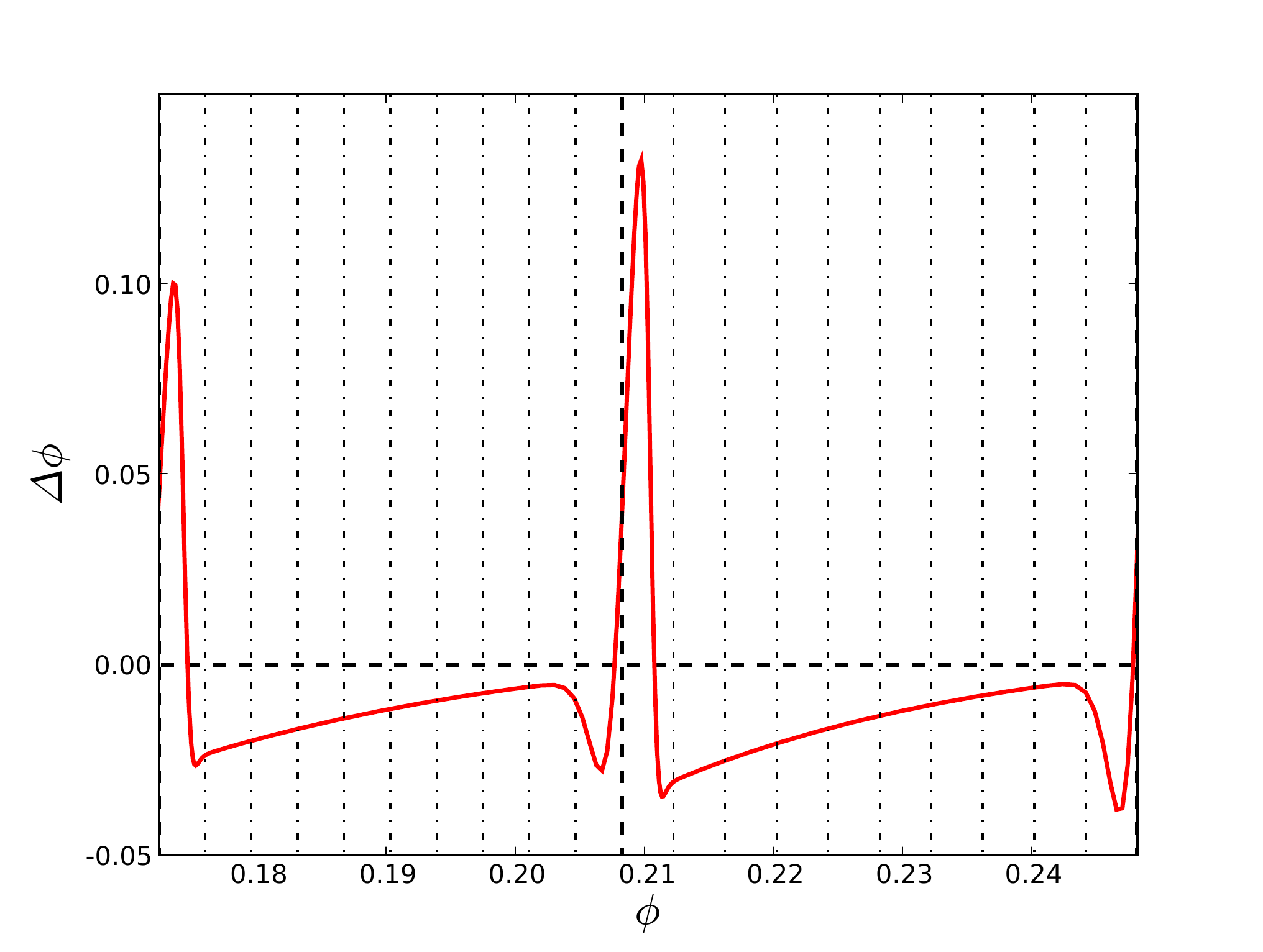}
    \end{center}
\caption[Late Active Segment Infinitesimal Burst Phase Response ]{Close-up of infinitesimal BPRC in the full HR model in the latter half of the active segment. The heavy dashed line marks the spike peak for spike closest to $h=1.95$ slow variable value. Light dash-dotted lines denote equally spaced interspike phases in the full system. Note that the peak of phase response is greater in magnitude than the peak shown in Figure \ref{fig:Iso_BPRC_18}, and it falls after the phase of the proximal spike peak in the full HR system.}
\label{fig:Iso_BPRC_195}
\end{figure}

\begin{figure}[!ht]
 \begin{minipage}{3in}
     \begin{center}
	\includegraphics[width= 3in]{\figpath/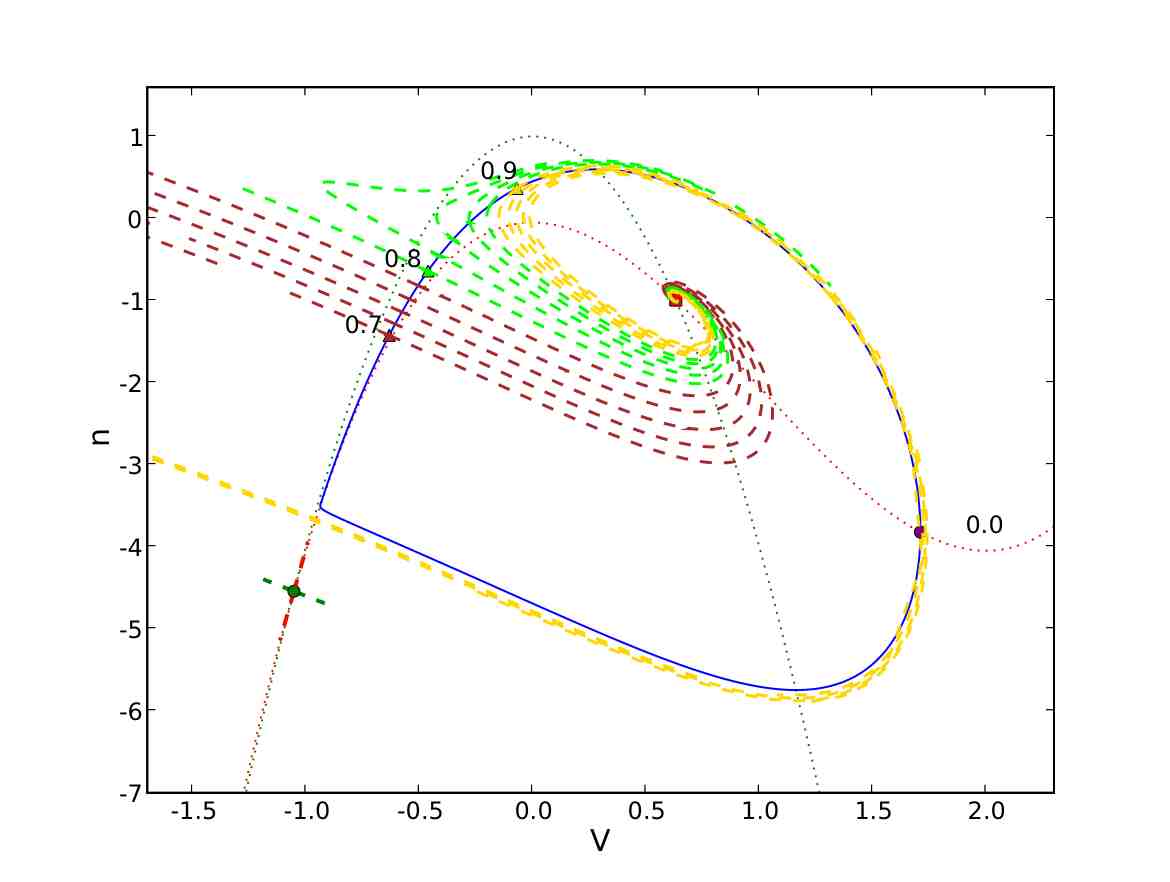}
	(a)
	    \end{center}
  \end{minipage}
  \begin{minipage}{3in}
      \begin{center}
	  \includegraphics[width= 3in]{\figpath/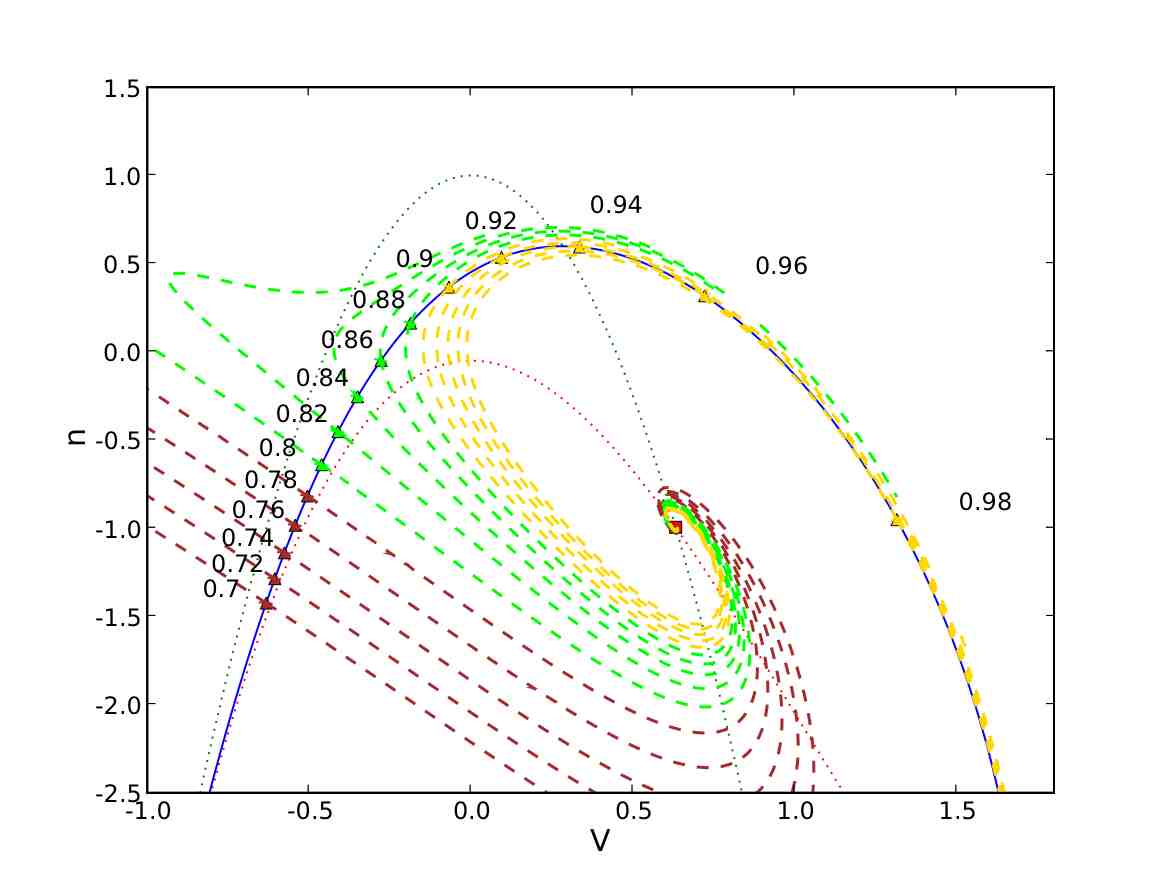}
	  (b)
    \end{center}
  \end{minipage}%

\caption[Late Active Segment Isochron Portrait (Late Phases)]{Phase portrait of the HR fast subsystem and isochrons for phases in $[0.7, 0.0)$, $h=1.95$. Isochrons with phases in $[0.7, 0.8)$ are colored brown; $[0.8, 0.9)$, green; $[0.9, 0.0)$, yellow. (a) Full orbit. (b) Close-up of late phase region.}
\label{fig:PP_Iso_195upper}
\end{figure}

\begin{figure}[!ht]
 \begin{minipage}{3in}
     \begin{center}
	\includegraphics[width= 3in]{\figpath/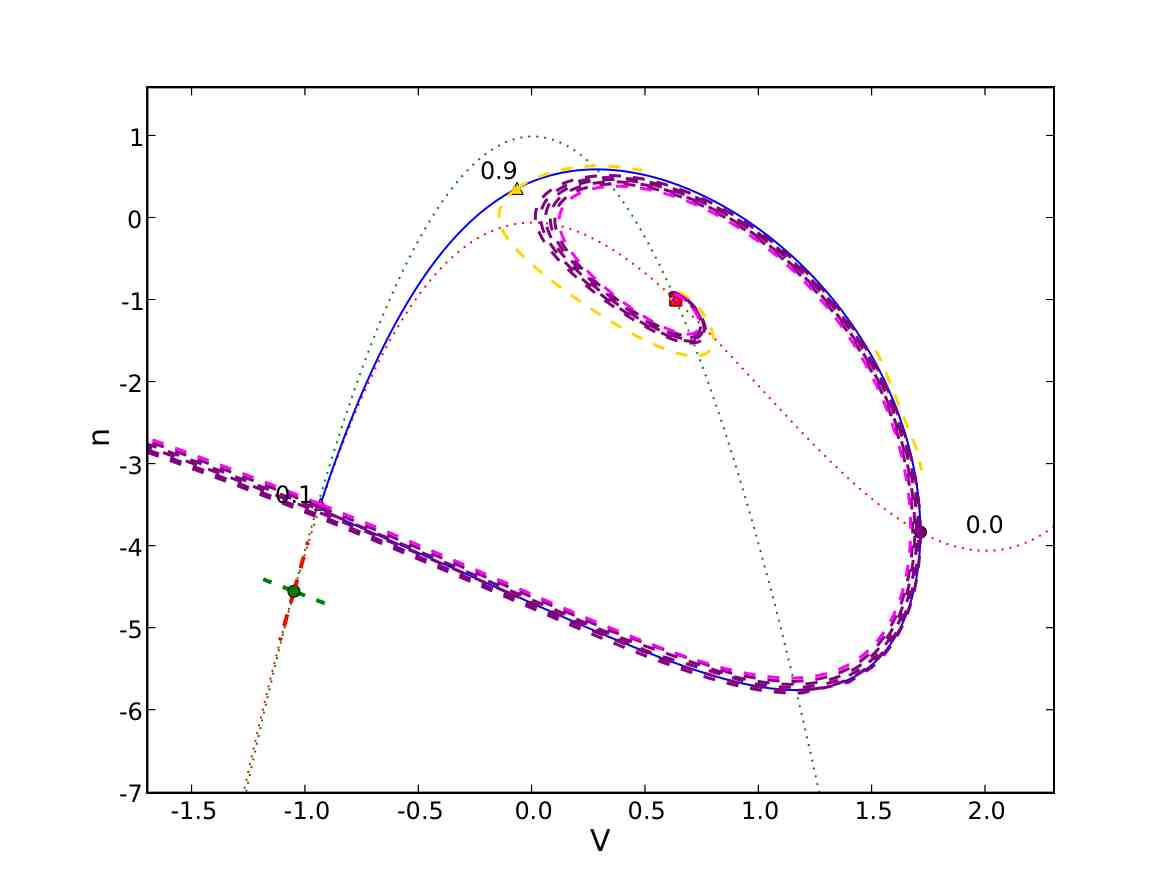}
	(a)
	    \end{center}
  \end{minipage}
  \begin{minipage}{3in}
      \begin{center}
	  \includegraphics[width= 3in]{\figpath/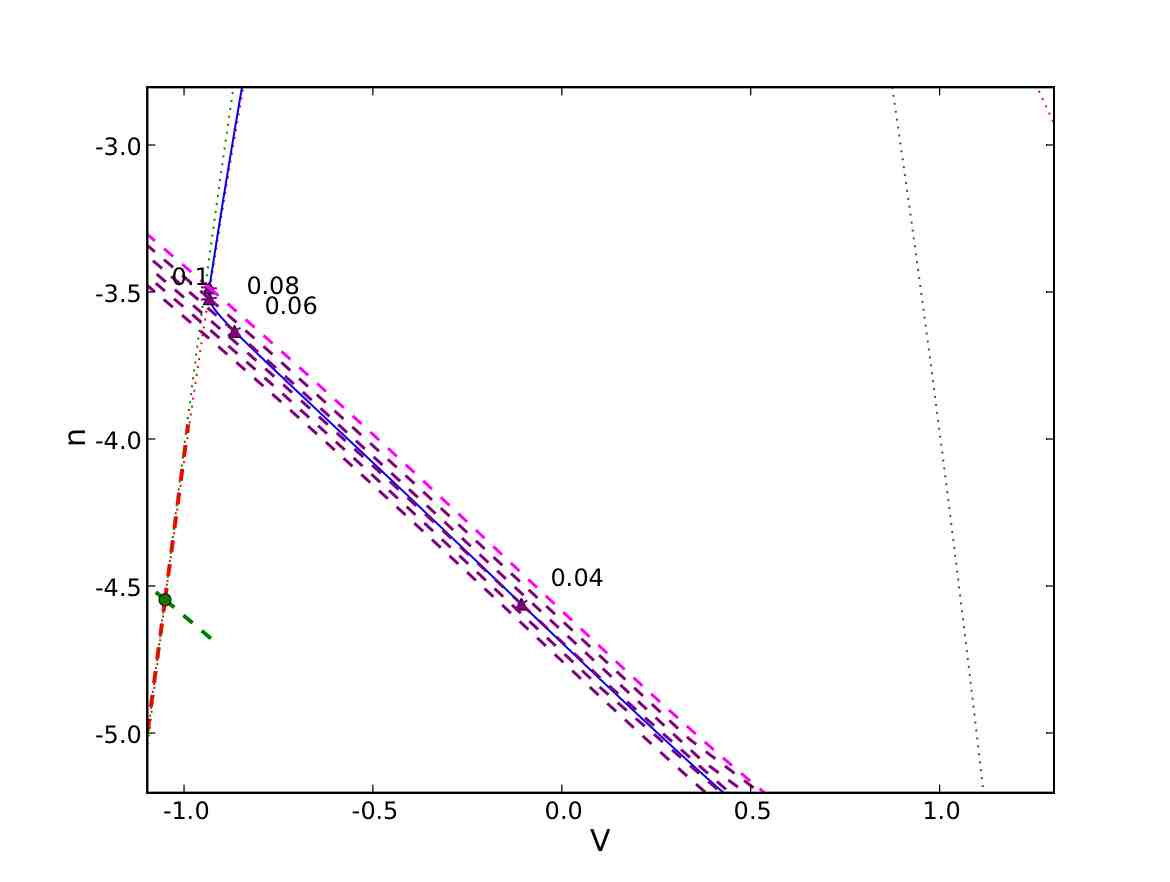}
	  (b)
    \end{center}
  \end{minipage}%
\caption[Late Active Segment Isochron Portrait (Early Phases)]{Phase portrait of the HR fast subsystem and isochrons for phases in $[0.9, 0.1]$, $h=1.95$. Isochrons with phases in $[0.0, 0.1)$ are colored purple; the 0.9 isochron, yellow; the 0.1 isochron, pink. (a) Full orbit. (b) Close-up of early phase region.}
\label{fig:PP_Iso_195lower}
\end{figure}

The configurations of the isochrons in the middle and later portions of the active segment largely resemble that of the early portion, and the corresponding regions of the infinitesimal  BPRC also have similar shapes. Figure \ref{fig:PP_Iso_195}, the phase portrait for $h = 1.95$, shows the same basic spiral pattern of inner isochrons as Figure \ref{fig:PP_Iso_18}. The section of the infinitesimal BPRC corresponding to $h \in [1.92341235175,1.99914946946]$ (including the sixth spike of the full system), shown in Figure \ref{fig:Iso_BPRC_195}, also has a shape very similar to the one seen in Figure \ref{fig:Iso_BPRC_18} for $h=1.8$. The isochrons at phases in the upswing and downswing of the spike depolarization (see Figures \ref{fig:PP_Iso_195upper} and \ref{fig:PP_Iso_195lower}) are arranged quite similarly to their counterparts in the fast subsystem at $h=1.8$. 

The most significant change in the isochron geometry between $h=1.8$ and $h=1.95$ is that the isochrons lie closer together, particularly near the depolarized edge of the periodic orbit. This compression is visible in Figure \ref{fig:PP_Iso_195} for the isochrons at phases 0.1 to 0.6. The closeness does not simply stem from the tighter arrangement of these phase points along the hyperpolarized edge of the periodic orbit, as can be seen by comparing it with Figure \ref{fig:PP_Iso_18}. 

%%%%%%%%%%%%%%%%%%
%	Isochrons 1.95
%%%%%%%%%%%%%%%%%%  

%%%%%%%%%%%%%%%%%%
%	PRC close-up 1.95
%%%%%%%%%%%%%%%%%%  

%\clearpage

%%%%%%%%%%%%%%%%%%
%	Upper Isochrons 1.95
%%%%%%%%%%%%%%%%%%  

%%%%%%%%%%%%%%%%%%
%	Lower Isochrons 1.95
%%%%%%%%%%%%%%%%%%  

A similar pattern of compression can be seen for the isochrons at phases between 0.8 and 0.1 in Figures \ref{fig:PP_Iso_195upper} and \ref{fig:PP_Iso_195lower}. To the right of the $n$ nullcline, the inner isochrons hug the inside of the periodic orbit. The outer isochrons remain close to the outer rim of the periodic orbit, passing between the saddle point and the periodic orbit and extending linearly beyond the left branches of the nullclines. The isochrons of the fast subsystem exist only within the basin of attraction of the periodic orbit; the emergence of the saddle point in the fast subsystem introduces a nearby boundary (the stable manifold of the saddle point) to the region in which the isochrons are confined. The saddle point moves closer to the periodic orbit as $h$ increases and the homoclinic bifurcation approaches. This forces the isochrons closer together in order to squeeze through the narrowing of the gap between the saddle and the limit cycle. 

Since the arrangement of the isochrons is roughly the same for fast subsystem at $h=1.95$ and $h=1.8$, the full system phase response curves have very similar forms, but with two notable differences. First, the magnitude of the phase response is larger for $h=1.95$. This follows from the closer proximity of isochrons of a given phase to the periodic orbit. If we consider the horizontal line segment representing a voltage perturbation, as above, then for a fixed line length (perturbation strength), the perturbation point for a perturbation at a given phase will lie at nearby locations in phase space for $h=1.8$ and $h=1.95$. However, since the isochrons for $h=1.95$ lie closer to the periodic orbit, the perturbation point for $h=1.95$ lies on an isochron with a greater phase difference (further away in phase) than does the perturbation point for $h=1.8$. This difference is reflected in the greater magnitude of phase response in the full subsystem, as seen in Figure \ref{fig:Iso_BPRC_195}. This pattern of isochron compression progresses as $h$ increases, so that the isochrons lie increasingly close together over the course of the active segment, and the magnitude of phase response in the full system grows accordingly over the course of Segment I.

The second difference between the linear PRCs is that the peak of phase response is shifted further past the phase of the spike peak. The maximum $n$ value of the periodic orbit now lies near 0.94, rather than 0.92 for $h=1.8$, so that the zero crossing of phase response in the full system is shifted closer to the  spike maximum. The point of maximum proximity of the outer isochrons to the periodic orbit is shifted further past phase 0.0 as well. In addition, the $h$ value of the spike peak pictured in Figure \ref{fig:Iso_BPRC_195} is 1.96048585711, so that the exactly corresponding fast subsystem has even greater isochron compression than the one drawn for $h=1.95$. These differences in the isochron configuration and $h$ value produce a phase response peak that lies further to the right of the spike peak. 

\subsubsection{Isochron geometry at the end of active spiking}
\label{subsec:active_term}

\begin{figure}[!ht]
     \begin{center}
	\includegraphics[width= 3.5in]{\figpath/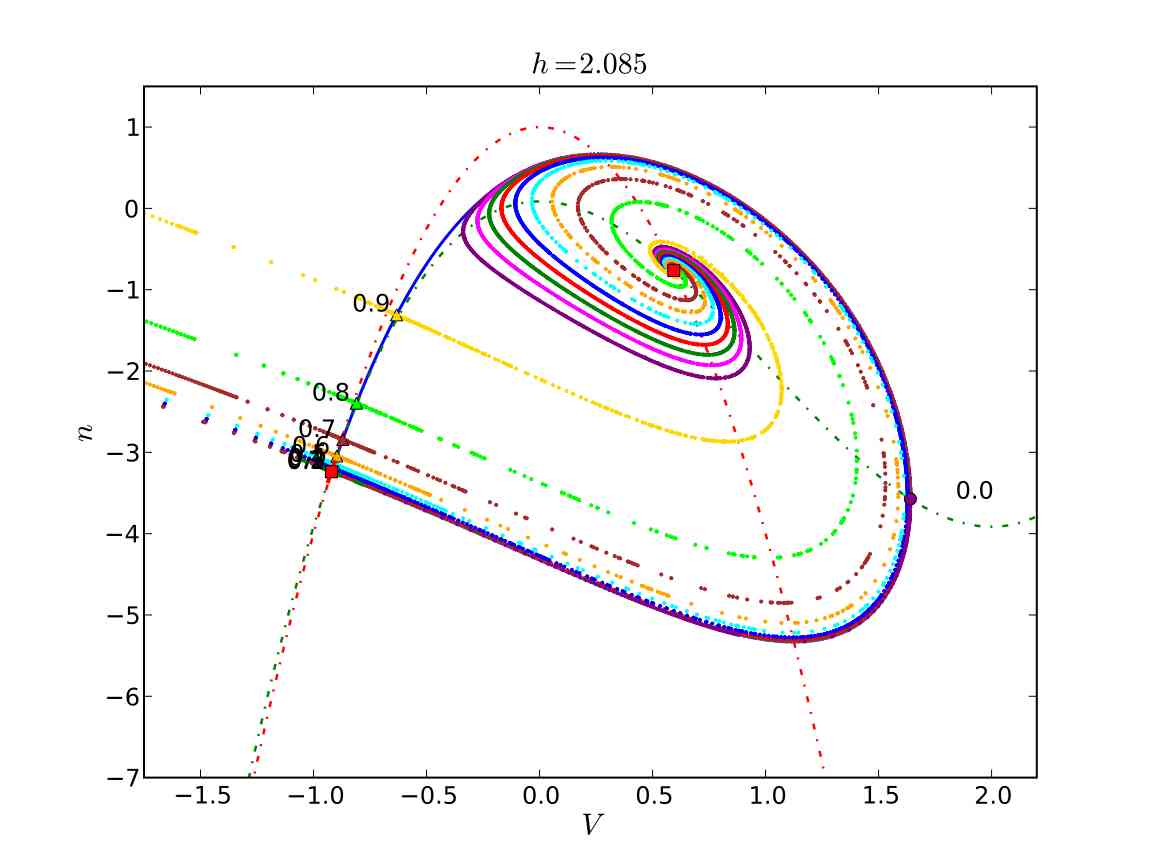}
	    \end{center}
\caption[Active Segment Termination Isochron Portrait]{Phase portrait and isochrons in the HR fast subsystem near the end of the active segment of the full system, $h = 2.085$.}
\label{fig:PP_Iso_2085}
\end{figure}

\begin{figure}[!ht]
 \begin{minipage}{3in}
     \begin{center}
	\includegraphics[width= 3in]{\figpath/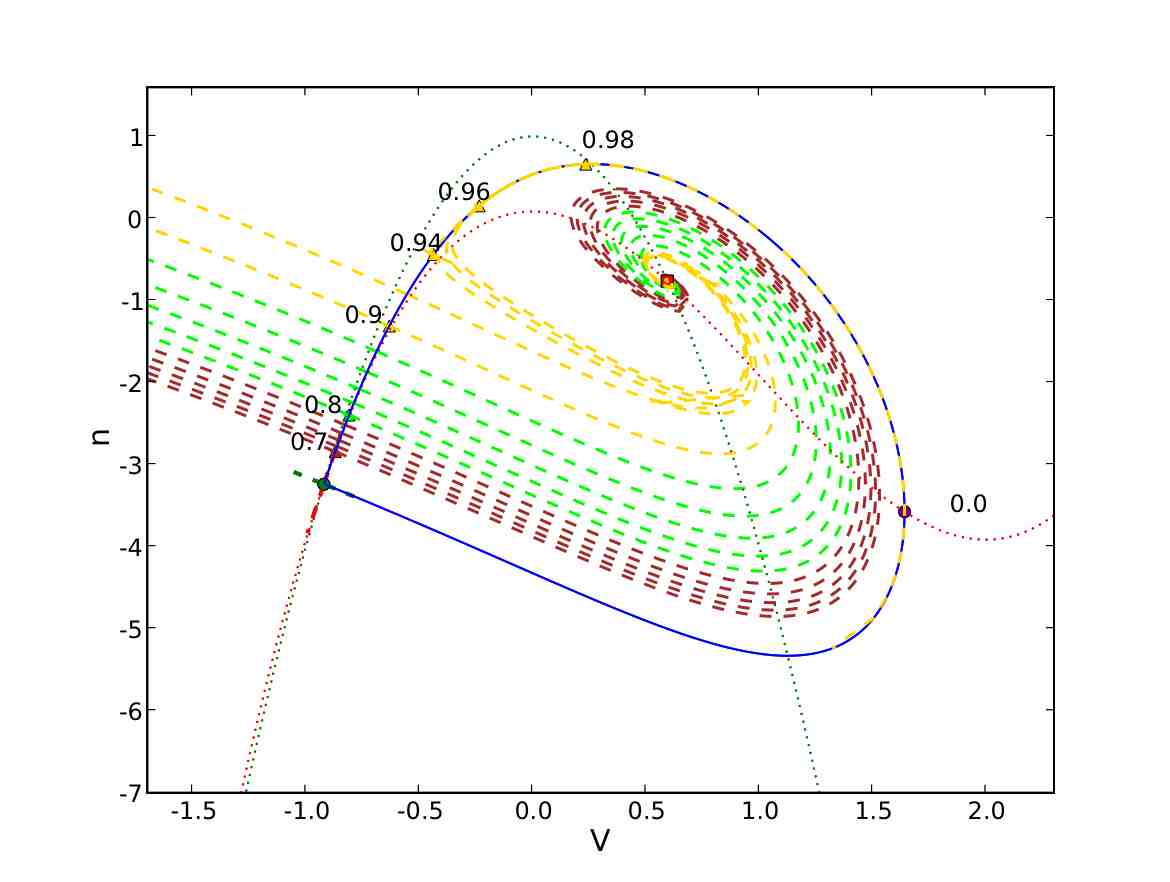}
	(a)
	    \end{center}
  \end{minipage}
  \begin{minipage}{3in}
      \begin{center}
	  \includegraphics[width= 3in]{\figpath/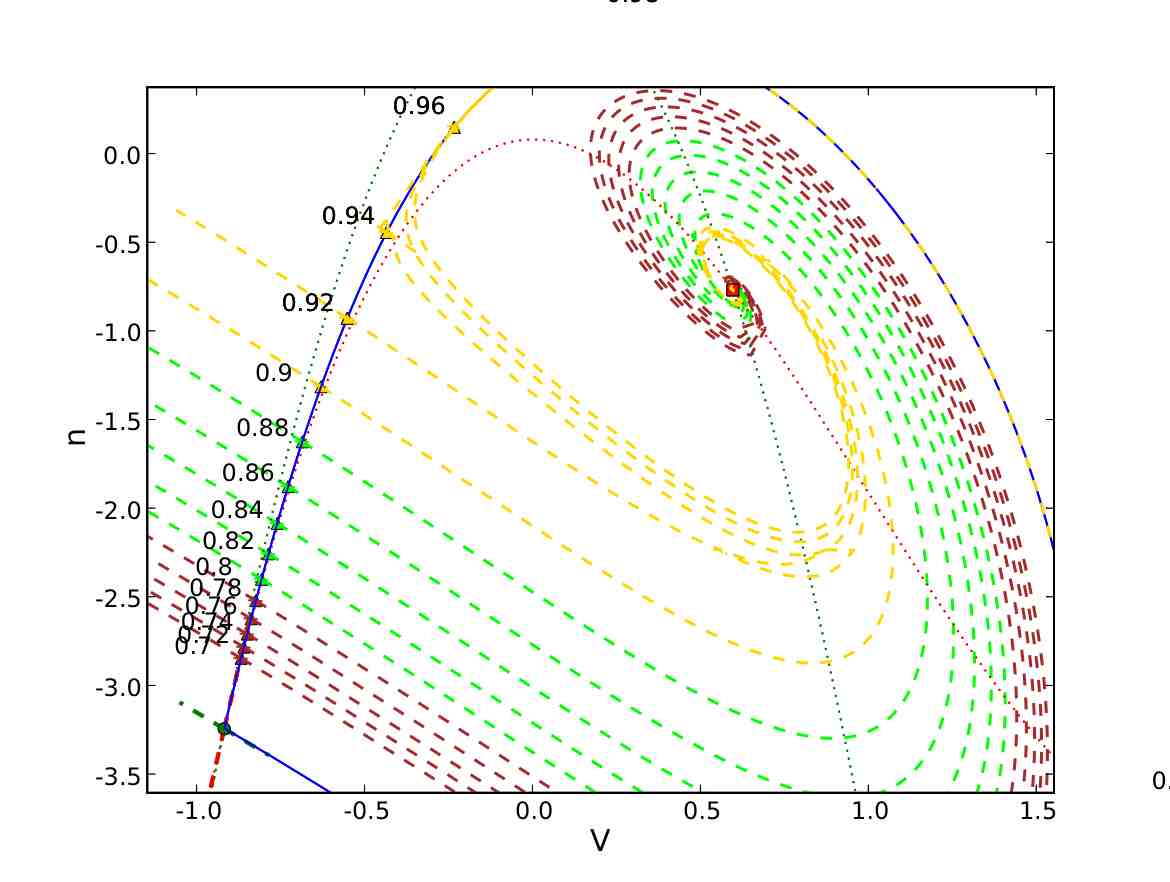}
	  (b)
    \end{center}
  \end{minipage}%
\caption[Active Segment Termination Isochron Portrait (Late Phases)]{Phase portrait of the HR fast subsystem and isochrons for phases in $[0.7, 0.0)$, $h=2.085$. Isochrons with phases in $[0.7, 0.8)$ are colored brown; $[0.8, 0.9)$, green; $[0.9, 0.0)$, yellow. (a) Full orbit. (b) Close-up of late phase region.}
\label{fig:PP_Iso_2085upper}
\end{figure}

\begin{figure}[!ht]
 \begin{minipage}{3in}
     \begin{center}
	\includegraphics[width= 3in]{\figpath/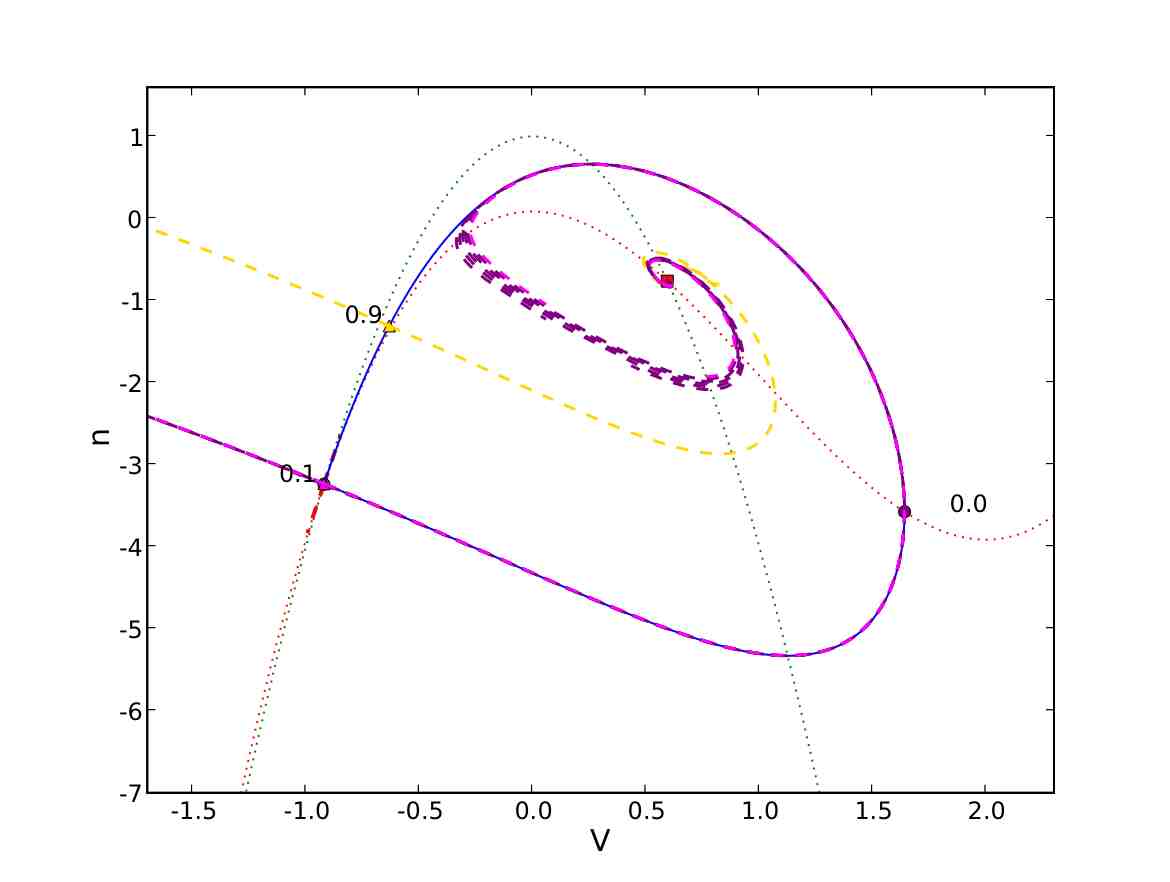}
	(a)
	    \end{center}
  \end{minipage}
  \begin{minipage}{3in}
      \begin{center}
	  \includegraphics[width= 3in]{\figpath/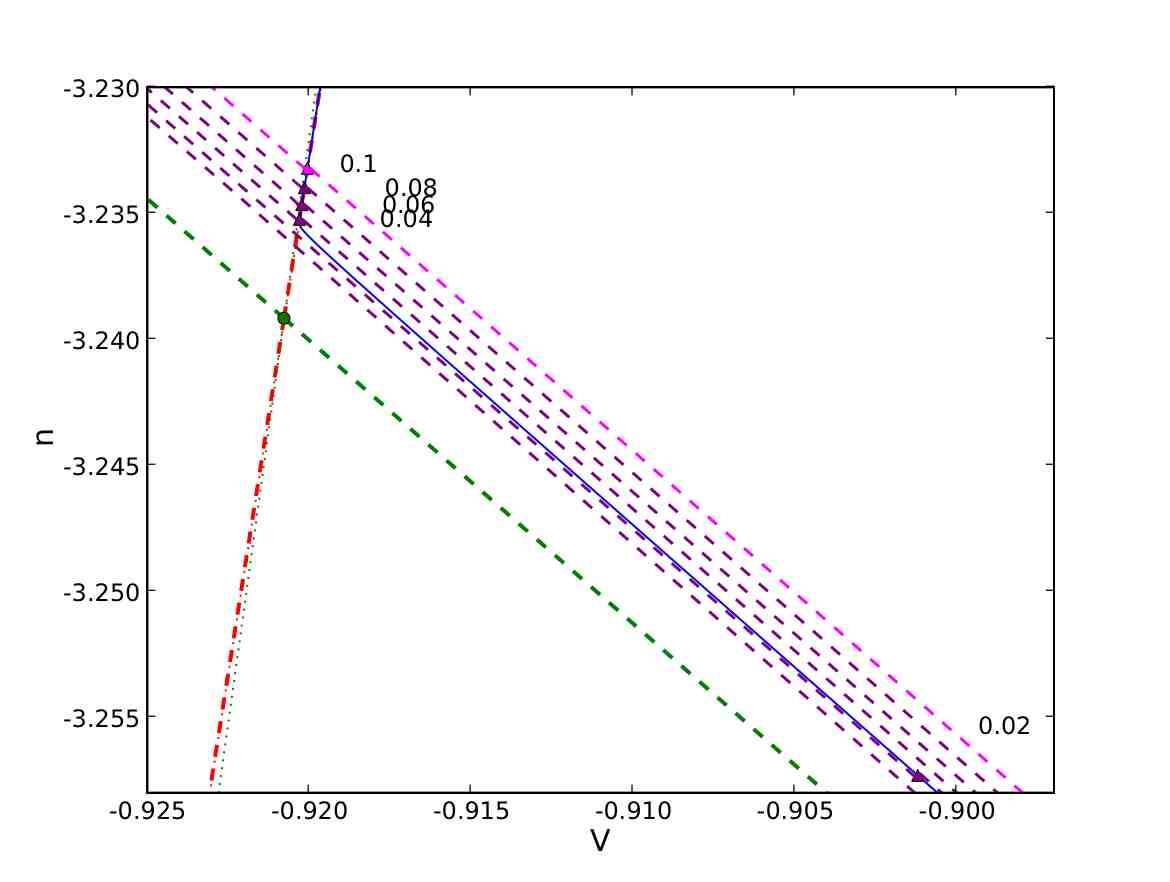}
	  (b)
    \end{center}
  \end{minipage}%
\caption[Active Segment Termination Isochron Portrait (Early Phases)]{Phase portrait of the HR fast subsystem and isochrons for phases in $[0.9, 0.1]$, $h=2.085$. Isochrons with phases in $[0.0, 0.1)$ are colored purple; the 0.9 isochron, yellow; the 0.1 isochron, pink. (a) Full orbit. (b) Close-up of early phase region near saddle point.}
\label{fig:PP_Iso_2085lower}
\end{figure}

The active segment of the HR burst cycle ends in a homoclinic bifurcation in which the saddle point and the periodic orbit merge. This seals the gap through which the outer isochrons for phases near 0.0 must pass, and it increases the period of the orbit towards infinity. As the homoclinic point approaches, the inner isochrons press closer to the inner rim of the periodic orbit, and outer isochrons along the depolarized portion of the periodic orbit are forced closer to the orbit's outer edge. Outer isochrons in the hyperpolarized part of the orbit straighten, extending further and more linearly in the $-V$ direction. These changes can be seen in Figures \ref{fig:PP_Iso_2085}, \ref{fig:PP_Iso_2085upper}, and \ref{fig:PP_Iso_2085lower}, which show the isochrons and phase portraits for the fast subsystem at $h=2.085$, very close to the bifurcation point at $h \approx  2.08560088198$.

\begin{figure}[!ht]
      \begin{center}
	  \includegraphics[width= 3.5in]{\figpath/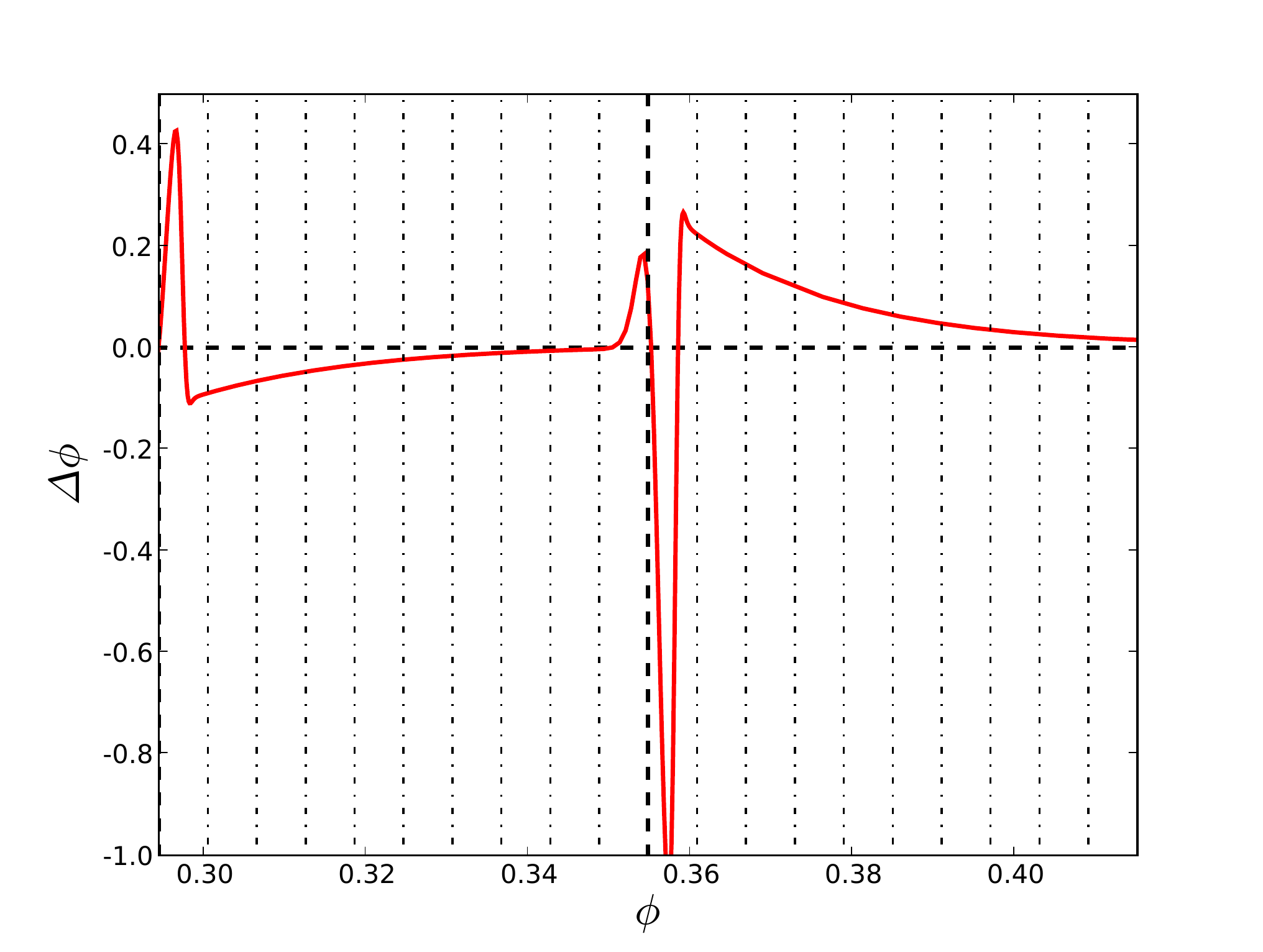}
   \end{center}
\caption[Active Segment Termination Linear Burst Phase Response ]{Close-up of linear burst phase response curve in the full HR model at the end of the active segment. The heavy dashed line marks the spike peak for spike closest to $h=2.085$ slow variable value. Light dash-dotted lines denote equally spaced phases in the full system.}
\label{fig:Iso_BPRC_2085}
\end{figure}

The periodic orbit is destroyed in the homoclinic bifurcation; afterwards, isochrons do not exist and phase has no meaning in the fast subsystem. Thus we expect a significant change in the phase response of the full system near the end of the active segment. Figure \ref{fig:Iso_BPRC_2085} shows the infinitesimal BPRC for the full system at the end of the active spiking segment. The left (penultimate) spike peak of the figure occurs at $h \approx 2.04023469758$, and the final spike, marked in the center of the figure, occurs at $h \approx 2.08623316653$. The phase response near the final spike is indeed quite different from the phase response near preceding spikes. 

%It is possible to derive analytic approximations to the PRC for spiking neural models near a homoclinic bifurcation \cite{Brown:2004}, but it is perhaps more informative to consider the phase portraits and isochron geometry of the fast subsystem in order to understand the shape of the PRC in the full system. 

%%%%%%%%%%%%%%%%%%
%	Isochrons 2.085
%%%%%%%%%%%%%%%%%%  

%%%%%%%%%%%%%%%%%%
%	PRC close-up 2.085
%%%%%%%%%%%%%%%%%%  

%\clearpage

%%%%%%%%%%%%%%%%%%
%	Upper Isochrons 2.085
%%%%%%%%%%%%%%%%%%  

%%%%%%%%%%%%%%%%%%
%	Lower Isochrons 2.085
%%%%%%%%%%%%%%%%%%  

Before the final spike peak, the full HR system responds to perturbation with relatively large phase delay, which reaches a maximum just before the spike peak. This phase delay occurs for the same reasons as the spike upswing phase delays at earlier spikes in the active segment.  Immediately after the spike peak, however, there is an extremely large, abrupt switch to phase advancement, after which there is an abrupt return to phase delay. This phase delay gradually declines; the peak of this last phase delay marks the transition between BPRC Segments I and II. The sharp phase advancement is due to the qualitative change in the fast subsystem at the homoclinic point, which occurs at a $(V,n,h)$ phase space point in the full system that projects to a point $(V,n)$ lying on the periodic orbit at a phase just before 0.0 in the fast subsystem. At the homoclinic point and beyond, the isochrons no longer exist, but one can think of their `ghosts' as merging with the periodic orbit during the homoclinic bifurcation. The outer isochrons along the depolarized portions of the periodic orbit merge with the right edge of the orbit, and their straight segments, which passed between the saddle and the periodic orbit, merge with the lower half of the periodic orbit and the stable manifold of the saddle. 

Loosely speaking, each of the isochrons become tangent everywhere to the periodic orbit, and this tangency occurs first near the spike peak at phase $0.0$. This means that the slight perturbation used to construct the BPRC no longer records phase delay around 0.0; its perturbation point no longer touches any isochrons after the homoclinic bifurcation. Instead, any small perturbation in the $V$ direction crosses the stable manifold of the saddle and leaves the perturbation point in the basin of attraction of the stable fixed point. This stable fixed point is the one that the full system tracks during quiescence; the attraction of the perturbed trajectory to the stable fixed point corresponds to advancing the onset of quiescence and therefore substantially advancing the phase of the full burst cycle.

The large magnitude of the phase advancement recorded in Figure \ref{fig:Iso_BPRC_2085} is in some sense a numerical artifact. As the homoclinic point nears, the resolvable phases on the fast subsystem periodic orbit are compressed along the hyperpolarized edge of the orbit, except for 0.0, by definition. The corresponding compression of the isochrons against the periodic orbit and the unbounded increase of the period cause a large, rapid increase in the size of the numerical estimates of phase delay given by the adjoint method. The discontinuity in phase response after the homoclinic point is resolved numerically as a huge phase advancement; the magnitude of this peak grows dramatically as the resolution of the numerical calculations is increased.

The final region of phase delay recorded after the homoclinic point can also be understood in terms of `ghosts' of the structures in the phase plane before the bifurcation. Perturbations after the spike peak occur at $(V,n,h)$ phase space points in the full system that project to points $(V,n)$ on the periodic orbit at phases just after 0.0 in the fast subsystem. Prior to the homoclinic bifurcation, the perturbation points corresponding to these phases would have lain inside the periodic orbit and above the stable manifold for the saddle. The trajectories followed by the system from the perturbation points remain inside the stable manifold, tracking it to the vicinity of the saddle, and then follow the unstable manifold of the saddle until it joins the stable fixed point. This extended excursion, mimicking the addition of a spike, is recorded as a large phase delay in the infinitesimal BPRC.

%These differences, and the ability of our isochron calculations to explain the salient features of the burst phase response curve
%\footnote{Calculated using the adjoint method advocated and used in \cite{Brown:2004, Brown:2004a}} for the full system, indicate that the assumptions and reductions used in studying the phase response dynamics of neural oscillators presented in \cite{Brown:2004, Brown:2004a, Brown:2004b} may not be valid for more realistic neuron models, especially those which endogenously burst. In particular, our results contradict an assumption used in \cite{Brown:2004, Brown:2004a} to derive analytic approximations of phase response near bifurcations, namely that phase response is approximately 0 at the peak of a spike.}
 
\section{Strong perturbation regime}
\label{sec:strong_perturbation}

Strong perturbations to the HR model, \ie $\gSyn \geq 1$, may affect the spike structure of the perturbed burst in ways that substantially influence phase response: alteration of the intraburst interspike interval, deletion of spikes from the active segment, or addition of spikes to the active segment. Perturbations that change the spike number typically induce a large change in phase, while changes in spike timing that leave the spike number unchanged have a weaker effect on phase response. Very strong and/or precisely timed perturbations may prematurely silence the active spiking segment of a burst, or prematurely initiate a new round of spiking from a hyperpolarized, quiescent state. In this section we examine the  effects of large perturbations, in particular changes in spike number, with some illustrative examples, and we relate these particular phenomena to the multiple time-scale dynamics and fast-subsystem bifurcation structure. 

In the three-dimensional figures that follow, the unperturbed (reference) burst orbit is drawn as a solid black line, and the perturbed orbit is drawn with three colors: blue for the portion of the trajectory prior to perturbation, red for the portion during the perturbation, and orange for the portion after perturbation. The same reference and perturbed orbits are depicted in two-dimensional time series representations that show at least two burst cycles; the traces for the voltage variable $V$ follow the same color key as the three dimensional figures. Also shown is the time series evolution of the slow variable $h$, which is colored cyan before the perturbation, magenta during the perturbation, and green after the perturbation. We denote the synaptic strength of excitatory and inhibitory perturbations as $\gSyne$ and $\gSyni$, respectively.

%Along with the full system trajectories, the figures show lines of fixed points from the fast subsystems for the relevant $h$ values. The line of stable fixed points is drawn as a solid green line, and a solid light blue line represents the line of saddle points; these join at a saddle-node bifurcation near the left front edge of each figure. The stable eigendirections of the saddle points are drawn in green, and the unstable eigendirections are drawn in red. Note that the values of the voltage variable, $V$, for the fixed points and saddles are always lower than the minimum voltage values of the fast-subsystem periodic orbits and of the full system trajectory during its active segment.

\subsection{Spike shift}
\label{subsec:spike_shift}

At many phases, particularly those falling in BPRC Segment II and Segment III, strong perturbations do not change the spike number of the burst. The explanation of the response to strong perturbation at these phases is essentially the same as for weak perturbations: For most of the quiescent segment, the full system trajectory tracks \MQ, the curve of stable fixed point of the fast subsystem, returning rapidly to the fixed point and thus incurring little phase change after voltage perturbations. Perturbations near the end of the quiescent segment, where the stability of the quiescent fixed point is waning, may move the trajectory into the domain of attraction for the stable periodic orbit and so accelerate the onset of the active segment of the next burst cycle, thereby inducing phase advancement. In neither of these two cases does the spike number change in the current burst, nor is the intraburst interspike timing affected, since the active segment of the current burst cycle is already complete. 

%%%%%%%%%%%%%%%%%%%%%%%%%%%%%%%%%%%%
%	Spike shift with inhibition and excitation
%%%%%%%%%%%%%%%%%%%%%%%%%%%%%%%%%%%%  
\begin{figure}[!ht]
     \begin{minipage}{3in}
     \begin{center}
	\includegraphics[height= 2.5in]{\figpath/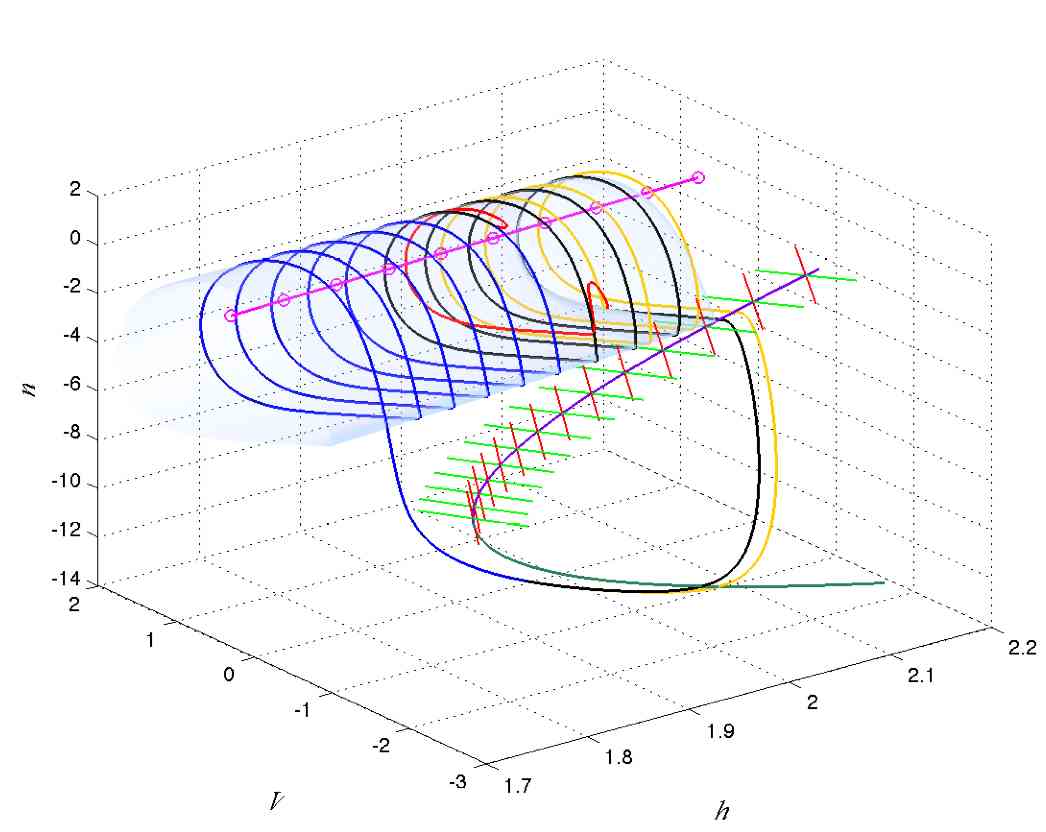}
	(a) %$\gSyni = 0.1$, $\theta = 0.265$
	    \end{center}
  \end{minipage}
  \begin{minipage}{3in}
      \begin{center}
	  \includegraphics[height= 2.5in]{\figpath/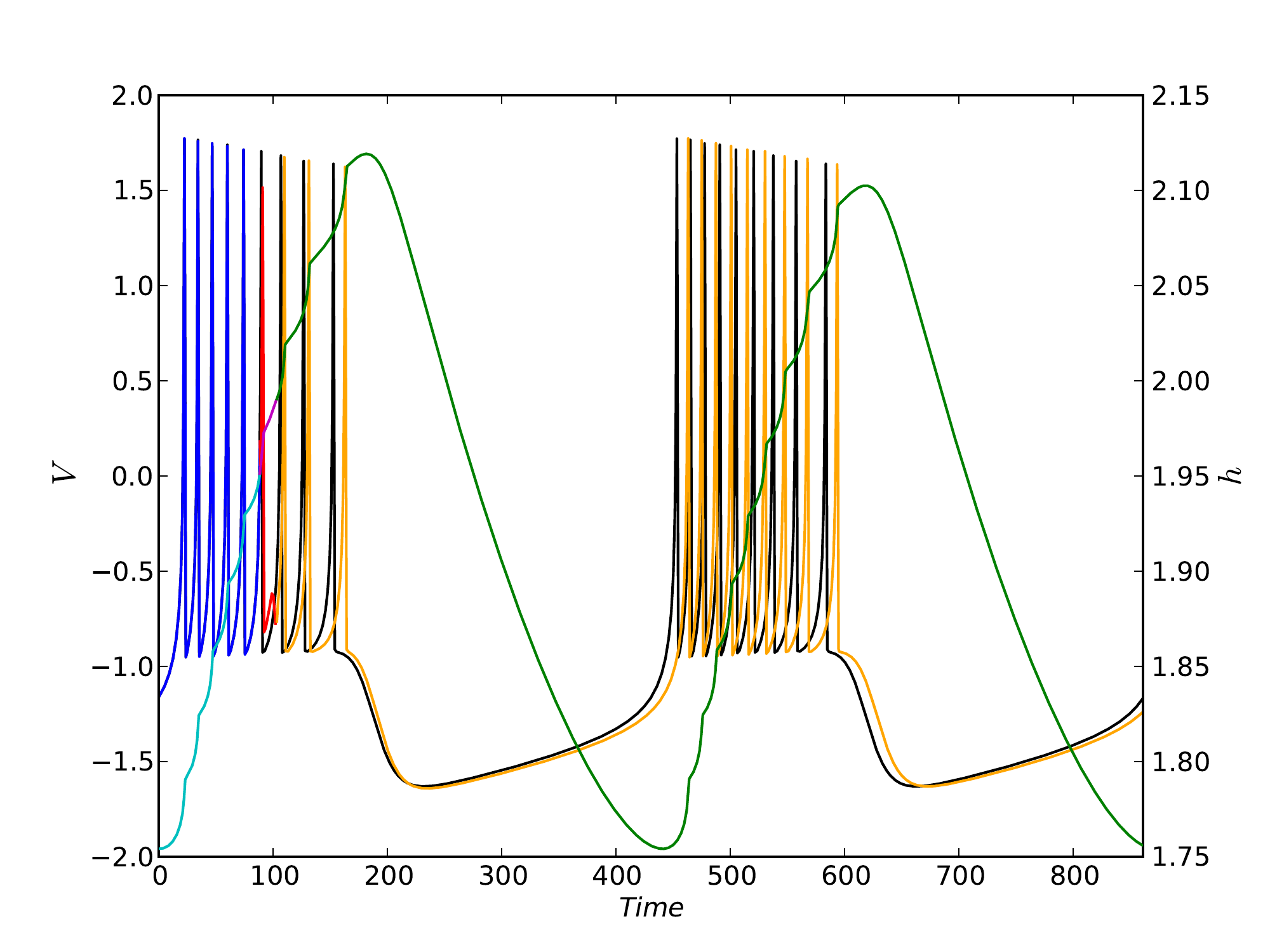}
	  (b) %$\gSyni = 0.1$, $\theta = 0.265$
    \end{center}
  \end{minipage}%
  
\caption[Spike Shift]{Shift in intraburst interspike timing due to inhibitory perturbation: $\gSyni = 1.0, \theta = 0.205538462$. The unperturbed reference burst orbit is drawn in black. The segment of the perturbed burst trajectory prior to spike injection is drawn in blue; during spike injection in red; after spike injection in orange. Note that the perturbation shifts the trajectory away from the line of saddle points and inside the manifold of fast subsystem periodic orbits, \MP.}
\label{fig:Spike_Shift_3D}
\end{figure}

Perturbations at phases in BPRC Segment I---during the active segment of the burst---which do not affect the spike number typically alter the interspike intervals for the spikes following the perturbation, as illustrated in Figure \ref{fig:Spike_Shift_3D} for inhibitory spike injection. (Perturbations at the end of Segment III may also alter interspike intervals in the subsequent burst cycle.) Note that the inhibitory perturbation depicted spans the downswing of one spike and the trough at the base of the following (anticipated) spike, when the voltage value of the full system trajectory brings it closest to the line of fast subsystem saddle points, \MS. The trajectory is initially perturbed inside the tube-shaped manifold of fast subsystem periodic orbits, \MP\ (the unperturbed trajectory circles the outside of this manifold), then recovers to track the outside of the tube. The voltage variable $V$ is depressed from its usual value during the spike upswing, and the inhibition also deflects the perturbed trajectory farther away than usual from \MS. (This is typical for both inhibitory and excitatory perturbations resulting spike shifts: the perturbation does not bring the trajectory close enough the saddles for there to be any interaction.) By decreasing $V$, the inhibitory perturbation also slows the change in $h$, relative to its value in the reference trajectory. This delays the next spike peak, a timing change which propagates to the subsequent spike times, as seen in Figure \ref{fig:Spike_Shift_3D} (b). For perturbation strengths $\gSyne \geq 1$, excitatory perturbations in Segment I always alter the spike number of the perturbed burst, but excitatory perturbations in the transitional range between weak and strong perturbation strengths, \eg at $\gSyne = 0.1$, may simply shift intraburst interspike intervals. The typical effect of such an excitatory perturbation is to increase $V$ and accelerate the change in $h$, triggering an early spike peak and shifting the immediately following spikes to earlier times, thus inducing phase advance. 

A critical factor in determining magnitude of phase advance or delay, however, is the position of perturbed trajectory with respect to the homoclinic point. In the homoclinic bifurcation, the fast subsystem periodic orbit merges with the stable and unstable manifolds of the coexistent saddle point. The final spike of the unperturbed full system begins at an $h$ value  where the periodic orbit still exists in the family of fast subsystems, and it finishes at an $h$ value past the homoclinic point, so that it is drawn to the stable fixed point corresponding to quiescence. Perturbations prior to the homoclinic point may change the perturbed orbit's distance to the homoclinic point for the final spike. The closer the end of the full system spike is to the homoclinic point, the longer the final spike lasts because it tracks periodic orbits with periods diverging to infinity; the final interspike interval depends on this positioning relative to the homoclinic point. In addition, the perturbed trajectory may cross \MS\ at a point with a different value of the slow variable  $h$, and this difference may change the $h$ value at which the trajectory subsequently arrives at the line of quiescent fixed points. Arrival at \MQ\ at an increased $h$ value (assuming that $V$ is the same as in the unperturbed case) results in a longer recovery time between bursts, and thus (usually) phase delay, while arrival at a reduced $h$ value normally has the opposite result. Figure \ref{fig:Spike_Shift_3D} shows an example of the former case: the perturbed trajectory leaves \MP\ and crosses the line of saddles at a point further beyond the homoclinic point than usual, and it arrives at the line of quiescent fixed points at a higher $h$ value than normal. Recovery along \MQ\ to the saddle node bifurcation demarcating the start of the next burst cycle takes longer, hence the result is phase delay. 

%Even without changing spike number, the effect of excitation can be phase delay, although it typically accelerates the advent of the subsequent spike (see figure \ref{fig:HR_excite_BPRC_0} (c) and (d)). As explained below in subsection \ref{subsec:spike_add}, this paradoxical effect is due to changes in the interspike intervals near the end of the active segment, where the fast subsystem is approaching the homoclinic bifurcation point.

\subsection{Spike addition}
\label{subsec:spike_add}

Spike addition may occur in a number of different ways, each critically related to the proximity of the perturbed orbit to the line of saddle points, the behavior of the fast subsystem near the homoclinic  bifurcation at the end of the active segment, or both. Both excitation and inhibition may act to add one or more spikes to the burst, and the result may be phase advance or phase delay, depending on the exact circumstances.

\begin{figure}[!ht]
 \begin{minipage}{3in}
     \begin{center}
	\includegraphics[height= 2.5in]{\figpath/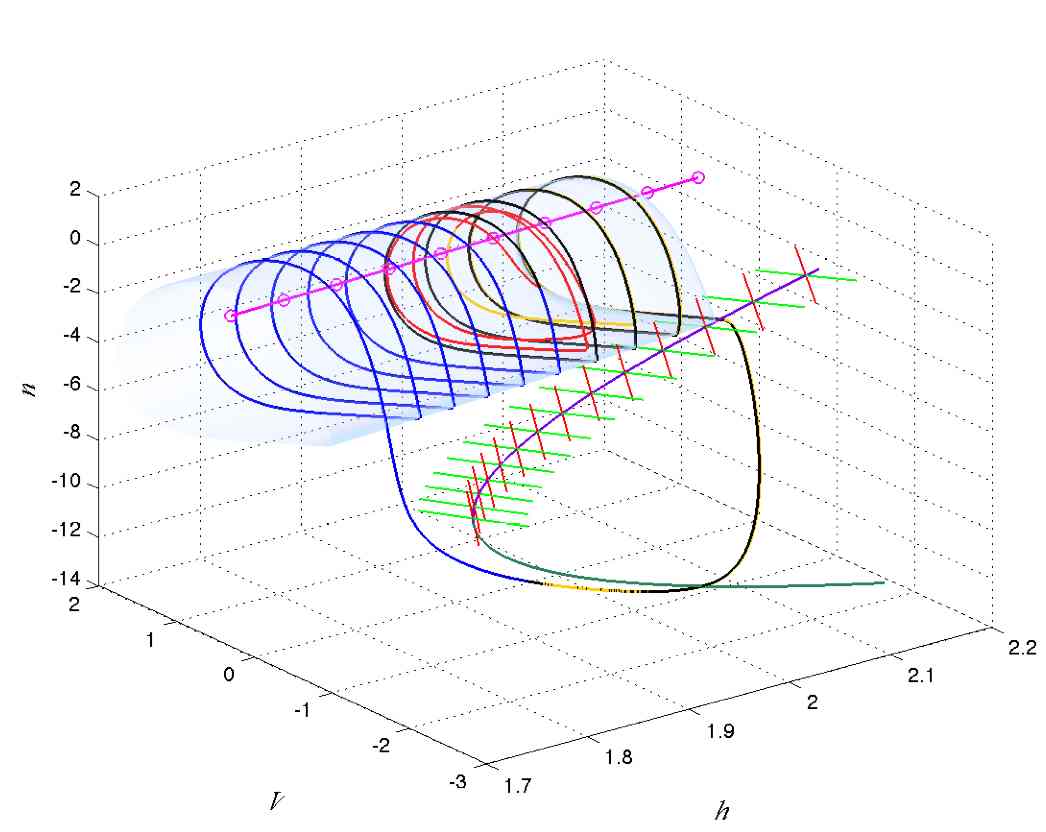}
	(a)
	    \end{center}
  \end{minipage}
  \begin{minipage}{3in}
      \begin{center}
	  \includegraphics[height= 2.5in]{\figpath/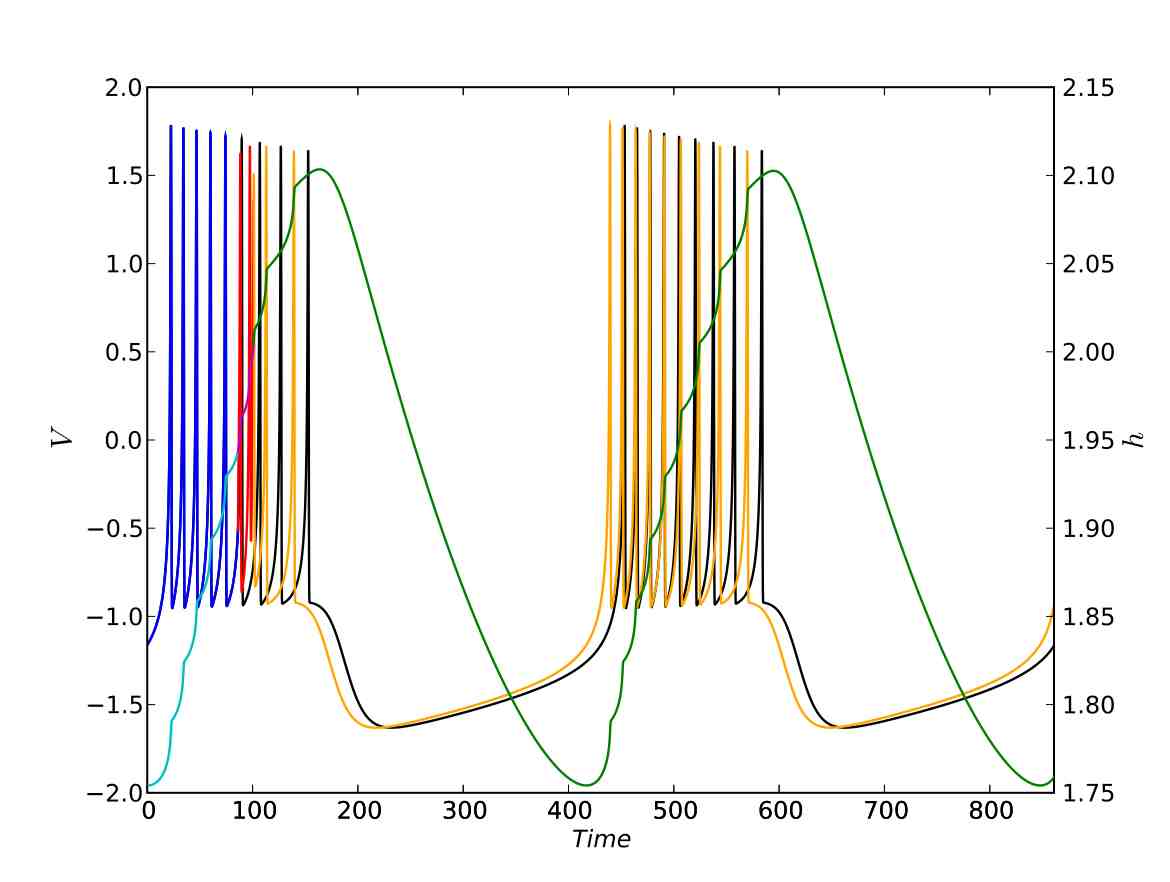}
	  (b) 
    \end{center}
  \end{minipage}%
  
\caption[Spike Addition, Excitation]{Addition of spikes via excitation far from the homoclinic point, $\gSyne = 1.0$, $\theta = 0.2$. During the perturbation, the trajectory moves inside the manifold of fast subsystem periodic orbits, and two spikes are emitted in the usual time span of a single spike. The net effect is slight phase advancement.}
\label{fig:Spike_Add_Excite1}
\end{figure}

The first case we consider is illustrated in Figure \ref{fig:Spike_Add_Excite1}, when spike addition is triggered far from the homoclinic point by excitation. Here the perturbation occurs in the middle of the active spiking segment but during the upswing of an anticipated spike. The excitatory perturbation pushes the trajectory inside \MP\ for the duration of that first expected spike, and holds the trajectory relatively far away from the hyperpolarized edge of \MP, and also away from \MS. The evolution of the variables $V$ and $n$ is more rapid in this region of phase space, and the trajectory passes quickly through the relatively hyperpolarized portion of this first spike loop, beginning a second spike loop while the excitatory perturbation is still operative. The second (new) spike loop is similar to the first, and the perturbation ends near the start of the upswing portion of a new spike. Without the extra perturbatory input, the trajectory recovers to track the outside of \MP, and the remaining two spikes trace nearly the same path through phase space as the final two spikes of the unperturbed reference burst orbit. In effect, excitation has inserted a spike between two `regular' spikes of the unperturbed burst. The perturbed trajectory crosses the homoclinic point slightly earlier than the unperturbed trajectory, and the end result is a slight phase advancement, despite the additional spike.

%The period of the orbit increases as the bifurcation approaches, becoming infinite at the homoclinic point. The periodic orbit no longer exists after the homoclinic bifurcation, and the branch of the unstable manifold of the saddle which formerly wrapped around the periodic orbit now terminates at the stable fixed point. For the full system trajectory tracking these fast subsystem objects, this means that the interspike interval grows as the end of the active segment nears. 

%%%%%%%%%%%%%%%%%%%%%%%%%%%%%%%%%%%%
%	Spike addition with excitation 1
%%%%%%%%%%%%%%%%%%%%%%%%%%%%%%%%%%%%  

\begin{figure}[!ht]
 \begin{minipage}{3in}
     \begin{center}
	\includegraphics[height= 2.5in]{\figpath/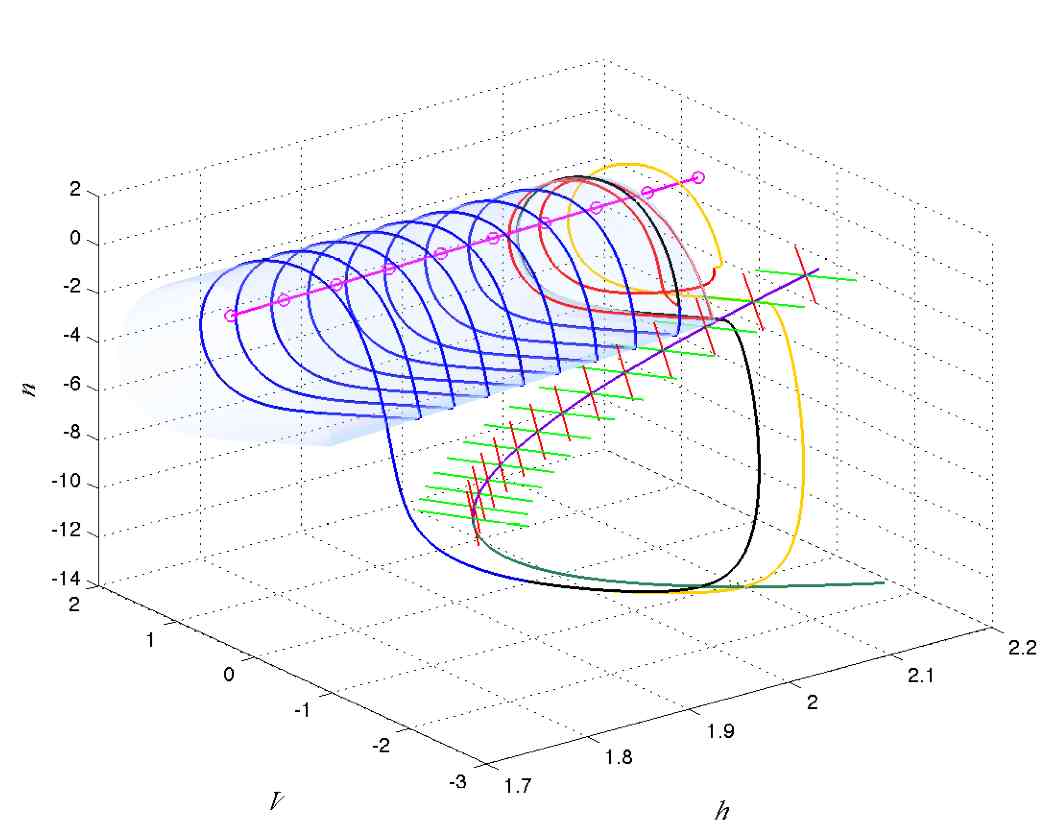}
	(a)
	    \end{center}
  \end{minipage}
  \begin{minipage}{3in}
      \begin{center}
	  \includegraphics[height= 2.5in]{\figpath/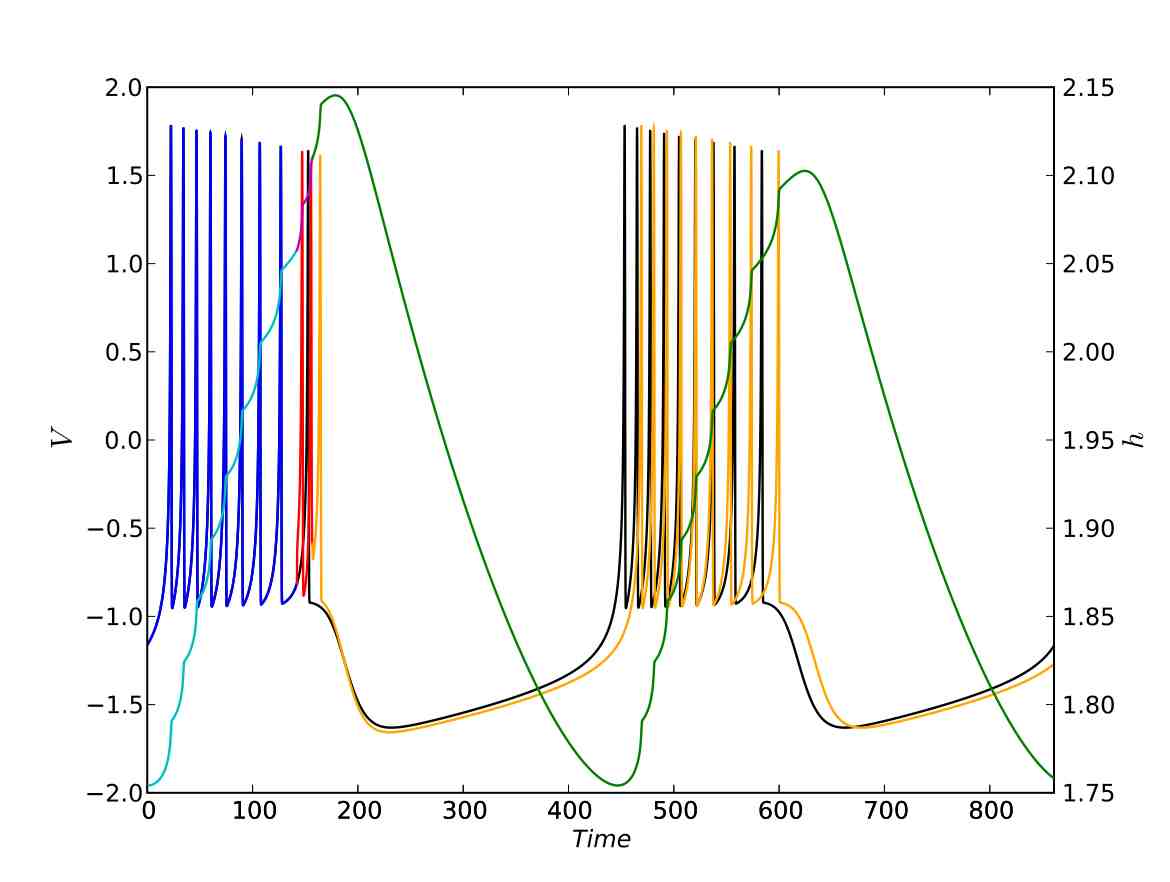}
	  (b) 
    \end{center}
  \end{minipage}%
  
\caption[Spike Addition, Excitation]{Addition of spikes via excitation far from the homoclinic point, $\gSyne = 1.0$, $\theta = 0.33$. During the perturbation, the trajectory moves inside \MP and rapidly emits two spikes; past the homoclinic point, the trajectory loops a third time around the line of unstable fixed points. The end result is slight phase delay.}
\label{fig:Spike_Add_Excite2}
\end{figure}

The example in Figure \ref{fig:Spike_Add_Excite2} depicts an excitatory perturbation near the beginning of the final spike before the homoclinic point. As in the previous case, the perturbation pushes the trajectory inside the manifold of fast subsystem periodic orbits and prevents the trajectory from crossing the line of saddles during the downswing of the first (anticipated) spike. The perturbed trajectory follows the upper edge of \MP\ through the homoclinic point, emitting a second spike. The perturbation ends after the homoclinic point, where \MP\ no longer exists, but the perturbed trajectory makes a third loop around the line of unstable fixed points before crossing \MS. The perturbation effectively appends two spikes to the end of the active spiking segment. The perturbed trajectory enters the basin of attraction for the quiescent fixed point at a higher $h$ but lower $V$ value, and the result is a moderate phase delay.

%%%%%%%%%%%%%%%%%%%%%%%%%%%%%%%%%%%%
%	Spike addition with excitation 2
%%%%%%%%%%%%%%%%%%%%%%%%%%%%%%%%%%%%  

Perturbations may also coincide with the final spike loop at a point past the homoclinic bifurcation in the fast subsystem. At this point, the trajectory is no longer following the edge of \MP\ and may be about to cross or have just crossed \MS. A perturbation may then displace the trajectory such that it either does not cross \MS\ or it recrosses \MS\ and reenters the region of phase space around the unstable fixed point bounded by the stable and unstable manifolds of the saddles. In this region, the trajectory completes one or more loops about the unstable fixed point before finally crossing the line of saddles and entering the basin of attraction for the quiescent stable fixed points. 

%%%%%%%%%%%%%%%%%%%%%%%%%%%%%%%%%%%%
%	Spike addition with excitation after homoclinic
%%%%%%%%%%%%%%%%%%%%%%%%%%%%%%%%%%%%  
\begin{figure}[!ht]
 \begin{minipage}{3in}
     \begin{center}
	\includegraphics[height= 2.5in]{\figpath/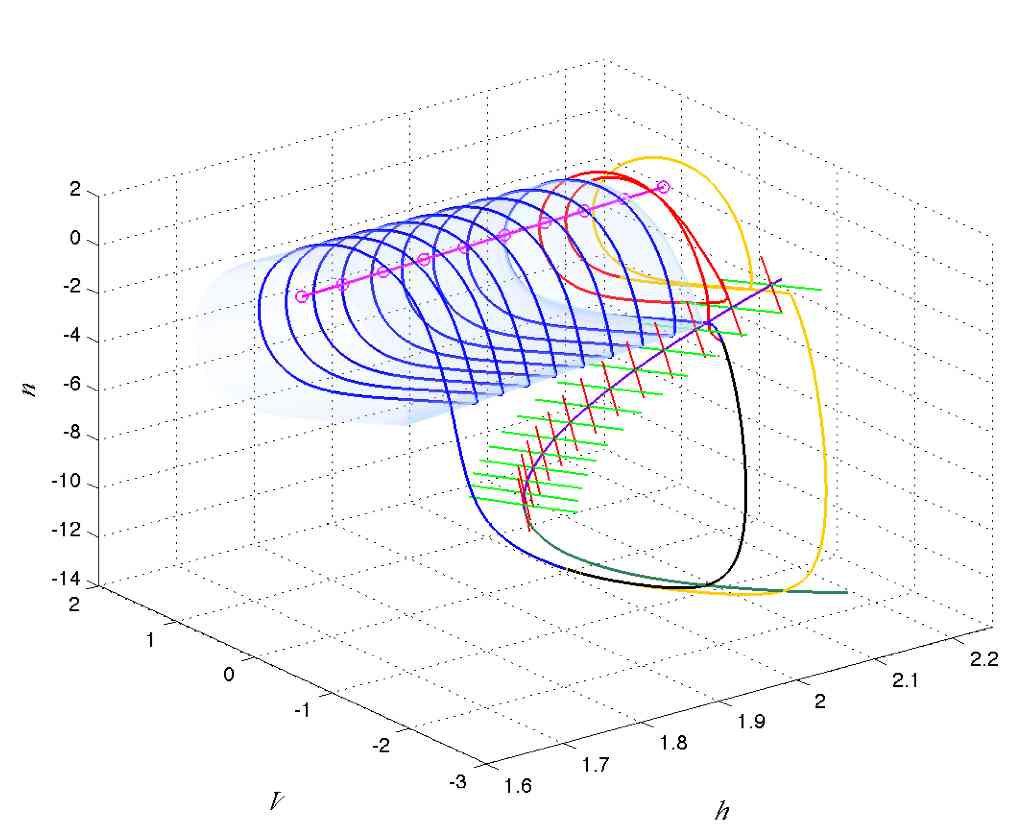}
	(a) 
	    \end{center}
  \end{minipage}
     \begin{minipage}{3in}
     \begin{center}
	\includegraphics[height= 2.5in]{\figpath/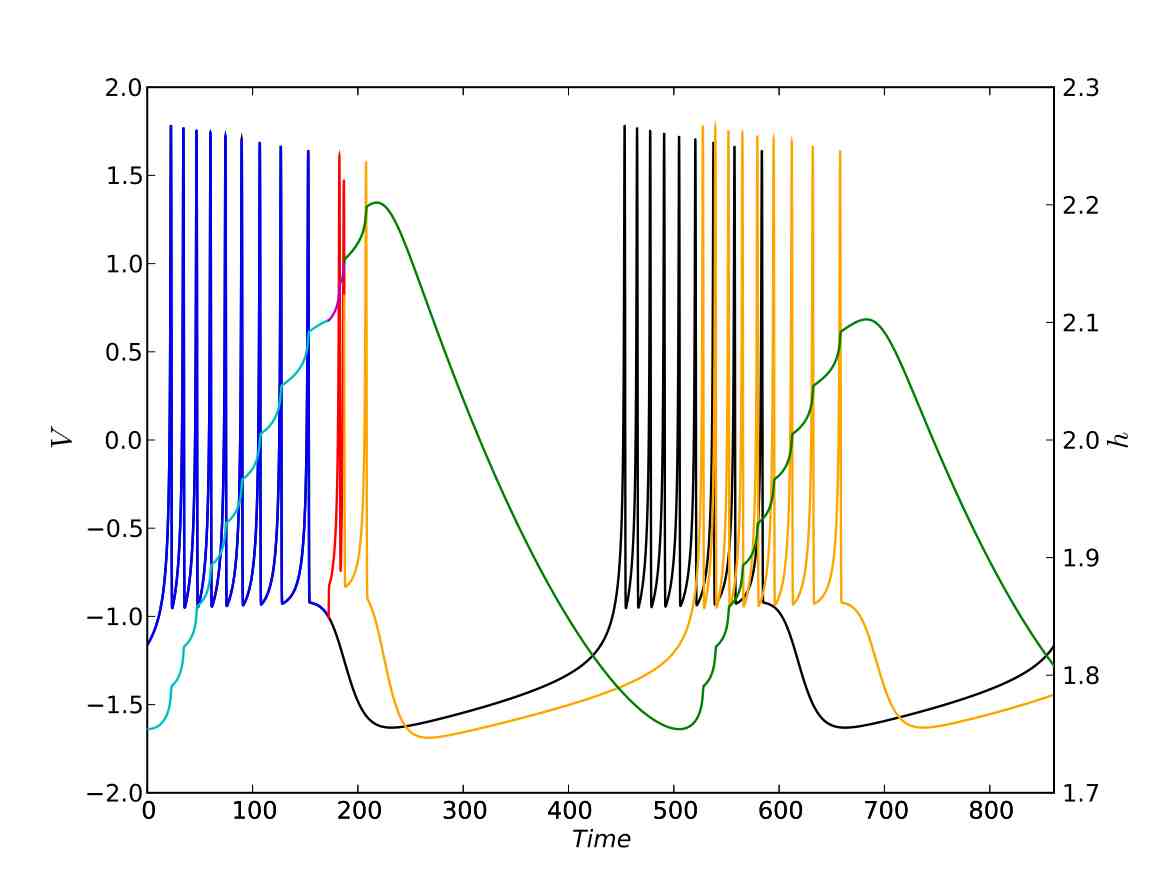}
	(b)
	    \end{center}
  \end{minipage}  
\caption[Spike Addition, Excitation]{Addition of spikes via an excitatory perturbation occurring after the homoclinic point, $\gSyne = 1.0$, $\theta = 0.4$. Note that prior to perturbation, the trajectory had just crossed the line of saddle points.}
\label{fig:Spike_Add_Excite_Late}
\end{figure}

\begin{figure}[!ht]
 \begin{minipage}{3in}
     \begin{center}
	\includegraphics[height= 2.5in]{\figpath/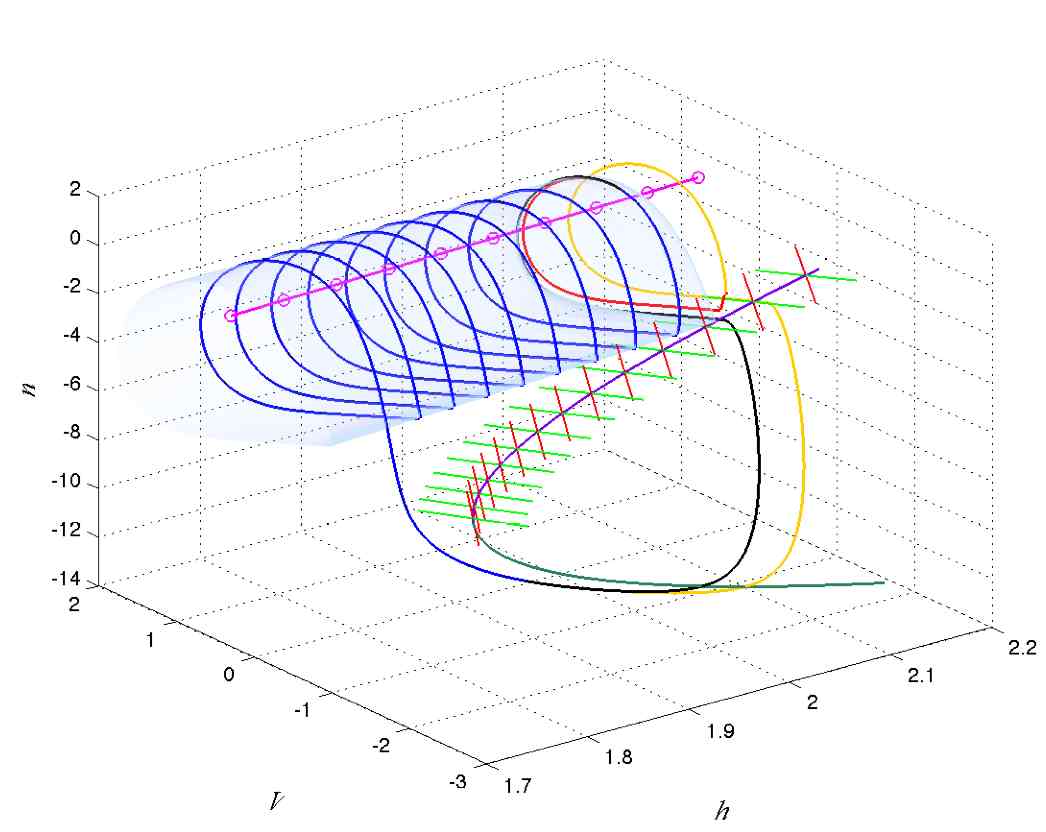}
	(a) 
	    \end{center}
  \end{minipage}
     \begin{minipage}{3in}
     \begin{center}
	\includegraphics[height= 2.5in]{\figpath/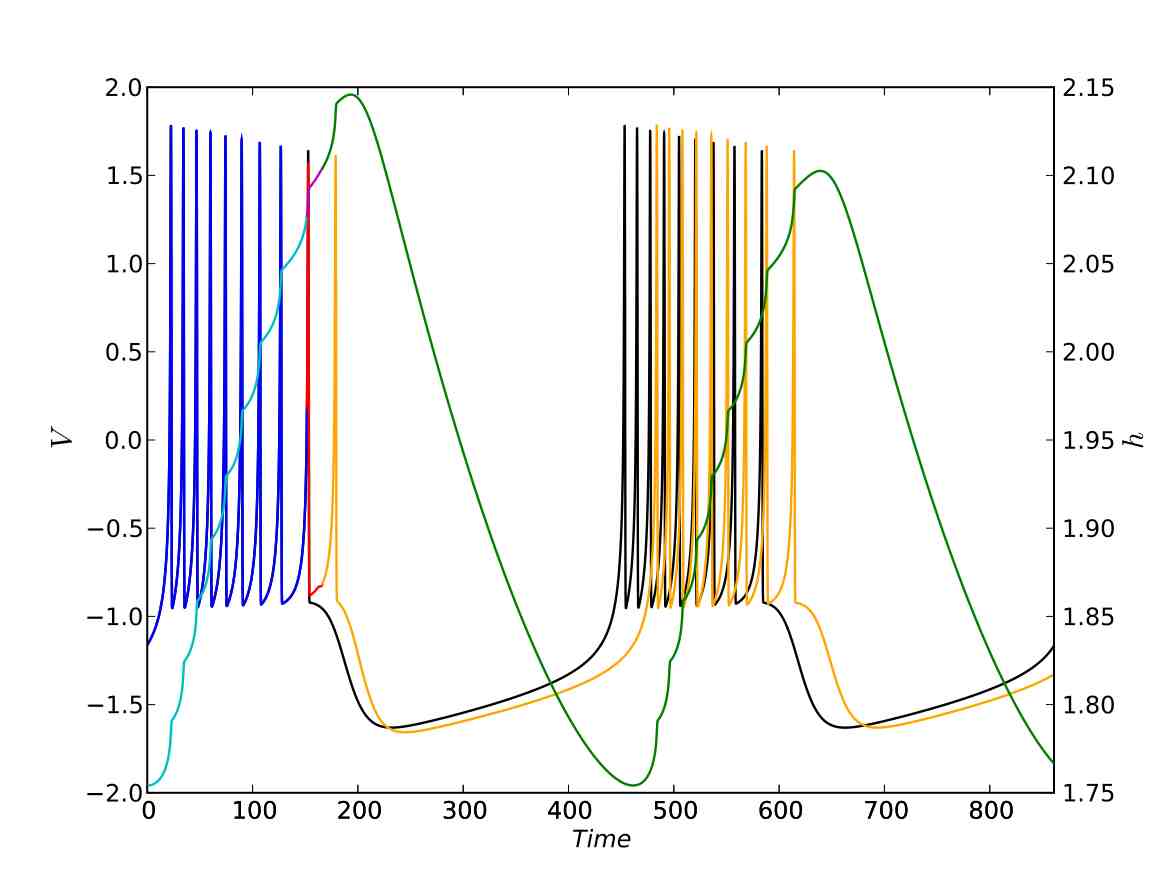}
	(b)
	    \end{center}
  \end{minipage}  
\caption[Spike Addition, Inhibition]{Addition of spikes via inhibition occurring immediately prior to the homoclinic point, $\gSyni = 1.0$, $\theta = 0.4$. Note that prior to perturbation, the trajectory has not crossed the line of saddle points.}
\label{fig:Spike_Add_Inhib}
\end{figure}

Figure \ref{fig:Spike_Add_Excite_Late} shows an example for an excitatory perturbation that begins just after the trajectory has crossed the line of saddles. The depolarizing action of the excitation holds the perturbed trajectory away from \MS\ long enough for two extra spike loops to be appended to the burst. Figure \ref{fig:Spike_Add_Inhib} depicts an inhibitory perturbation occurring after the homoclinic point; this perturbation begins during the upswing of the final spike loop but before the trajectory has crossed \MS. Again, the effect of the perturbation is to deflect the trajectory such that it remains within the unstable manifold of the saddle point and makes an additional loop about the unstable fixed point to append a spike. Both examples show phase delay, and as before, the values of $h$ at which the perturbed trajectory finally crosses \MS\ and subsequently arrives at \MQ\ is a key determinant of the magnitude of the phase shift.

%%%%%%%%%%%%%%%%%%%%%%%%%%%%%%%%%%%%
%	Spike addition with inhibition
%%%%%%%%%%%%%%%%%%%%%%%%%%%%%%%%%%%%  

\subsection{Spike deletion}
\label{subsec:spike_delete}

\begin{figure}[!ht]
 \begin{minipage}{3in}
     \begin{center}
	\includegraphics[width= 3in]{\figpath/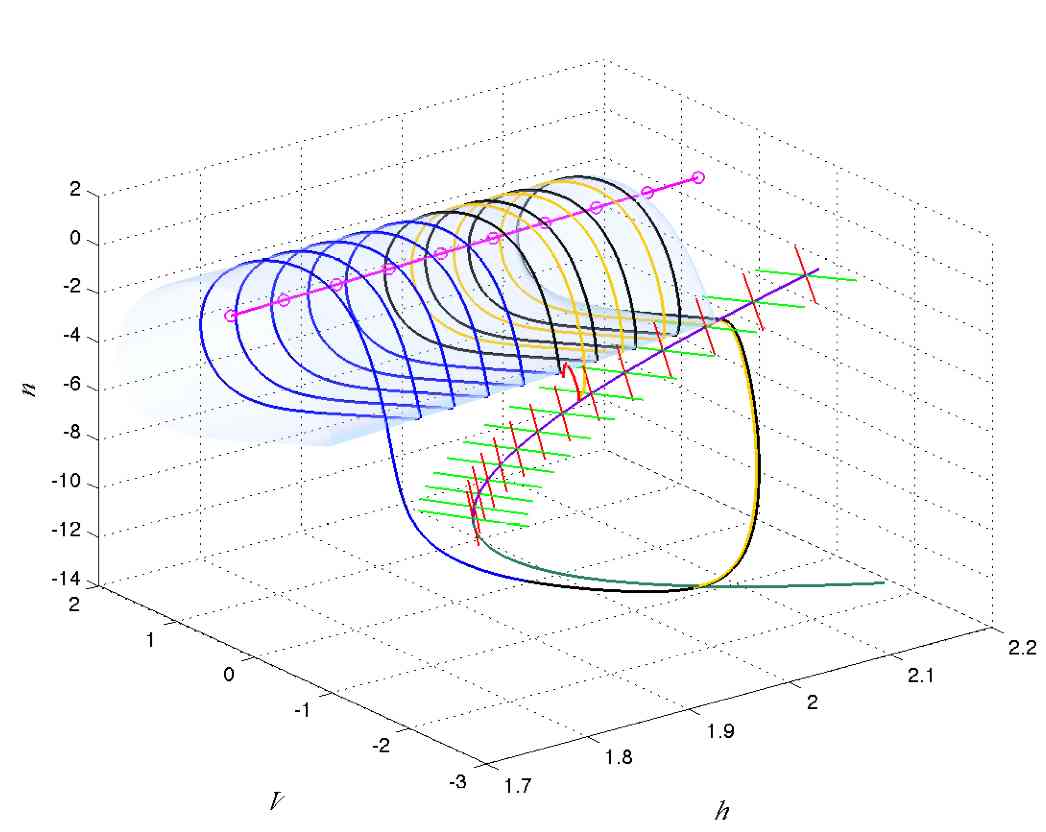}
	(a) 
	    \end{center}
  \end{minipage}
     \begin{minipage}{3in}
     \begin{center}
	\includegraphics[width= 3in]{\figpath/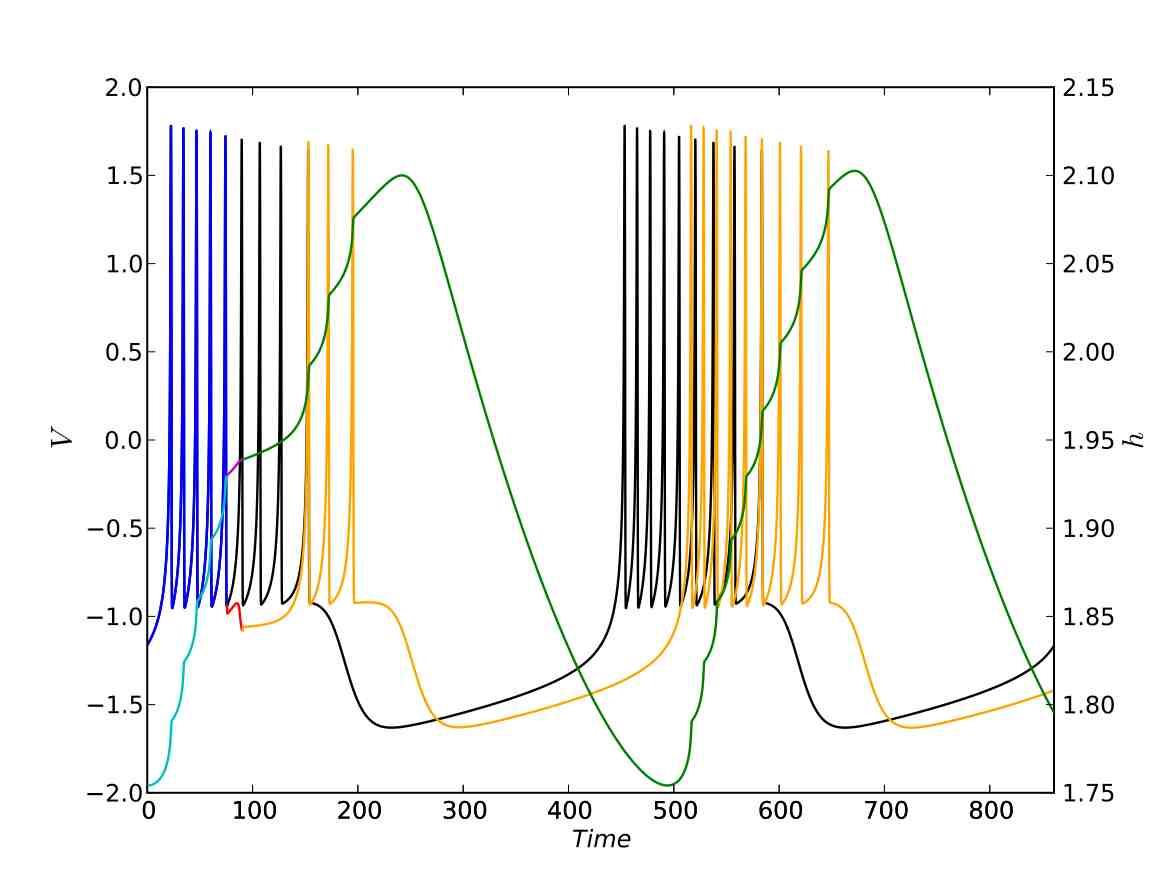}
	(b)
	    \end{center}
  \end{minipage}  
\caption[Spike Deletion, Inhibition]{Deleting a spike via inhibition, $\gSyni = 1.0$, $\theta = 0.176$.  The perturbation deflects the trajectory a small distance in phase space, away from \MP\ and very close to the line of saddles. Following perturbation,  trajectory transiently tracks the unstable manifold of the saddles in a canard-like fashion, so that $V$ and $h$ change only slowly and slightly for some time. Active spiking resumes after the transient and terminates in a nearly normal fashion, but a large phase delay is incurred. }
\label{fig:Spike_Delete_Inhib}
\end{figure}

Spike deletion without premature termination of the burst may happen due to an inhibitory perturbation away from the homoclinic point, \ie occurring one or more spikes before the final spike loop of the unperturbed burst orbit. Figure \ref{fig:Spike_Delete_Inhib} presents an example. The perturbation begins at the base of a spike, at the `corner' of the fast subsystem periodic orbit where $V$ and $n$ are both near their lowest values. The hyperpolarizing action of the inhibition is enough to deflect the trajectory away from \MP\ to the vicinity of the saddle line, but it does not bring the trajectory across the unstable manifold of the local saddle point. Instead, the trajectory begins to track the unstable sheet emanating from \MS. The result is a slow, canard-like motion during which there is only relatively slight change in $V$ and slow change in $h$. The full HR system appears to be silent, but depolarized, for some period of time, before resuming spiking. Geometrically, the perturbed trajectory eventually escapes from the unstable manifold of the saddles and resumes circling \MP\ at an $h$ value midway between where the two proximal spike loops of the reference burst orbit would lie. As the perturbed trajectory finishes out the active spiking segment of its burst cycle, it loops around \MP\ but completes one loop fewer than the reference burst, essentially skipping the spike nearest the time of the perturbation. The perturbed trajectory finally crosses \MS\ at a point just slightly ahead of where the unperturbed reference orbit crosses. In this scenario, the $h$ and $V$ values at which the perturbed trajectory terminates spiking and arrives at \MQ\ are nearly identical to those for the unperturbed orbit, but there is nevertheless a large phase delay before the onset of the next burst cycle. This is due to the extended transient period associated with the perturbed trajectory's canard-like tracking of the unstable manifold of the saddles near the point of perturbation.

%Spike deletion occurs when the full system trajectory is knocked away prematurely from the basin of attraction around the stable periodic orbit  and enters the basin of attraction of the stable fixed point. In this case, the burst cycle ends early as one or more spikes are dropped from the end of the active segment. Inhibitory perturbation drives the voltage value to a hyperpolarized level that is more negative than the line of saddle points in the fast subsystem. The perturbed trajectory rapidly approaches the saddles along the stable eigendirection and then quickly leaves along the unstable eigendirection once it has crossed the line of saddle points, as shown in figure \ref{fig:Spike_Del_3D}. 

%%%%%%%%%%%%%%%%%%%%%%%%%%%%%%%%%%%%
%	Spike deletion with inhibition
%%%%%%%%%%%%%%%%%%%%%%%%%%%%%%%%%%%%  

\subsection{Early burst initiation and termination}
\label{subsec:burst_change}

The mechanisms of spike number change so far considered involve changes to the intraburst spike structure in which the path taken by the perturbed trajectory is altered, but the trajectory continues to circle the line of  unstable fixed points and (possibly) track the manifold of stable fast subsystem periodic orbits until crossing the line of saddle points to approach the line of quiescent fixed points and finally end active spiking. More drastic alterations may occur that force the perturbed trajectory either into or out of the basin of attraction of the quiescent fixed point, thereby terminating spiking early or starting a new round of spiking (and a new burst cycle) prematurely.

%%%%%%%%%%%%%%%%%%%%%%%%%%%%%%%%%%%%
%	Burst termination with inhibition
%%%%%%%%%%%%%%%%%%%%%%%%%%%%%%%%%%%%  
 \begin{figure}[!ht]
 \begin{minipage}{3in}
     \begin{center}
	\includegraphics[width= 3in]{\figpath/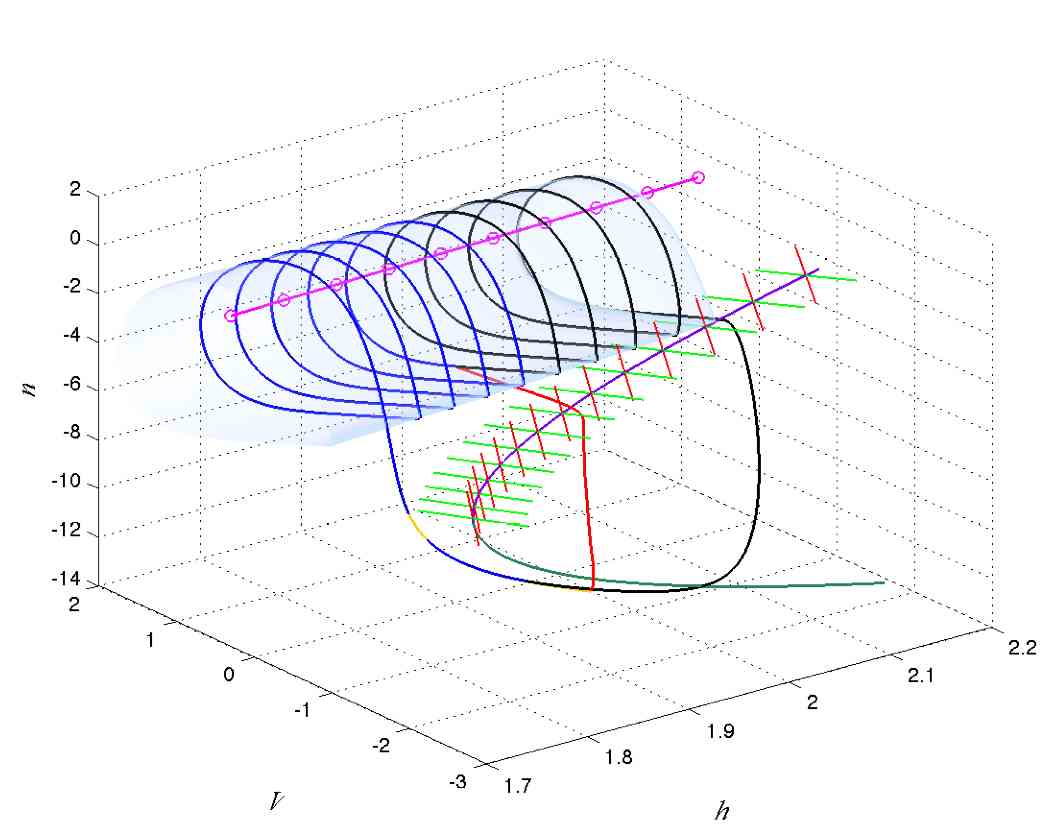}
	(a) 
	    \end{center}
  \end{minipage}
     \begin{minipage}{3in}
     \begin{center}
	\includegraphics[width= 3in]{\figpath/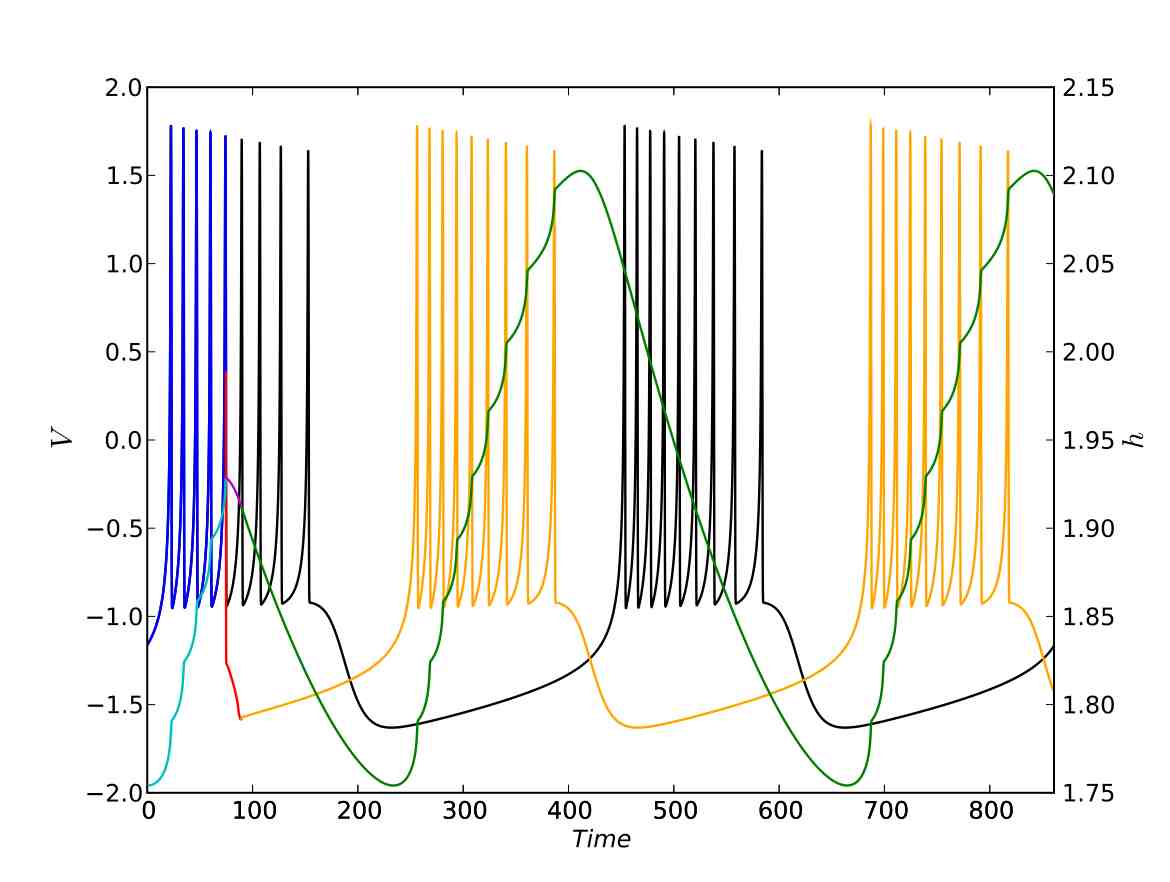}
	(b)
	    \end{center}
  \end{minipage}  
\caption[Burst Termination, Inhibition]{Early burst termination via inhibition, $\gSyni = 10.0$, $\theta = 0.17353846$. The subsequent burst cycle tracks the reference burst orbit very closely.}
\label{fig:Burst_Termination}
\end{figure}

Figure \ref{fig:Burst_Termination} presents an example of early burst termination due to an inhibitory perturbation at a point distant from the homoclinic bifurcation in the fast subsystem. The perturbation, which occurs during the downswing of a spike, pushes the trajectory along the stable eigendirection of the saddle and across \MS. Once past the saddles, the trajectory is strongly and immediately attracted to the line of stable quiescent fixed points. The remaining spikes that would be present  in the unperturbed orbit are omitted. The full system trajectory then tracks \MQ\ as $h$ recovers. Though the figure shows the trajectory traveling towards the quiescent fixed points during the perturbation (red line), it is not necessary for the perturbation to persist until the trajectory reaches \MQ, but rather only until the trajectory has completely crossed \MS\ and entered the basin of attraction for the quiescent fixed point. The duration of the quiescent segment after perturbation is shorter than normal in proportion to the difference in $h$ values at the end of active segments of the perturbed and reference orbits. By the onset of the next burst cycle, the perturbed trajectory recovers to follow a path nearly identical to the reference burst orbit. The net effect of the perturbation is to shut off the active segment of the burst early, deleting any remaining spikes, thus resetting the trajectory for a new burst cycle and causing a large phase advancement.

Though not shown in Figure \ref{fig:Burst_Termination}, the absolute values of both eigenvalues for the saddle point grow by about one and a half orders of magnitude as $h$ increases and the homoclinic point approaches, though the stable eigenvalue is always about three orders of magnitude larger than the unstable eigenvalue. Similarly, the eigenvalues for the stable fixed point also increase in magnitude towards the homoclinic point, with the stronger eigenvalue typically one or two orders of magnitude larger than the stable eigenvalue of the saddle at the same $h$ value. At the same time, the stable manifolds of the saddle points lie very close to the periodic orbits near the homoclinic point. This increase in the strength of attraction to and repulsion from the fixed points, along with the greater proximity of  periodic orbits to the saddle points and their stable manifolds means that perturbations can more easily push the full system trajectory past \MS, across the saddle points' stable manifolds, and towards \MQ at $h$ values closer to the homoclinic bifurcation that terminates the active segment. Said differently, perturbations of a given strength that occur near the end of the active segment are more likely to drive the system into quiescence, prematurely silencing the burst, than those near the start of the active segment. Thus large phase advancement due to early burst termination first emerges closer to the homoclinic point as perturbation strength increases, as can be seen by comparing the BPRCs in Figure \ref{fig:HRIi_BPRCs} (b), (c), and (d). 
 
%%%%%%%%%%%%%%%%%%%%%%%%%%%%%%%%%%%%
%	Burst initiation with excitation
%%%%%%%%%%%%%%%%%%%%%%%%%%%%%%%%%%%%  
 \begin{figure}[!ht]
 \begin{minipage}{3in}
     \begin{center}
	\includegraphics[width= 3in]{\figpath/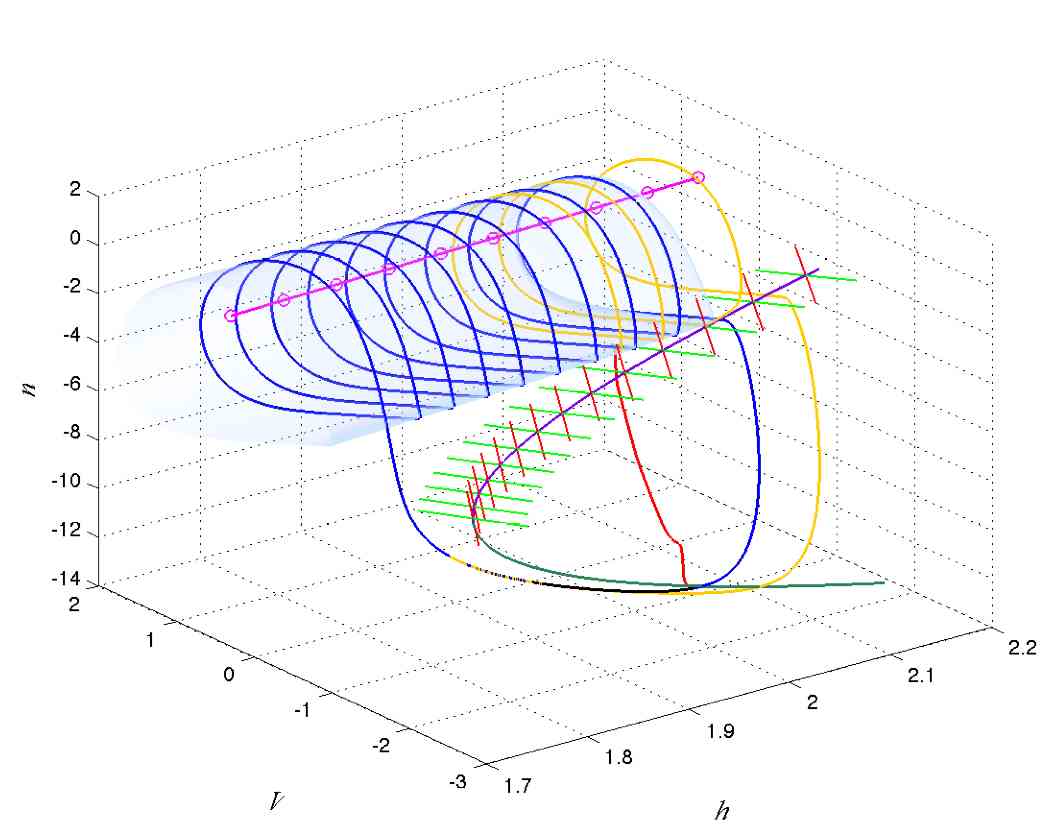}
	(a) 
	    \end{center}
  \end{minipage}
     \begin{minipage}{3in}
     \begin{center}
	\includegraphics[width= 3in]{\figpath/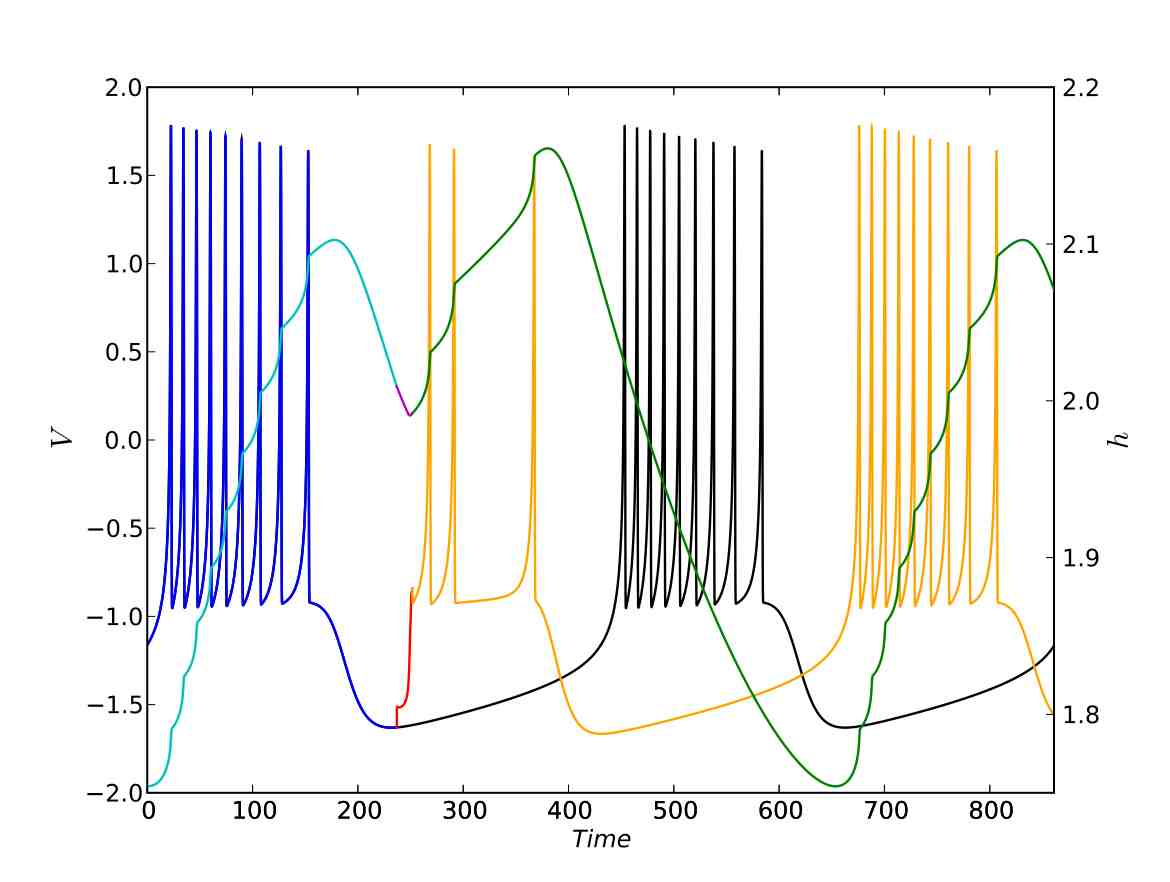}
	(b)
	    \end{center}
  \end{minipage}  
\caption[New Burst, Excitation]{Premature initiation of a new burst via excitation, $\gSyne = 1.0$, $\theta = 0.4$. The second burst cycle following the perturbation tracks the reference burst orbit very closely.}
\label{fig:Spike_NewBurst_Excite}
\end{figure}

Early burst initiation follows a sequence of events reversed from early burst termination. As shown in Figure \ref{fig:Spike_NewBurst_Excite}, strong excitation applied during the quiescent segment of the burst cycle may depolarize the perturbed trajectory far enough away from the line of hyperpolarized, quiescent fixed points that it crosses the line of saddle points and enters the basin of attraction of the fast subsystem periodic orbits. Once there, the perturbed trajectory traces out spike loops as it tracks \MP\ to the homoclinic point, where spiking terminates in an essentially normal fashion. Notice that this prematurely initiated round of spiking is truncated in comparison to a normal burst, and the intraburst spike structure (interspike intervals and spike shapes) differs from the normal active spiking segment. Subsequent burst cycles closely track the original, unperturbed reference orbit. 

The stability of the quiescent fixed point wanes as the  saddle-node bifurcation marking the start of a new burst cycle approaches, and hence the strength of excitation required to start spiking prematurely decreases closer to the saddle-node bifurcation, and the BPRCs for even weak excitation record an increase in phase advancement in Segment III. As the strength of excitatory perturbation increases, the left endpoint of the region of phases in Segment III at which early burst initiation is possible decreases from $\theta \approx 1.0$,  broadening the region of large phase advancement at the end of the burst cycle. At high excitatory perturbation strengths, perturbation anywhere in the quiescent segment of the burst initiates a new burst cycle, and thus the magnitude of phase advancement is proportional to the phase of perturbation, as seen in the linear phase advancement (for perturbations at $\theta \geq 0.42$) recorded in the BPRCs in Figure \ref{fig:HRIe_BPRCs} (c), and (d).

\section{Conclusions and future work}
\label{sec:discussion}

The present study demonstrates that the phase response structure of bursting neuronal models is significantly different and more complicated than the usual Type I and Type II phase response found in spiking models. Fast-slow dissection, phase plane analysis, and isochron calculations help to illuminate the mechanisms underlying complex phase responses to both small and large perturbations. Our results attest to the importance of multiple time-scale dynamics in shaping the phase response features of bursting neural models. 

We are unaware of any other calculations of isochrons for bursting neuronal models, either for a full bursting system (with phase space dimension $\geq 3$) or for fast subsystem cross-sections. Our portraits of the fast subsystem phase plane and isochrons differ substantially from those of either Type I or Type II spiking neuronal models \cite{Guillamon:2009}, as well as from the isochron portrait presented in \cite{Brown:2004b, Josic:2006} for a planar reduction of the Rose-Hindmarsh model.\footnote{That model is different from the standard HR Hindmarsh-Rose model we use in this paper, but it is related in form. It is the planar reduction of a three-dimensional thalamic cell model \cite{Rose:1989}, which is itself a reduction of a model that modifies the original Hodgkin-Huxley equations to include an A-current \cite{Connor:1977}.} In these other models, the isochrons radiate nearly linearly from the central fixed point, crossing the periodic orbit almost normal to it. One consequence of this geometry is that phase response is approximately zero at spike maxima in these systems, a fact used in the derivation of analytic formulae describing phase response in for some reduced neural models in \cite{Brown:2004, Brown:2004b}. In contrast, the isochrons for the fast subsystem given by Equations (\ref{eqn:HR_full_V})--(\ref{eqn:HR_full_n}) have a complicated geometry, bending around the periodic orbit in a very different fashion; their geometry changes significantly as the homoclinic point approaches, and the magnitude of phase response is maximum near spike peaks. 

Both infinitesimal and direct BPRCs show significant phase response sensitivity: perturbation at slightly different phases may result in a large, sudden switch from phase advance to phase delay, or vice versa. This phase sensitivity is strongly associated to spike times in the burst orbit and to interactions near the homoclinic bifurcation in the fast subsystem. As might be expected, the phase response behavior of the HR model differed substantially in the weak and strong perturbation regimes, and the phase response of the model to strong perturbations cannot be inferred by simply scaling the phase response to weak perturbations.  The phase sensitivity in the strong perturbation regime is particularly dramatic, closely connected as it is with changes in the number of spikes in the perturbed burst.

Some work has been done to adapt standard coupled PRC techniques to the case of coupled bursting neurons \cite{Canavier:1999, Dror:1999, Oprisan:2005, Canavier:2005, Maran:2007}.  These efforts define a  `burst'  as a single continuous period of repetitive spiking Associated with a burst is its `phase response curve' , which determines the change in timing of the beginning of the burst in response to perturbation, and its `burst resetting curve', which specifies the change in burst duration due to perturbation. Linear stability analysis of compositions of phase response and burst resetting curves is used to predict stable phase configurations for various network architectures \cite{Oprisan:2005, Canavier:2005, Maran:2007}. Though some work considers perturbations' effects persisting over multiple burst cycles, analyses in this vein do not treat changes to the internal structure of bursts due to perturbation. Furthermore, the predictions of the coupled map analyses are supported by computer simulations of non-bursting neural models, \eg Type II Morris-Lecar neurons \cite{Oprisan:2005, Canavier:2005}. Our results suggest that care should be taken when extrapolating from the phase response properties of spiking models to those of bursting models. Even the infinitesimal BPRC for the HR model is much more complex than typical PRCs for spiking models, and the sensitive changes in spike number exhibited by the HR model indicate that developing more accurate burst resetting curves for realistic bursting models may be quite difficult. One avenue for analyzing the phase relationships in networks of bursting neurons using maps is to treat burst interaction on spike-by-spike basis and to map changes in phase, burst duration, and spike structure \cite{Sherwood:2008}.

Our analysis focused on a model saddle/homoclinic burster with just three variables. Do more realistic bursting neuronal models  (\eg biophysically derived, multiple compartments) exhibit similarly complex phase responses? What is the phase response structure and isochron geometry of other burster types with different fast subsystem bifurcation structures? We have investigated the former question and arrived at an affirmative answer; we will report our findings in another paper. The latter question remains a topic for future work. 

A limitation of the work presented here, which does not significantly affect our conclusions, is the need to compute isochrons in the fast subsystem and then invoke theorems of geometric singular perturbation theory to reason about the isochron geometry of the full system. Our methods for computing isochrons are crude, but sufficient for our purposes; much more sophisticated and accurate methods for computing planar isochrons have been developed recently \cite{Guillamon:2009}. In order to study phase response in bursting models more completely, it is desirable to be able to compute and visualize directly the isochrons of the full systems, which are manifolds of dimension two and higher. Application and further development of robust methods for multiparameter continuation and computation of higher dimensional manifolds \cite{Guckenheimer:2004, Henderson:2002, Henderson:2005} may be useful in addressing this problem.

Another limitation of the present study is its reliance on geometrical arguments and numerical calculations, rather than formal proofs. The geometric reasoning and fast-slow analysis we use are very similar to the analysis in Terman's studies of mechanisms of spike number change and transitions between tonic spiking and bursting \cite{Terman:1991, Terman:1992}, and we expect that the arguments presented here could be translated into rigorous proofs in an analogous manner. 

%This section has not offered mathematical proofs of the validity of the explanations offered. However, an approach similar to that of Terman's papers, appears to be promising, and it seems likely that the numerical and intuitive explanations offered here can be translated into rigorous proofs following that example.

%\section*{Acknowledgments}

%\bibliographystyle{plainnat}      % basic style, author-year citations
\bibliographystyle{spmpsci}      % mathematics and physical sciences
\bibliography{BPRC_Isochron_Paper}   % name your BibTeX data base 
 
 \end{document}